\documentclass[11pt]{article}
\usepackage{booktabs}
\usepackage{etoolbox}
\usepackage{tikz}
\usepackage{graphicx}
\usetikzlibrary{arrows.meta,calc}
\usetikzlibrary{calc}
\patchcmd{\thebibliography}
  {\settowidth}
  {\setlength{\itemsep}{-3pt }\settowidth}
  {}{}
\apptocmd{\thebibliography}
  {\small}
  {}{}
\usepackage[a4paper,top=3cm,bottom=2cm,left=3cm,right=3cm,marginparwidth=1.75cm]{geometry}

\usepackage{scalerel}

\usepackage{enumitem}
\usepackage[T1]{fontenc}

\usepackage{mathtools,amsthm}
\usepackage{newtxtext}
\usepackage[subscriptcorrection]{newtxmath}
\usepackage{mathrsfs}
\usepackage{graphicx,xcolor}
\usepackage{booktabs,multirow}
\usepackage{caption,subcaption}
\usepackage{tikz}
\usetikzlibrary{calc}
\usepackage{pgfplots}
\pgfplotsset{compat=1.18}
\usepackage[section]{algorithm}
\usepackage{algpseudocode}
\usepackage{ragged2e,enumitem,changepage,comment,scalerel}
\usepackage[authoryear,round]{natbib}
\usepackage[colorlinks=true,allcolors=blue]{hyperref}

\allowdisplaybreaks

\newlist{steps}{enumerate}{1}
\setlist[steps,1]{label=Step \arabic*:}

\theoremstyle{plain}
\newtheorem{theorem}{Theorem}

\newtheorem{lemma}{Lemma}[section]

\newtheorem{assumption}{Assumption}[section]

\theoremstyle{definition}

\theoremstyle{remark}
\newtheorem{remark}{Remark}[section]

\numberwithin{equation}{section}
\counterwithin{figure}{section}
\counterwithin{table}{section}

\newcommand{\R}{\mathbb{R}}

\newcommand{\sX}{\mathsf{X}}
\newcommand{\sZ}{\mathsf{Z}}

\newcommand{\E}{\mathbb{E}}

\newcommand{\OB}{\text{OBABO}}
\newcommand{\BA}{\text{BAOAB}}

\newcommand{\der}{\mathrm{d}}

\title{Stochastic counterpart of Nesterov's method for sampling 
\\and provable non-acceleration of samplers  with \\ heavy-ball structure}

\author{A. Sharma\thanks{Institut f\"{u}r Mathematik, Universit\"{a}t Potsdam; akash.sharma@uni-potsdam.de}}

\date{}

\begin{document}

\maketitle

\begin{abstract}
Polyak's heavy-ball method and Nesterov's method are two
momentum-based methods that have inspired many optimization
techniques used in applied fields, including machine learning.
Momentum also provides a connection with sampling, and interestingly, many kinetic
Langevin samplers reduce to the heavy-ball recurrence when
velocity is eliminated and noise is removed. This observation motivates a natural question: what is the corresponding sampler for Nesterov's method?  We propose a sampler whose deterministic part is exactly
Nesterov's recurrence, driven by temporally correlated Gaussian noise and requires one
new gradient evaluation per step. Noise correlation provides $\mathcal O(h^2)$ bias for smooth observables. We establish geometric ergodicity for possibly non-convex potentials. For strongly convex potentials with Lipschitz
continuous gradients, we obtain contraction rates and show that, up to an additive constant, the energy decays at a rate of order $\kappa^{-1/2}$, where $\kappa$ is the condition number. For Gaussian targets,  we show that relaxation has square-root dependence on $\kappa$. We further prove that any sampler based on
heavy-ball dynamics with additive noise cannot achieve
accelerated convergence uniformly over the class of strongly convex functions. Finally, we compare the proposed sampler with OBABO, BAOAB, HMC, MALA, BPS, and NUTS in numerical experiments and demonstrate the better performance of proposed sampler. 

\end{abstract}

\section{Introduction}

Computing expectations under Gibbs measures is central to Bayesian
inference, Monte Carlo maximum likelihood estimation
\citep{mcculloch1997maximum} and molecular dynamics \cite{leimkuhler2015molecular,lelievre2010free}.
Reliable inference, whether in statistics, physics or machine learning, require these expectations to be
computed efficiently without any bias or with controlled bias.
This is particularly difficult when the potential is ill-conditioned, that is,
when target distribution has very different scales across parameter
directions, because steps small enough to resolve narrow directions may explore
wider directions slowly. In sampling, several strategies have been proposed, in continuous as well as discrete time setting, to deal with such landscapes, including ensemble preconditioning \cite{leimkuhler2018ensemble,garbuno2020affine,ringh2025kalman}, implicit methods \cite{hodgkinson2021implicit}, stabilized explicit methods \cite{pereyra2020accelerating}, adaptive step sizes \cite{leroy2024adaptive}, and non-reversible perturbations of the dynamics \cite{lelievre2013optimal}.  

A similar situation arises in optimization, where different scales make a problem ill-conditioned and slow down gradient descent, and momentum methods are the classical way to address this. Two popular momentum methods, heavy-ball \citep{polyak1964some} and Nesterov's method \citep{nesterov1983method} both use the previous displacement, but Nesterov's method evaluates the gradient at the extrapolated point, and this changes its behaviour beyond quadratic objectives \citep{nesterov2013introductory,lessard2016analysis,goujaud2025provable}.
 Momentum also provides a natural meeting point between optimization and sampling methods. 
In sampling, discretizations of kinetic Langevin dynamics and its variants have
been widely adopted in computational statistics and machine
learning. The interesting observation is that several of these discretizations have heavy ball structure (see Section~\ref{subsec_2.2_nesterov_sharma}) in their deterministic component. One popular class of discretization of kinetic Langevin dynamics is splitting techniques
\citep{leimkuhler2016computation} (also see \citep{bou2010longaccuracy, monmarche2021high}), which concatenates flows corresponding
to linear transport, force and refreshment updates.
Treating velocity as an auxiliary variable and eliminating it
shows a structural similarity with the heavy-ball method \cite{bou2026provable}. For the two splitting methods \OB\ and \BA, the resulting deterministic position recurrences coincide, while the noise
components differ. The same heavy-ball structure appears in one-step unadjusted generalized HMC with partial refreshment
\citep{gouraud2025hmc}, the Br{\"u}nger--Brooks--Karplus (BBK)
scheme \citep{brunger1984stochastic}, the Langevin impulse
(LI) scheme \citep{skeel2002impulse}, and the
stochastic leap-frog scheme \citep{van1988leap}. 
For UBU \citep{sanz2021wasserstein} and stochastic position
Verlet (SPV) \citep{melchionna2007design}, this correspondence
holds at their force-evaluation positions.
It also holds for the inertial Langevin algorithm (ILA)
\citep{falk2025inertial}, which is built directly on the
heavy-ball recurrence.

\begin{figure}
    \centering
    \includegraphics[width=1\linewidth, height = 0.16\textheight]{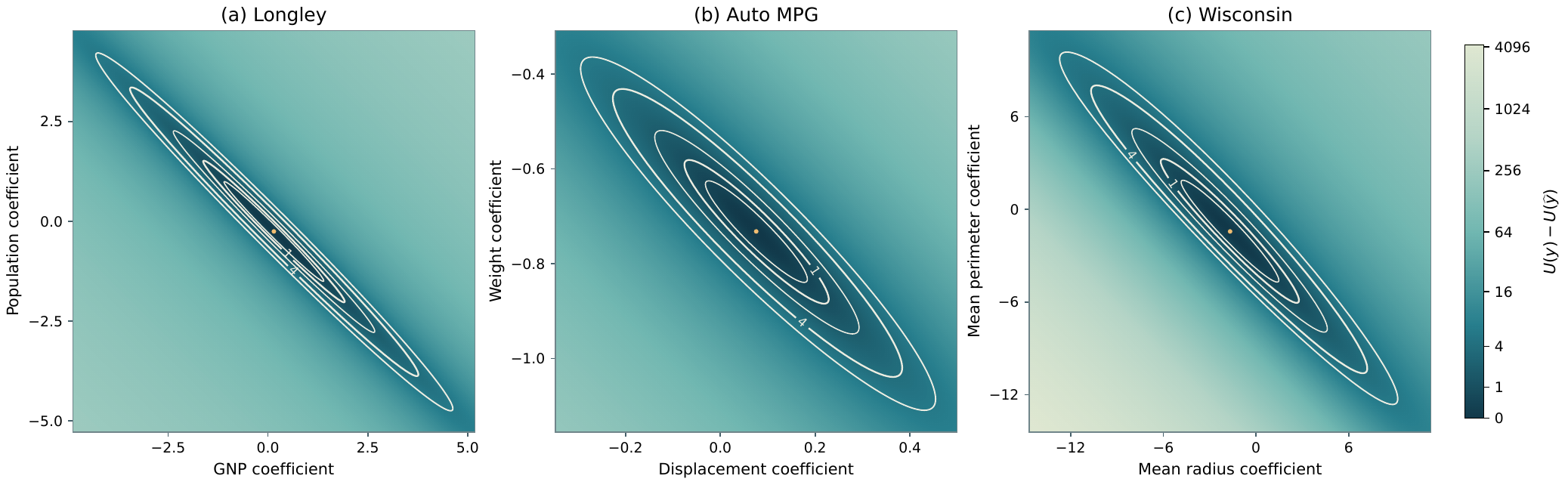}
    \caption{\small Contours of the negative log posterior density for
(a) Longley \citep{longley1967appraisal},
(b) Auto MPG \citep{auto_mpg_9} and
(c) Wisconsin breast cancer
\citep{breast_cancer_wisconsin_diagnostic_17} datasets.
Two coefficients vary in each panel, while the others are
fixed at the posterior mode. The elongated contours illustrate
ill-conditioning, where gradient descent can oscillate across
narrow directions while progressing slowly along wider directions.}
    \label{fig:placeholder}
\end{figure}

\paragraph{Contributions.} If the unadjusted overdamped Langevin algorithm is a stochastic
counterpart of gradient descent and, as mentioned above, so many kinetic Langevin
samplers (\citep{leimkuhler2016computation, bou2010longaccuracy,monmarche2021high,gouraud2025hmc,brunger1984stochastic,skeel2002impulse,van1988leap,alamo2016technique,melchionna2007design,falk2025inertial}) are stochastic counterparts of the heavy-ball method, then these connections raise a natural question: what is the stochastic
counterpart of Nesterov's position recurrences?
Adding independent Gaussian noise to the deterministic Nesterov recursion without a specific structure does
not, by itself, answer this question. The coefficients of the noise vectors, and
their relation to the deterministic components within an iteration must be taken into account as they affect the
numerical accuracy of computed expectations.

Therefore, the aim here is to construct a sampler which should have one additional gradient computation per step, Nesterov 's recurrence as its deterministic component and maintains $\mathcal{O}(h^2)$ bias in long run accuracy. For this purpose, we choose temporally correlated Gaussian
perturbations such that the leading error terms are canceled in computing expectations. We clearly see benefits of   constructing a sampler on the basis of  Nesterov recursion. On a family of nonlinear strongly convex potentials based on the constructions of \citet{goujaud2025provable} and
\citet{bou2026provable}, we compare the proposed sampler with \OB, \BA, HMC and MALA through synchronous coupling, and with BPS through errors in Gibbs expectations. The
relaxation cost of NS grows approximately as $\kappa^{1/2}$ at all 16 tested settings, compared with 3 and 4 settings for \OB\ and \BA, respectively. The coupling for HMC does not resolve in most settings and MALA shows clearly linear growth in cost with respect to $\kappa$.  For the proposed sampler, we provide the following theoretical results:
\begin{itemize}
\item  We establish geometric ergodicity of the proposed sampler under dissipativity assumption on potential which also allows for non-convex potentials (Theorem~\ref{nesterov_sharma_theorem_geometric_ergodicity}).  We establish an $\mathcal{O}(h^2)$ bias for
smooth position observables (Theorem~\ref{nesterov_sharma_theorem_2}).

\item  For strongly convex potentials with Lipschitz gradient,
we prove contraction of a synchronous coupling in a modified Euclidean norm under a
stepsize restriction (Theorem~\ref{nesterov_sharma_contraction_thrm}).  For the same class of functions, we also prove an energy decay at rate of order $\kappa^{-1/2}$ up to an additive noise term. Energy decay up to a constant is not a mixing guarantee.
Without noise, however, no such decay is possible for heavy-ball recursions at
tunings that admit a cycle  \citep{goujaud2025provable}.

\item For Gaussian targets and deterministic initialization,
we obtain a bound in KL divergence with square-root dependence on $\kappa$ in the relaxation cost
(Theorem~\ref{nesterov_acc:thm:gaussian}). With discretization error at stationarity, total variation error $\varepsilon$ is reached using $\widetilde{\mathcal O}(\sqrt{\kappa}d^{1/4}\varepsilon^{-1/2})$ gradient evaluations. This matches, up to logarithmic factors, the bound of \citet[Proposition~7]{apers2024hamiltonian}
for unadjusted HMC with Gaussian targets and random
integration times.  

\item  We perform three experiments. The first two compare relaxation cost in terms of the normalized root mean square coupling distance maximized over three deterministic starting points. In first experiment, the proposed sampler demonstrates square-root growth with $\kappa$ for all 16 settings, while \OB\ and \BA\ show linear growth for some settings and MALA does so for all settings (see  Section~\ref{nesterov_sha_subsection4.1}).
 For smoothed total variation potential function in dimension $256$, our sampler
reaches the coupling threshold for more parameter choices
and achieves lower best observed costs at larger $\kappa$
(Section~\ref{nesterov_sharma_tv_results}). 
In the third experiment, we compare seven samplers for a $1035$ dimensional Bayesian logistic regression problem on real-world mortgage application data. The proposed sampler gives the smallest coefficient error (Section~\ref{nesterov_sha_subsection4.2}).
\end{itemize}
A further  contribution is a provable non-acceleration result
for samplers whose deterministic position recurrences have
heavy-ball structure.
\begin{itemize}
\item Using the deterministic non-acceleration result of \citet{goujaud2025provable}, we rule out a uniform $\mathcal{W}_2$ convergence bound with rate $e^{-cn/\sqrt{\kappa}}$, where $c>0$ is independent of
$\kappa$, over the class of strongly convex potentials with
globally Lipschitz gradients and all deterministic initial
states
(Theorem~\ref{nesterov_sharma_thm_hb_wasserstein_nonacceleration}).
The result applies in dimension $d\geq2$ to any fixed-parameter sampler with heavy-ball deterministic structure satisfying the stated Lipschitz and scaling conditions, and is therefore not restricted to a particular sampling method.  
\end{itemize}

We discuss related literature at appropriate places in the paper highlighting and comparing the theoretical and numerical results.

\paragraph{Notation and overview.}
We use $|\cdot|$ for the Euclidean norm, $\|\cdot\|$ for
the induced matrix norm and $I_d$ for the identity matrix.
We write $F=-\nabla U$ and $\kappa=L/m$ when
$mI_d\preceq\nabla^2U\preceq LI_d$, with $0<m\leq L$.
The symbols $\mathrm{KL}$ and $\mathcal W_2$ denote relative
entropy and the $2$-Wasserstein distance, respectively.

Section~\ref{nesterov_sharma_section_2} describes the
heavy-ball structure of several kinetic Langevin samplers.
Section~\ref{nesterov_sharma_section_3} introduces the
proposed sampler and establishes its convergence properties.
Section~\ref{nesterov_sharma_sec_hb_nonacceleration} proves
the non-acceleration result for samplers with heavy-ball
structure. Section~\ref{nesterov_sharma_numerical_exp_sec_5}
presents the numerical comparisons.

\section{Momentum methods and kinetic Langevin samplers}\label{nesterov_sharma_section_2}

\subsection{Eliminating the auxiliary variable in optimization methods}

Let $F : = -\nabla U$ where $U:\mathbb{R}^d \rightarrow \R$. For a constant step size $\alpha>0$ and momentum parameter
$\eta\in[0,1)$, Nesterov's method \citep{nesterov1983method} is
\begin{align}
Y_n=X_n+\eta(X_n-X_{n-1}),\qquad
X_{n+1}=Y_n+\alpha F(Y_n)
\label{nestorovsha_nestorov_2_rec}
\end{align}
with initialization $X_{-1}=X_0$, i.e., zero initial momentum. If we eliminate $Y_n$, we get 
\begin{align}
X_{n+1}
=X_n+\eta(X_n-X_{n-1})
+\alpha F\big(X_n+\eta(X_n-X_{n-1})\big).
\label{nestorovsha_eq:nesterov-corrected}
\end{align}
Alternatively, if we eliminate $X_n$, we get, for $n\geq1$,
\begin{align}
Y_{n+1} = Y_n+\eta(Y_n-Y_{n-1})
+\alpha\big((1+\eta)F(Y_n)-\eta F(Y_{n-1})\big).
\label{nesterov_sharma_eqn_nesterov-extrapolated}
\end{align}
Here $Y_0=X_0$ and
$Y_1=Y_0+\alpha(1+\eta)F(Y_0)$. Saving the previous gradient evaluation in
\eqref{nesterov_sharma_eqn_nesterov-extrapolated} requires one new force evaluation per step.
The standard tuned parameters for strongly convex function are $\alpha=L^{-1}$ and $\eta=((\kappa^{1/2}-1)/(\kappa^{1/2}+1)$ 
\citep[Chap.~2]{nesterov2013introductory}, where $\kappa := L/m$ denotes the condition number for $m$-strongly convex function $U$ with $L$-Lipschitz gradient vector field $\nabla U$.

The heavy-ball method \citep{polyak1964some} instead uses
$D_{n+1}=\eta D_n+\alpha F(X_n)$ and $X_{n+1}=X_n+D_{n+1}$.
Eliminating the displacement $D_n=X_n-X_{n-1}$ yields
\begin{align}
X_{n+1}=X_n+\eta(X_n-X_{n-1})+\alpha F(X_n).
\label{eq:heavy-ball}
\end{align}
Therefore, the recursion \eqref{nestorovsha_eq:nesterov-corrected} changes the point of force evaluation and the recursion \eqref{nesterov_sharma_eqn_nesterov-extrapolated} adds a term of difference of gradients
correction to \eqref{eq:heavy-ball}, where as heavy ball recursion simply has one gradient term evaluated at previous iterate. The quadratic-optimal heavy-ball
parameters are
$\alpha=4/(L^{1/2}+m^{1/2})^2$ and
$\eta=\{(\kappa^{1/2}-1)/(\kappa^{1/2}+1)\}^2$.
However, this optimality on quadratics does not extend to all nonquadratic
potentials
\citep{goujaud2025provable}.

\subsection{Eliminating velocity from the kinetic samplers and heavy-ball structure}\label{subsec_2.2_nesterov_sharma}

The kinetic Langevin dynamics is
\begin{align}
\der X(t)&=V(t)\der t,\qquad
\der V(t)=F(X(t))\der t-\gamma V(t)\der t
+(2\gamma)^{1/2}\der W(t),
\label{eq:kinetic-langevin}
\end{align}
where $\gamma>0$ and $W(t)$ is standard Brownian motion in
$\R^d$. Its invariant density is proportional to
$\exp\{-U(x)-|v|^2/2\}$, with positional marginal denoted as $\mu$.
The transport, force and Ornstein-Uhlenbeck subflows are
\begin{align}
A_h:\ (x,v)&\longmapsto(x+hv,v),\qquad
B_h:\ (x,v)\longmapsto(x,v+hF(x)),
\nonumber\\
O_h:\ (x,v)&\longmapsto
\bigl(x,e^{-\gamma h}v+\sqrt{1-e^{-2\gamma h}}\,\xi\bigr),
\label{eq:kinetic-subflows}
\end{align}
where $\xi$ is a standard Gaussian vector, sampled independently
at each application of $O_h$. The schemes \OB\ and \BA\
apply these subflows in the order indicated by their names,
with step $h$ for the middle stage and $h/2$ for the others.
Eliminating velocity from consecutive steps, for $n\geq1$, we get the following recursion \cite{bou2026provable}:
\begin{align}
X_{n+1}
&=X_n+e^{-\gamma h}(X_n-X_{n-1})
+(h^2/2)(1+e^{-\gamma h})F(X_n)+\varepsilon_{n+1},
\label{nesterov_sharma_eq:splitting-position}\\
\varepsilon_{n+1}^{\mathrm{OBABO}}
&=h\sqrt{1-e^{-\gamma h}}\,
\bigl(e^{-\gamma h/2}\zeta_n+\xi_{n+1}\bigr),
\qquad
\varepsilon_{n+1}^{\mathrm{BAOAB}}
=(h/2)\sqrt{1-e^{-2\gamma h}}\,(\xi_{n+1}+\xi_n).
\label{nesterov_sharma_eq:sBAOAB_noises}
\end{align}
The vectors $\xi_n$ and $\zeta_n$ in \OB\ and $\xi_n$ in \BA\
are independent standard Gaussians. In \OB, $\xi_{n+1}$ and
$\zeta_{n+1}$ are used in the first and last refreshment of
the update from $n$ to $n+1$, respectively. Let $\hat{\mu}_h$ denote
the positional invariant measure of either splitting.
Under appropriate regularity and ergodicity assumptions,
sufficiently smooth position observables satisfy
$|\hat{\mu}_h(\varphi)-\mu(\varphi)|=\mathcal O(h^2)$
for sufficiently small $h$ \citep{leimkuhler2016computation}.

The important observation for us is that the deterministic part of
\eqref{nesterov_sharma_eq:splitting-position} is the heavy-ball
recurrence \eqref{eq:heavy-ball}, with
$\eta_{\mathrm{OBABO}}=\eta_{\mathrm{BAOAB}}=e^{-\gamma h}$
and
$\alpha_{\mathrm{OBABO}}=\alpha_{\mathrm{BAOAB}}
=(h^2/2)(1+e^{-\gamma h})$.
Therefore, \OB\ and \BA\ can be considered as two stochastic
generalizations of the same heavy-ball dynamics.
The same deterministic structure also appears in the
Br{\"u}nger--Brooks--Karplus (BBK) scheme
\cite{brunger1984stochastic}.
Expressing the one-step procedure by eliminating the velocity,
we get
\begin{align}
X_{n+1}
&=X_n+\frac{1-\gamma h/2}{1+\gamma h/2}(X_n-X_{n-1})
+\frac{h^2}{1+\gamma h/2}F(X_n)
+\frac{\sqrt{2\gamma h^3}}{1+\gamma h/2}\,\xi_{n+1},
\label{nesterov_sharma_eq:bbk-position}
\end{align}
where $(\xi_n)$ are independent standard Gaussian vectors.
Therefore, BBK also has the heavy-ball deterministic recurrence
with
$\eta_{\mathrm{BBK}}=(1-\gamma h/2)/(1+\gamma h/2)$
and $\alpha_{\mathrm{BBK}}=h^2/(1+\gamma h/2)$,
with $0<\eta_{\mathrm{BBK}}<1$ when $0<\gamma h<2$.

The Langevin impulse (LI) scheme \cite{skeel2002impulse}
and the stochastic leap frog scheme by van Gunsteren--Berendsen
\citep{van1988leap} also have this structure.
For both schemes, if we eliminate velocity, we ascertain
\begin{align}
X_{n+1}
&=X_n+\eta_s(X_n-X_{n-1})
+\alpha_s F(X_n)+\varepsilon_{n+1}^{s},
\qquad s\in\{\mathrm{LI},\mathrm{VGB}\},
\label{nesterov_sharma_eq:li-vgb-position}
\end{align}
where
\begin{align}
\eta_{\mathrm{LI}}
=\eta_{\mathrm{VGB}}=e^{-\gamma h},\qquad
\alpha_{\mathrm{LI}}
=\frac{h(1-\eta_{\mathrm{LI}})}{\gamma},
\qquad
\alpha_{\mathrm{VGB}}
=\frac{(1-\eta_{\mathrm{VGB}})^2}
{\gamma^2\sqrt{\eta_{\mathrm{VGB}}}}.
\label{nesterov_sharma_eq:li-vgb-parameters}
\end{align}
Their Gaussian noise terms can be expressed using the exact integration of transport, friction and noise. In fact, LI applies
$B_{h/2}$, $\mathcal U_h$ and $B_{h/2}$, with $\mathcal U_h$ defined below, and is therefore also called BUB \citep{sanz2021wasserstein}.
Let $\mathcal U_h$ denote the exact step of
\begin{align*}
\der X(t)&=V(t)\der t,\qquad
\der V(t)=-\gamma V(t)\der t+(2\gamma)^{1/2}\der W(t).
\end{align*}
Writing $\eta_{\mathrm{BUB}}=e^{-\gamma h}$, this step has the form
\begin{align}
\mathcal U_h:\ (x,v)\longmapsto
\left(x+\frac{1-\eta_{\mathrm{BUB}}}{\gamma}v+\zeta,\,
\eta_{\mathrm{BUB}}v+\omega\right),
\label{nesterov_sharma_eq:free-kinetic-subflow}
\end{align}
where $(\zeta,\omega)$ is centred Gaussian, with
\begin{align*}
\E[\zeta\zeta^\top]
=\left\{\frac{2h}{\gamma}
-\frac{3-4\eta_{\mathrm{BUB}}+\eta_{\mathrm{BUB}}^2}
{\gamma^2}\right\}I_d,\;
\E[\omega\omega^\top]
=(1-\eta_{\mathrm{BUB}}^2)I_d,
\;
\E[\zeta\omega^\top]
&=\frac{(1-\eta_{\mathrm{BUB}})^2}{\gamma}I_d.
\end{align*}
Independent Gaussian pairs are used in successive
$\mathcal U_h$ steps.
Using Gaussian pairs with these covariances, the noise terms
in \eqref{nesterov_sharma_eq:li-vgb-position} can be written as
\begin{align}
\varepsilon_{n+1}^{s}
&=\zeta_{n+1}-\eta_s\zeta_n
+\frac{1-\eta_s}{\gamma}\omega_n,
\qquad s\in\{\mathrm{LI},\mathrm{VGB}\},
\label{nesterov_sharma_eq:li-vgb-noise}
\end{align}
where each scheme uses its own sequence of independent
Gaussian pairs. The UBU scheme \citep{sanz2021wasserstein} applies $\mathcal U_{h/2}$,
$B_h$ and $\mathcal U_{h/2}$. Its intermediate positions $Z_n$, where
the force is evaluated, satisfy
\eqref{nesterov_sharma_eq:li-vgb-position} with $s=\mathrm{LI}$ and
$X_n$ replaced by $Z_n$, so that $\eta_{\mathrm{UBU}}=\eta_{\mathrm{LI}}$
and $\alpha_{\mathrm{UBU}}=\alpha_{\mathrm{LI}}$. Hence, the heavy-ball
correspondence holds for the full-step positions of LI and the
intermediate force-evaluation positions of UBU.

The same heavy-ball structure arises for unadjusted HMC
\citep{neal2011mcmc} when we use one leapfrog step and allow
partial momentum refreshment, as in generalized HMC
\citep{gouraud2025hmc}. In this case, before each velocity-Verlet
step, we replace $V_n$ by
$\eta V_n+\sqrt{1-\eta^2}\,\xi_{n+1}$, where
$\eta=e^{-\gamma h}$ and $(\xi_n)$ are independent standard
Gaussian vectors. Expressing the resulting updates in terms
of position alone gives
\eqref{nesterov_sharma_eq:splitting-position}, with
$\varepsilon_{n+1}=h\sqrt{1-\eta^2}\,\xi_{n+1}$.
The inertial Langevin algorithm (ILA) \citep{falk2025inertial} uses the recurrence, $X_{n+1}
=X_n+\beta(X_n-X_{n-1})+\tau F(X_n)
+\sqrt{2\tau(1-\beta)}\,\xi_{n+1}$, 
where $\tau>0$ is the step parameter and $\beta\in[0,1)$
is the momentum parameter. Taking $\beta=\eta$ and
$\tau=h^2(1+\eta)/2$ matches both the deterministic and noise
coefficients of one step unadjusted gHMC \citep{gouraud2025hmc}. Therefore, with compatible initialization,
ILA  is simply a position formulation of the one-step unadjusted
generalized HMC scheme \cite{gouraud2025hmc}.

Higher order Runge-Kutta-type schemes
\citep{milstein2003quasi,milstein2004stochastic,scott2026underdamped} and midpoint methods
\citep{shen2019randomized,altschuler2026shifted} use forces evaluated at intermediate positions and do not generally reduce to the heavy-ball or Nesterov recurrences
above, and they often require more than one new gradient evaluation per step.. The exponential integrator of \citet{cheng2018underdamped},
also investigated by \citet{dalalyan2022bounding,zhang2023improved},
gives a different deterministic position recurrence, in which the current and previous forces both have positive coefficients unlike the Nesterov sampler proposed in the next section.

\section{Stochastic counterpart of Nesterov's method}\label{nesterov_sharma_section_3}
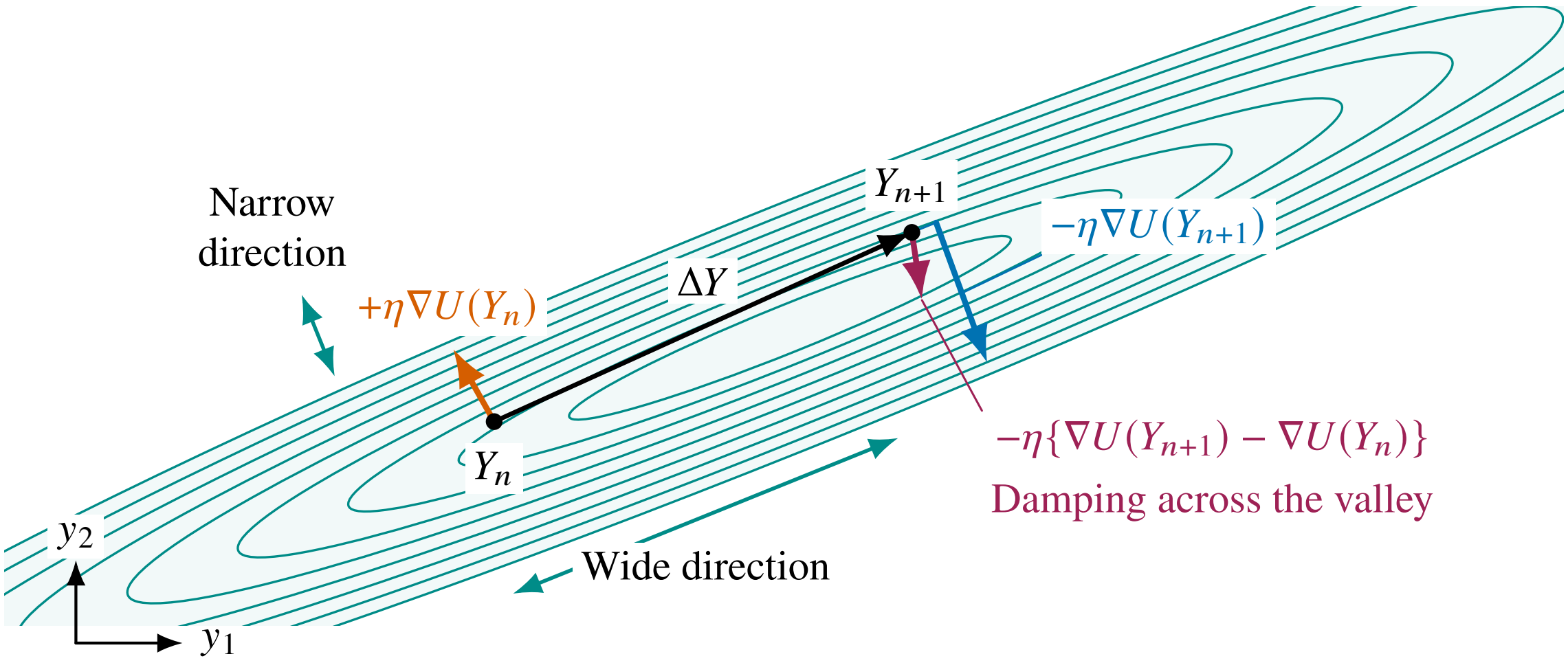
\begin{figure}[t]
\centering
\definecolor{nscontour}{HTML}{008C88}
\definecolor{nscurrent}{HTML}{0072B2}
\definecolor{nsprevious}{HTML}{D55E00}
\definecolor{nscorrection}{HTML}{9E2155}
\resizebox{\linewidth}{!}{%
\begin{tikzpicture}[
    >=Latex,
    line cap=round,
    line join=round,
    font=\small,
    force/.style={->,line width=1.4pt},
    annotation/.style={fill=white,inner sep=2pt}
]
\begin{scope}
    \clip (-6.6,-2.5) rectangle (6.5,2.7);
    \begin{scope}[rotate=22]
        \fill[nscontour!5]
            (0,0) ellipse [x radius=10,y radius=1];
        \foreach \r in {10,9,8,7,6,5,4,3,2}
            \draw[nscontour,line width=0.55pt]
                (0,0) ellipse
                [x radius=\r,y radius={\r/10}];
    \end{scope}
\end{scope}

\begin{scope}[rotate=22]
    \coordinate (Yn) at (-2.6,0.20);
    \coordinate (Ynext) at (1.25,0.36);
    \coordinate (PreviousTip) at (-2.691,0.90);
    \coordinate (CurrentStart) at (1.47,0.36);
    \coordinate (CurrentTip) at (1.42625,-0.90);

    \draw[<->,nscontour,line width=0.9pt]
        (-3,-1.20) -- (0.5,-1.20)
        node[midway,below=6pt,annotation,text=black]
        {Wide direction};
    \draw[<->,nscontour,line width=0.9pt]
        (-3.7,1.05) -- (-3.7,1.80);
\end{scope}

\coordinate (CorrectionTip) at
    ($(Ynext)+(CurrentTip)-(CurrentStart)+(PreviousTip)-(Yn)$);

\draw[->,line width=1.1pt]
    (Yn) -- (Ynext)
    node[midway,above=5pt,annotation] {$\Delta Y$};

\draw[force,nsprevious] (Yn) -- (PreviousTip);
\draw[nscurrent,line width=0.5pt]
    (Ynext) -- (CurrentStart);
\draw[force,nscurrent] (CurrentStart) -- (CurrentTip);
\draw[force,nscorrection] (Ynext) -- (CorrectionTip);

\fill (Yn) circle (2pt);
\fill (Ynext) circle (2pt);

\node[below=5pt,annotation] at (Yn) {$Y_n$};
\node[above=5pt,annotation] at (Ynext) {$Y_{n+1}$};
\node[align=center] at (-4.35,0.85)
    {Narrow\\direction};

\node[annotation,text=nsprevious,anchor=south east]
    at ($(PreviousTip)+(0.80,0.10)$)
    {$+\eta\nabla U(Y_n)$};

\node[annotation,text=nscurrent]
    (NewLabel) at (3.1,0.85)
    {$-\eta\nabla U(Y_{n+1})$};
\draw[nscurrent,line width=0.5pt]
    (NewLabel.south west) --
    ($(CurrentStart)!0.5!(CurrentTip)$);

\node[annotation,text=nscorrection,align=center]
    (CorrectionLabel) at (3.55,-1.20)
    {$-\eta\{\nabla U(Y_{n+1})-\nabla U(Y_n)\}$\\[3pt]
     Damping across the valley};
\draw[nscorrection,line width=0.5pt]
    (CorrectionTip) -- (CorrectionLabel.north west);

\draw[->,line width=0.7pt]
    (-6.0,-2.65) -- (-5.1,-2.65)
    node[right] {$y_1$};
\draw[->,line width=0.7pt]
    (-6.0,-2.65) -- (-6.0,-1.95)
    node[above,annotation] {$y_2$};
\end{tikzpicture}%
}
\caption{Illustration of the force correction
$-\eta\{\nabla U(Y_{n+1})-\nabla U(Y_n)\}$,
which opposes the displacement and gives stronger damping, where the curvature is larger.}
\label{nesterov_sharma_fig_force_correction}
\end{figure}

As is clear from our discussion in previous section that there is a remarkable  correspondence between heavy-ball dynamics and several widely known kinetic samplers. This motivates to construct a sampler whose
deterministic component has Nesterov structure. To the best of our knowledge, we are not aware of any such sampler. We below present a
position-velocity formulation sampler with correlated Gaussian
increments to cancel the leading terms in error expansion with respect to $h$ and therefore resulting in $\mathcal{O}(h^2)$ bias.
For $h>0$ and $\gamma>0$, consider
\begin{align}
Y_{n+1}
&=Y_n+hV_n+
\left((2h)/{\gamma}\right)^{1/2}
\theta_0\xi_{n+1},
\label{nesterov_sharma_position_update}\\
hV_{n+1}
&=he^{-\gamma h}V_n
+h(1-e^{-\gamma h})
\bigl\{(1+e^{-\gamma h})F_{n+1}-e^{-\gamma h}F_n\bigr\}/\gamma
\nonumber\\
&\quad +\left((2h)/\gamma\right)^{1/2}
(e^{-\gamma h}\theta_0+\theta_1)\xi_{n+1},
\label{nesterov_sharma_velocity_update}
\end{align}
where $F_n=F(Y_n)$, the vectors $(\xi_n)$ are independent
$\mathcal N(0,I_d)$, and
\begin{align}
\theta_0 =\frac{1-e^{-\gamma h}}{2}
+\left\{\frac{e^{-\gamma h}(1-e^{-2\gamma h})}{2}\right\}^{1/2},
\quad \theta_1 =\frac{1-e^{-\gamma h}}{2}
-\left\{\frac{e^{-\gamma h}(1-e^{-2\gamma h})}{2}\right\}^{1/2}.
\label{eq:first-noise-coefficients}
\end{align}
Note the opposite signs of the two force vectors in \eqref{nesterov_sharma_velocity_update}.  We can write the force term in
\eqref{nesterov_sharma_velocity_update} as
$F_{n+1}+e^{-\gamma h}(F_{n+1}-F_n)$,
which combines the force at the new position with a correction based on its change over the step.
For convex potentials, this correction acts like additional friction in the deterministic update, with stronger damping in directions of larger curvature. It can therefore reduce oscillations across narrow directions while retaining momentum along wider directions, using only
force evaluations. This interpretation is consistent with the role of the gradient correction in Nesterov's method described by \citet{shi2022understanding}.  Figure~\ref{nesterov_sharma_fig_force_correction} illustrates how the difference between successive gradients introduces extra damping.

Note that the same Gaussian vector is used in both updates \eqref{nesterov_sharma_position_update}-\eqref{nesterov_sharma_velocity_update}, and the force $F_{n+1}$ is reused in the next step. Eliminating $V_n$ variable, we get
\begin{align}
Y_{n+1}
&=Y_n+e^{-\gamma h}(Y_n-Y_{n-1})
+h(1-e^{-\gamma h})
\bigl\{(1+e^{-\gamma h})F_n-e^{-\gamma h}F_{n-1}\bigr\}/\gamma
\nonumber\\
&\quad+
\left((2h)/\gamma\right)^{1/2}
\bigl(\theta_0\xi_{n+1}+\theta_1\xi_n\bigr),
\label{nesterov_sharma_ns-position-recurrence}
\end{align}
and therefore, if the noise
$(2h/\gamma)^{1/2}
(\theta_0\xi_{n+1}+\theta_1\xi_n)$ is removed, we recover
Nesterov's recurrence
\eqref{nesterov_sharma_eqn_nesterov-extrapolated}, with
$\eta=e^{-\gamma h}$ and
$\alpha = h(1-e^{-\gamma h})/\gamma$. This is exactly the basis of derivation of sampler, i.e., we keep Nesterov's deterministic structure and choose coefficients in \eqref{eq:first-noise-coefficients} to ensure second order accuracy (see Theorem~\ref{nesterov_sharma_theorem_2}).
The sampler has the continuous time limit
\eqref{eq:kinetic-langevin} and hence the
desired invariant measure $\mu \propto e^{-U} \der x$ to approximate.  A corresponding stochastic version of \eqref{nestorovsha_eq:nesterov-corrected} is
\begin{align}
X_{n+1}
&=Y_n+
h(1-e^{-\gamma h})
F(Y_n)/\gamma
+ \left((2h)/\gamma\right)^{1/2}
\widetilde\theta_0\xi_{n+1},
\label{eq:second-correction}\\
Y_{n+1}
&=X_{n+1}+e^{-\gamma h}(X_{n+1}-X_n)
+\left((2h)/\gamma\right)^{1/2}
(\widetilde\theta_1\xi_{n+1}+\widetilde\theta_2\xi_n),
\label{eq:second-extrapolation}
\end{align}
where $\widetilde\theta_0=(1-e^{-\gamma h})^{1/2}$,
$\widetilde{\theta}_1=\theta_0-(1+e^{-\gamma h})\widetilde\theta_0$
and $\widetilde\theta_2=\theta_1+e^{-\gamma h}\widetilde\theta_0$.
Eliminating $X_n$ gives the same recurrence for $Y_n$.
With compatible initialization, the formulations therefore generate
identical positions.

We now establish the geometric ergodicity of \eqref{nesterov_sharma_position_update}. For that purpose, we impose the following assumption (see \cite{mattingly2002ergodicity}). 
\begin{assumption}
\label{nesterov_sharma_assumption_geometric_ergodicity}
 The potential $U\in C^2(\R^d)$ is
nonnegative and has a globally Lipschitz gradient.
There exist $\beta\in(0,1)$ and $b\geq0$ such that
\begin{align*}
\frac12\big(y\cdot\nabla U(y)\big)
\geq
\beta U(y)
+\frac{\gamma^2\beta(2-\beta)}{8(1-\beta)}|y|^2-b,
\qquad y\in\R^d.
\end{align*}
\end{assumption}

Let $P_h$ denote the one-step transition kernel of the chain
$(Y_n,V_n)$ defined by
\eqref{nesterov_sharma_position_update}-\eqref{nesterov_sharma_velocity_update}.

\begin{theorem}
\label{nesterov_sharma_theorem_geometric_ergodicity}
Under Assumption~\ref{nesterov_sharma_assumption_geometric_ergodicity},
there exists $h_0>0$ such that, for every $0<h\leq h_0$,
the sampler \eqref{nesterov_sharma_position_update}-\eqref{nesterov_sharma_velocity_update} has a unique invariant probability measure
$\nu_h$ on $\R^{2d}$.
Moreover, there exist constants $C_h,\lambda_h>0$ such that
\begin{align}
\left|P_h^n f(y,v)-\nu_h(f)\right|
\leq
C_h e^{-\lambda_h n}(1+|y|^2+|v|^2),
\label{nesterov_sharma_eq_geometric_ergodicity}
\end{align}
for all $n\geq0$, $(y,v)\in\R^{2d}$ and measurable functions
$f$ satisfying $|f(y,v)|\leq1+|y|^2+|v|^2$.
The constants are independent of $n$, $(y,v)$ and $f$,
but may depend on $h$.
\end{theorem}

The theorem is proved once we prove the following two lemmas \cite{mattingly2002ergodicity}.

\begin{lemma}
\label{nesterov_lem_uniform-moment-bounds_sharma}
Let Assumption~\ref{nesterov_sharma_assumption_geometric_ergodicity} be satisfied. 
Then, there exist $h_0,c>0$ such that, for every integer $p\geq1$,
there is a constant $C<\infty$ satisfying
\begin{align}
\E\bigl(|Y_n|^2+|V_n|^2\bigr)^p
&\leq
C e^{-cnh/2}
\E\bigl(1 + |Y_0|^2+|V_0|^2\bigr)^p
+C
\label{eq:moment-estimate}
\end{align}
for all $n\geq0$ and $0<h\leq h_0$, where $C$ and $c$ are independent of $h$ and $n$.
\end{lemma}

\begin{lemma}
\label{nesterov_sharma_lem_minorization}
Let $U\in C^2(\R^d)$ have a globally Lipschitz
gradient. For all sufficiently small $h>0$ and every compact
set $K\subset\R^{2d}$, there exist $\varepsilon_{h,K}>0$ and
a probability measure $\rho$ such that
\begin{align}
P_h^2((y,v),A)\geq\varepsilon_{h,K}\rho(A),
\qquad (y,v)\in K,\quad A\in\mathcal B(\R^{2d}),
\label{nesterov_sharma_eq_minorization}
\end{align}
where $P_h$ is the transition kernel  of  \eqref{nesterov_sharma_position_update}-\eqref{nesterov_sharma_velocity_update}.
\end{lemma}

There are sampling methods that preserve the target measure exactly, for example Metropolis-Hastings algorithms and their Langevin and Hamiltonian
variants \citep{hastings1970monte,roberts1996exponential,duane1987hybrid,lelievre2026unbiasing},
as well as more recent nonreversible and locally balanced constructions \citep{jansson2025rebalancing}. 
The accuracy attainable within a finite computational budget also depends on how rapidly the sampler relaxes from its initial state.
Allowing a small, controlled bias can therefore be
harmless if the sampler can deal with stiff landscapes and, therefore, reduces the
overall estimation error in finite computation budget in terms of gradient evaluations. In this scenario, one hopes for second or higher order accuracy in long time simulation. And this is our next result which establishes second order invariant accuracy and its proof
requires derivative bounds for the kinetic Langevin
semigroup, for which we need the following additional assumption.
\begin{assumption}\label{nesterov_sharma_assum_bias}
Let Assumption~\ref{nesterov_sharma_assumption_geometric_ergodicity} hold. Let $U\in C^\infty(\R^d)$ with higher derivatives having polynomial growth.
Suppose that $U=U_1+U_2$, where $U_1,U_2\in C^\infty(\R^d)$,
$U_1$ and all its derivatives are bounded, and $U_2$ is convex.
Also, assume that there exist
$c>0$ and $C\geq0$ such that
\begin{align*}
U(y)&\geq c|y|^2-C,
\qquad y\in\R^d.
\end{align*}
\end{assumption} 

Related derivative bounds were established by \citep{talay2002stochastic} under different assumptions. We use \citep[Proposition~2.15]{kopec2015weak}, since it gives explicit dependence on the test function and its
derivatives which we utilize in our proof. We do not pursue weaker assumptions here.

\begin{theorem} \label{nesterov_sharma_theorem_2}
Let Assumption~\ref{nesterov_sharma_assum_bias} hold.
Let $\mu_h$ represent $y$ marginal of an invariant
measure of \eqref{nesterov_sharma_position_update}-\eqref{nesterov_sharma_velocity_update}.
Then, for every smooth observable $\varphi$ with
polynomially growing derivatives, there exist $h_0>0$ such that
\begin{align}
|\mu_h(\varphi)-\mu(\varphi)|
\leq Ch^2,\quad 0<h<h_0,
\label{eq:first-invariant-order}
\end{align}
where $C$ is independent of $h$.
\end{theorem}

Theorem~\ref{nesterov_sharma_theorem_geometric_ergodicity} establishes geometric ergodicity under dissipativity, but does not give an explicit dependence of the convergence rate on the sampler parameters and curvature bounds. Under strong convexity and a Lipschitz condition on the
gradient,  Theorem~\ref{nesterov_sharma_contraction_thrm}
provides an explicit contraction rate and step size restriction.  To express the difference between the chains using forces evaluated at
their current positions, we need the following notation: 
\begin{align}
P_n=hV_n-\alpha(1+\eta)F(Y_n).
\label{nesterov_sharma:eq:transformed_velocity}
\end{align} 
Using
$\eta\alpha(1+\eta)-\alpha\eta=\alpha\eta^2$, we obtain
\begin{align}
Y_{n+1}
&=Y_n+P_n+\alpha(1+\eta)F(Y_n)
+((2h)/\gamma)^{1/2}\,\theta_0\xi_{n+1},
 \label{nesterov_sharmaeqn_transformed_updates_y} \\
P_{n+1}
&=\eta P_n+\alpha\eta^2F(Y_n)
+((2h)/\gamma)^{1/2}\,
(\eta\theta_0+\theta_1)\xi_{n+1}.
\label{nesterov_sharmaeqn_transformed_updates_p}
\end{align}
We consider the following modified Euclidean norm:
\begin{align}
\|(y,p)\|_h^2 =\left|y+\frac{p}{1-\eta}\right|^2
+\left|\frac{p}{\eta(1-\eta)}\right|^2,
\qquad (y,p)\in\R^{2d},
\label{nesterov_sharma_norm}
\end{align}
where $|\cdot|$ denotes the Euclidean norm. Contraction rates for kinetic Langevin schemes have been established using synchronous coupling in modified Euclidean norms \citep{monmarche2021high,sanz2021wasserstein,
leimkuhler2024contraction}. The next theorem gives a corresponding bound for our sampler \eqref{nesterov_sharma_position_update}-\eqref{nesterov_sharma_velocity_update}
in the norm \eqref{nesterov_sharma_norm}, yielding convergence
to its numerical invariant measure in the associated
Wasserstein distance. We prove it by working, at each step, in the eigenbasis of the averaged Hessian between the coupled positions. This  reduces the computations to scalar $2\times2$ recursions, so that only scalar inequalities have to be verified, instead of the block-matrix positive-definiteness criterion used by \citet{leimkuhler2024contraction}.

\begin{theorem}
\label{nesterov_sharma_contraction_thrm}
Let $U\in C^2(\R^d)$ satisfying 
$mI_d\preceq\nabla^2U(y)\preceq LI_d\,$ for every $y\in\R^d$,
where $0<m\leq L$. Let 
\begin{align}
h\leq \gamma(1-e^{-\gamma h})/(4L).
\label{nesterov_sharma:eq:step_condition}
\end{align}
Let $(Y_n,P_n)$ and $(\widetilde Y_n, \widetilde P_n)$ be two synchronously coupled chains governed by
\eqref{nesterov_sharmaeqn_transformed_updates_y}-\eqref{nesterov_sharmaeqn_transformed_updates_p}, the following contraction result hold:
\begin{align}
\|(Y_n-\widetilde Y_n,P_n-\widetilde P_n)\|_h^2
\leq
\left(1-\frac{mh}{\gamma}\right)^n
\|(Y_0-\widetilde Y_0,P_0-\widetilde P_0)\|_h^2,
\qquad n\geq0.
\label{nesterov_sharma:eq:contraction}
\end{align}
\end{theorem}

We compare this result with contraction bounds for other
kinetic Langevin samplers under strong convexity.
Since $1-e^{-\gamma h}<\gamma h$,
condition~\eqref{nesterov_sharma:eq:step_condition}
requires $\gamma>2\sqrt L$. For the same friction,
it permits larger steps than
$h\leq(1-e^{-\gamma h})/(2\sqrt L)$, required for
\BA\ and \OB\ by
\citet[Theorems~5.1--5.2]{leimkuhler2024contraction},
with strict inequality for \OB.

For $U(x)=x^\top Kx/2+G(x)$, where $K\succeq mI_d$
and $G$ is convex, the harmonic splitting of
\citet[Theorem~5]{schuh2026convergence} requires
\[
\gamma^2\geq2L_G,\qquad
h\leq\min\left\{\frac1{4\gamma},\frac{\gamma}{4L_K}\right\},
\]
where $L_K=\|K\|$ and $L_G$ is the Lipschitz constant
of $\nabla G$. The Wasserstein contraction factor is
$e^{-nh\min\{\gamma/8,m/(4\gamma)\}}$.
With $\gamma$ proportional to $\sqrt L$, these splitting
bounds and Theorem~\ref{nesterov_sharma_contraction_thrm}
permit steps of order $1/\sqrt L$ and give relaxation
costs of order $\kappa=L/m$, up to logarithmic factors.
For one-step gHMC, the sufficient conditions of
\citet[Theorem~2.11(3) and Corollary~2.12]{gouraud2025hmc}
instead give steps of order $1/(\kappa\sqrt L)$ and
relaxation costs of order $\kappa^2$ under the same friction
scaling. These comparisons concern sufficient contraction
bounds in different norms, rather than exact stability
limits or sampling costs including invariant bias.

Our next result shows the expected energy decay up to an additive constant for general strongly
convex potentials. The proof follows the energy estimates for
Nesterov's method \citep{hu2017dissipativity,bansal2019potential},
with transformation convenient for the analysis and with an additional term accounting for the Gaussian noise.
For strongly convex functions, Nesterov's parameter choice
relates strong convexity parameter $m$ with the gradient step $\alpha$ and the momentum coefficient
$\eta$ through $\eta=\frac{1-\sqrt{m\alpha}}{1+\sqrt{m\alpha}}$. 
In our sampler, these coefficients are determined by the time
step $h$ and friction $\gamma$ through
$\eta=e^{-\gamma h}$ and $\alpha=h(1-\eta)/\gamma$.
Therefore, to get the same parameter relation as in Nestorov's case, we first choose
$\alpha$ and then determine $h$ and $\gamma$ by
\begin{align}
0<\alpha\leq\frac1{4L},\quad
\tau=\sqrt{m\alpha},\quad
\eta=\frac{1-\tau}{1+\tau},\quad
h=\left(\frac{-\alpha\log\eta}{1-\eta}\right)^{1/2},
\quad
\gamma=\frac{-\log\eta}{h}.
\label{nesterov_sharma:eq:energy-parameters}
\end{align}
The time step can be made arbitrarily small by reducing
$\alpha$, with $\gamma$ adjusted accordingly.  These choices satisfy $\eta=e^{-\gamma h}$ and
$\alpha=h(1-\eta)/\gamma$.
Recall that $P_n=hV_n-\alpha(1+\eta)F(Y_n)$, and we need the following notation for analysis below:
\begin{align}
\mathsf X_n=Y_n-\frac{P_n}{\eta},
\qquad
\sZ_n=Y_n+\frac{P_n}{\eta\tau}.
\label{nesterov_sharma:eq:energy-coordinates}
\end{align}
We have chosen the coordinates in \eqref{nesterov_sharma:eq:energy-coordinates}
 so that the decrease in the potential energy cancels
the squared-force term arising from the update of $\mathsf Z_n$ as can be seen in the proof of Theorem~\ref{nesterov_sharma_thrm_energy_decay}.
As can be easily verified $Y_n=(\mathsf{X}_n+\tau\mathsf Z_n)/(1+\tau)$.
We denote the energy function as
\begin{align}
\mathcal{E}_n
=U(\sX_n)-U(x_\star)
+\frac m2|\mathsf Z_n-x_\star|^2,
\label{nesterov_sharma_eqn_nesterov-energy}
\end{align}
where $x_\star$ is the minimizer of $U$.

\begin{theorem}
\label{nesterov_sharma_thrm_energy_decay}
Let $U\in C^2(\R^d)$ satisfying 
$mI_d\preceq\nabla^2U(y)\preceq LI_d$ for every $y\in\R^d$,
where $0<m\leq L$.
With the parameters
\eqref{nesterov_sharma:eq:energy-parameters} and a deterministic
initial state, we have
\begin{align}
\mathbb E\mathcal E_n
\leq
(1-\tau)^n\mathcal E_0
+\frac{K_0}{\tau}\{1-(1-\tau)^n\},
\qquad n\geq0,
\label{nesterov_sharma_eq:strongly-convex-energy_resu}
\end{align}
where $K_0
=\frac{dh}{\gamma}
\left[\frac{L\theta_1^2}{\eta^2}
+m\left(
\theta_0+\frac{\eta\theta_0+\theta_1}{\eta\tau}
\right)^2
\right]$. 
In particular, taking $\alpha=1/(4L)$ gives
\begin{align}
\mathbb E\mathcal E_n
\leq
e^{-n/(2\sqrt\kappa)}\mathcal E_0
+21d\sqrt\kappa,
\qquad \kappa=L/m.
\label{nesterov_sharma:eq:accelerated-strongly-convex-energy}
\end{align}
\end{theorem}

For $\alpha=1/(4L)$ and $0<\varepsilon<1$, the contribution
of the initial energy is at most $\varepsilon\mathcal E_0$
whenever $n\geq2\sqrt\kappa\log(1/\varepsilon)$.
Since each iteration requires one new gradient evaluation,
this gives
\[
\mathcal O\left(\sqrt\kappa\log\frac{e}{\varepsilon}\right)
\]
gradient evaluations, including the initial gradient.
The accumulated contribution of the Gaussian noise is bounded
by $21d\sqrt\kappa$, independently of the initialization and $n$.

For general storngly convex function, Theorem~\ref{nesterov_sharma_thrm_energy_decay} keeps
the square-root dependence on $\kappa$ in the energy
estimates for Nesterov's method
\citep{hu2017dissipativity,bansal2019potential},
with an additive term accounting for the Gaussian noise. Its parameter choice uses friction of order $\sqrt m$, whereas Theorem~\ref{nesterov_sharma_contraction_thrm}
requires $\gamma>2\sqrt L$.
The energy estimate controls a single trajectory relative to the minimizer and does not by itself lead to accelerated
mixing. However, accelerated mixing can be provided for quadratic potentials which we consider next. Let 
\[
U(y)=U_\star+\frac12(y-x_\star)^\top A(y-x_\star),
\qquad mI_d\preceq A\preceq LI_d,\qquad 0<m\leq L.
\]
Therefore, the target measure is $\mu=\mathcal N(x_\star,A^{-1})$. Recall $\eta=e^{-\gamma h}$, $\alpha=h(1-\eta)/\gamma$ and let
\begin{align}
\vartheta : = \max\left\{
\frac{\gamma}{2\sqrt m}+h\sqrt m,\,
h\sqrt L+\frac{\gamma}{2\sqrt L}
\right\}.
\label{nsacc:vartheta}
\end{align}
Let $a_i$ and $u_i$ be the eigenvalues and orthonormal
eigenvectors of $A$. We also need the following notation
\[
E_0=|A^{1/2}(y_0-x_\star)|^2+|v_0|^2.
\]
Let $\epsilon_i=a_i\sigma_{\infty,i}^2-1$, where
$\sigma_{\infty,i}^2$ is the stationary variance of
$(u_i\cdot(Y_n-x_\star))$.

\begin{theorem}
\label{nesterov_acc:thm:gaussian}
Assume that $\vartheta<1$, $0<h\leq1/(2\sqrt L)$ and
$(Y_0,V_0)=(y_0,v_0)$ is deterministic. Then
\begin{align}
\mathrm{KL}\bigl(\operatorname{Law}(Y_n)\mid\mu\bigr)
\leq\sum_{i=1}^d\epsilon_i^2
+Cd\,e^{-2\gamma hn}+CE_0e^{-\gamma hn}, \qquad n\geq\frac{\log(2C)}{\gamma h}.
\label{nsacc:gauss-relax}
\end{align}
In particular, taking $h=1/(2\sqrt L)$ and
$\gamma=\sqrt m/2$, for which $\vartheta\leq3/4$, requires
\[
\mathcal O\!\left(
\sqrt\kappa\log\frac{2(d+E_0)}{\delta}
\right)
\]
gradient evaluations to obtain
$\mathrm{KL}(\operatorname{Law}(Y_n)\mid\mu)
\leq\sum_{i=1}^d\epsilon_i^2+\delta$, for $0<\delta<1$.
\end{theorem}

The exponential terms decrease with the number of iterations,
while $\sum_i\epsilon_i^2$ accounts for the discretization error at stationarity of chain. 
We can easily obtain a bound on discretization error, i.e., $\sum_i\epsilon_i^2\leq h^4\sum_i a_i^2$. For Gaussian targets, Theorem~\ref{nesterov_acc:thm:gaussian}
gives a finite-iteration bound in relative entropy which imply results in 2-Wasserstein and total variation.
The contraction results of \citet{leimkuhler2024contraction}
are in Wasserstein distance, and their Gaussian example for BAO scheme also exhibits square-root dependence on $\kappa$ in the asymptotic contraction rate.
Gaussian acceleration in $\mathcal W_2$ is also established for their scheme by \citet[Section~3.1]{falk2025inertial}.
 For HMC, \citet{apers2024hamiltonian} obtain bounds in total variation distance, which in our case is implied by Pinsker's inequality. Choosing $h$ to control the stationary error gives total variation error at most
$\varepsilon$ using
\[
\mathcal O\left(
\frac{\sqrt\kappa\,d^{1/4}}{\sqrt\varepsilon}
\log\frac{2(d+E_0)}{\varepsilon}
\right)
\]
gradient evaluations.
Therefore, we control relative entropy while matching, up to
logarithmic factors, the total variation gradient cost
of unadjusted HMC
\citep[Proposition~7]{apers2024hamiltonian}.

\section{Non-acceleration of samplers with heavy-ball structure}
\label{nesterov_sharma_sec_hb_nonacceleration}
In Section~\ref{nesterov_sharma_section_2}, we identified the
heavy-ball structure underlying several kinetic Langevin
samplers
\citep{leimkuhler2016computation,bou2010longaccuracy,
monmarche2021high,gouraud2025hmc,brunger1984stochastic,
skeel2002impulse,van1988leap,alamo2016technique,
melchionna2007design,falk2025inertial}.
This correspondence is also studied for the $A$, $B$ and $O$
splitting schemes in \citet{bou2026provable}.
To examine the consequences of heavy-ball structure, we consider
samplers of the form (see Section~\ref{nesterov_sharma_section_2})
\begin{align}
X_{n+1}
&=X_n+\eta(X_n-X_{n-1})
-\alpha\nabla U(X_n)+\varepsilon_{n+1},
\qquad \alpha>0,\quad 0\leq\eta<1,
\label{nesterov_sharma_eq_hb_noisy_recurrence}
\end{align}
where $(\varepsilon_n)$ are fixed linear combinations of independent, identically distributed centred random vectors with finite second moments. Therefore, the increments $\varepsilon_n$ may be correlated across steps and need not
be Gaussian. Throughout this section, we fix $d\geq2$ and consider
$U\in C^\infty(\R^d)$ satisfying
\begin{align}
mI_d\preceq\nabla^2U(x)\preceq LI_d,
\qquad x\in\R^d,\qquad 0<m\leq L,
\label{nesterov_sharma_eq_hb_curvature_bounds}
\end{align}
and write $\kappa=L/m$. By translating the minimizer and
subtracting the minimum value, we take $x_\star=0$ and
$U(0)=0$, without loss of generality. The sampler parameters
and noise distribution may depend on $m,L,d$, but not on
the particular potential or initialization.

Our aim is to determine whether these samplers can satisfy
a uniform accelerated convergence bound in $\mathcal W_2$,
the $2$-Wasserstein distance.
For Strang splitting schemes, 
\citet{bou2026provable} establishes non-acceleration through
quantitative mixing time lower bounds in total variation.
In comparison, we use a contrapositive argument. A uniform accelerated
$\mathcal W_2$ bound would imply, after dilation, accelerated
convergence of the deterministic recurrence but this result would contradict
the heavy-ball obstruction proved in \citet{goujaud2025provable}.
This gives a shorter proof of uniform non-acceleration, without requiring separate estimates of how long noisy trajectories remain near deterministic cycles. Moreover, under the stated assumptions, the argument applies beyond
Strang splittings to other kinetic samplers with heavy ball structure and also allows correlated, non-Gaussian noise with finite second moments.

To include possible dependence between successive increments, let $Z_n$ contain the position and the additional variables
needed for the next update and write
\begin{align}
Z_{n+1}=\Phi_U(Z_n,\xi_{n+1}),
\label{nesterov_sharma_eq_hb_markov_update}
\end{align}
where $(\xi_n)$ are independent, identically distributed
centred random vectors with finite second moments,
independent of the initial state.
We assume that $\Phi_U(0,0)=0$ and
\begin{align}
|\Phi_U(z,\xi)-\Phi_U(\widetilde z,\widetilde\xi)|
\leq a|z-\widetilde z|+b|\xi-\widetilde\xi|,
\label{nesterov_sharma_eq_hb_update_lipschitz}
\end{align}
with finite constants $a,b$ independent of $U$.
In particular, $\Phi_U$ is continuous, and
\eqref{nesterov_sharma_eq_hb_markov_update} defines a
time homogeneous Markov chain.

To compare the stochastic and deterministic recurrences,
we use the dilation
$U_R(x)=R^2U(x/R)$, $R\geq1$,
from \citet[Lemma~5.4]{goujaud2025provable}.
It preserves the curvature bounds in
\eqref{nesterov_sharma_eq_hb_curvature_bounds}.
We require the update \eqref{nesterov_sharma_eq_hb_markov_update} to satisfy 
\begin{align}
\frac1R\Phi_{U_R}(Rz,\xi)=\Phi_U(z,\xi/R).
\label{nesterov_sharma_eq_hb_update_scaling}
\end{align}
Therefore, dividing the state by $R$ gives the update for $U$
with the random input divided by $R$.
When the random inputs are set to zero, the positions
must follow the deterministic heavy-ball recurrence, and
every initial pair $(x_{-1},x_0)$ of that recurrence must
be represented by a suitable initial state $z$.

Every sampler in Section~\ref{subsec_2.2_nesterov_sharma} admits a representation \eqref{nesterov_sharma_eq_hb_markov_update}: OBABO, BAOAB, LI and one-step gHMC with $Z_n=(X_n,V_n)$, BBK with $Z_n=(X_n,V_{n-1/2})$, UBU with its position and velocity at the force-evaluation points, ILA with $Z_n=(X_n,X_{n-1})$, and VGB with $Z_n=(X_n,X_{n-1},\chi_n)$, where $\chi_n$ is the vector computed using all Gaussian vectors drawn in step $n$. In each case, $\Phi_U$ is a composition of linear maps, evaluations of $\nabla U$ and additive noise, therefore $\Phi_U(0,0)=0$, \eqref{nesterov_sharma_eq_hb_update_lipschitz} and \eqref{nesterov_sharma_eq_hb_update_scaling} follow from $\nabla U(0)=0$, the $L$-Lipschitz continuity of $\nabla U$ and $\nabla U_R(Rx)=R\nabla U(x)$.

Let $\Pi_U$ denote an invariant distribution of the full
state $Z_n$ with finite second moment. The uniform convergence guarantee considered here requires such a distribution for every
potential satisfying 
\eqref{nesterov_sharma_eq_hb_curvature_bounds}.

\begin{theorem}
\label{nesterov_sharma_thm_hb_wasserstein_nonacceleration}
Let $d\geq2$. Consider a sampler whose positions satisfy \eqref{nesterov_sharma_eq_hb_noisy_recurrence}, 
where the centred increments $\varepsilon_{n+1}$ are fixed
linear combinations of independent, identically distributed
centred random noise vectors with finite second moments, independent
of the initial state. Let its Markov representation
\eqref{nesterov_sharma_eq_hb_markov_update} satisfy
$\Phi_U(0,0)=0$,
\eqref{nesterov_sharma_eq_hb_update_lipschitz} and
\eqref{nesterov_sharma_eq_hb_update_scaling}.
Then no choice of fixed parameters and noise distribution based only on $m,L,d$ provides invariant distributions $\Pi_U$ with finite second moments and the bound 
\begin{align}
\mathcal W_2\bigl(\operatorname{Law}_z(Z_n),\Pi_U\bigr)
\leq
Ce^{-cn/\sqrt\kappa}\mathcal{W}_2(\delta_z,\Pi_U)
\label{nesterov_sharma_eq_hb_assumed_wasserstein}
\end{align}
for every normalized potential $U\in C^\infty(\R^d)$
satisfying \eqref{nesterov_sharma_eq_hb_curvature_bounds},
every deterministic initial state $z$ and every $n\geq0$,
for all sufficiently large $\kappa$.
Here, $c>0$ is independent of $\kappa$, and $C<\infty$
may depend on $m,L,d$ and the chosen parameters and noise
distribution, but not on $U,z,n$.
\end{theorem}

\begin{proof}
Let \eqref{nesterov_sharma_eq_hb_assumed_wasserstein} hold.
Since the minimizer is at the origin, we have $
|\nabla U(x)|\leq L|x|$ and $
\frac m2|x|^2\leq U(x)\leq\frac{L}{2}|x|^2$. We first bound the second moment of the stationary measure using the assumed
convergence bound. Starting from $Z_0=0$, using
\eqref{nesterov_sharma_eq_hb_update_lipschitz} and
$\Phi_U(0,0)=0$, we get
\begin{align*}
\bigl(\mathbb{E}|Z_{n+1}|^2\bigr)^{1/2}
\leq a\bigl(\E|Z_n|^2\bigr)^{1/2}
+b\bigl(\E|\xi_1|^2\bigr)^{1/2}.
\end{align*}
Therefore, we have
\begin{align}
\bigl(\E|Z_N|^2\bigr)^{1/2}
\leq
b\bigl(\E|\xi_1|^2\bigr)^{1/2}
\sum_{j=0}^{N-1}a^j
=:A_N.
\label{nesterov_sharma_eq_hb_finite_time_moment}
\end{align}
For each finite $N$, this bound is uniform over the potentials
under consideration and it does not require $a<1$. Let $
M_U : =\left(\int|z|^2\,\Pi_U(\der z)\right)^{1/2}$. 
We choose $N\geq1$ such that $Ce^{-cN/\sqrt\kappa}\leq1/2$.
Using triangle inequality and then using 
\eqref{nesterov_sharma_eq_hb_finite_time_moment} and
\eqref{nesterov_sharma_eq_hb_assumed_wasserstein}, we obtain
\begin{align*}
M_U &= \mathcal W_2(\delta_0,\Pi_U) 
\leq \mathcal W_2\bigl(\delta_0,\operatorname{Law}_0(Z_N)\bigr)
+
\mathcal{W}_2\bigl(\operatorname{Law}_0(Z_N),\Pi_U\bigr) \leq A_N+\frac12M_U.
\end{align*}
Therefore,
\begin{align}
M_U\leq2A_N.
\label{nesterov_sharma_eq_hb_stationary_moment}
\end{align}
We now use this bound to estimate the expected potential energy. Using Lipschitz continuity of gradient vector field, we get
\begin{align}
\E U(X_n)
\leq L \mathcal W_2^2\bigl(\operatorname{Law}_z(Z_n),\Pi_U\bigr)
+LM_U^2.
\label{nesterov_sharma_eq_hb_energy_from_wasserstein}
\end{align}
Also, $W_2^2(\delta_z,\Pi_U)
=\int|z-w|^2\,\Pi_U(\der w)
\leq2|z|^2+2M_U^2$.  
Combining these estimates with
\eqref{nesterov_sharma_eq_hb_assumed_wasserstein} and
\eqref{nesterov_sharma_eq_hb_stationary_moment}, we obtain
\begin{align}
\E_z U(X_n)
\leq
2LC^2e^{-2cn/\sqrt\kappa}|z|^2
+4L(2C^2+1)A_N^2.
\label{nesterov_sharma_eq_hb_expected_energy}
\end{align}
The last term is independent of the potential and
initialization. Consider the function 
$U_R(x):=R^2U(x/R)$, with $R\geq1$. The curvature  bounds of $U_R$ are
unchanged because
\begin{align*}
\nabla U_R(x)=R\nabla U(x/R), \qquad
\nabla^2U_R(x)=\nabla^2U(x/R).
\end{align*}
Hence the parameters and all constants in
\eqref{nesterov_sharma_eq_hb_expected_energy} remain fixed as $R$ increases. Let $(Z_n^R)_{n\geq0}$ denote the sampler with potential $U_R$,
starting from $Z_0^R=Rz$. We divide the state by $R$ to compare this sampler with the deterministic update for $U$. Using \eqref{nesterov_sharma_eq_hb_update_scaling}, we get
\begin{align*}
\frac{Z_{n+1}^R}{R}
=\Phi_U\left(\frac{Z_n^R}{R},\frac{\xi_{n+1}}{R}\right).
\end{align*}
Therefore, after rescaling, the potential is $U$ and the random
input is divided by $R$. Let $(z_n)_{n\geq0}$ denote the
corresponding deterministic sequence, i.e.,
\begin{align*}
z_{n+1}=\Phi_U(z_n,0),\qquad z_0=z.
\end{align*}
Since $Z_0^R/R=z_0$, using
\eqref{nesterov_sharma_eq_hb_update_lipschitz}, we obtain
\begin{align*}
\left(\E\left|\frac{Z_{n+1}^R}{R}-z_{n+1}\right|^2\right)^{1/2}
\leq
a\left(\E\left|\frac{Z_n^R}{R}-z_n\right|^2\right)^{1/2}
+\frac bR\bigl(\E|\xi_1|^2\bigr)^{1/2}.
\end{align*}
Applying this inequality repeatedly, we get
\begin{align}
\left(\E\left|\frac{Z_n^R}{R}-z_n\right|^2\right)^{1/2}
\leq
\frac{b(\E|\xi_1|^2)^{1/2}}{R}
\sum_{j=0}^{n-1}a^j.
\label{nesterov_sharma_eq_hb_deterministic_limit}
\end{align}
Therefore, for each fixed $n$, $Z_n^R/R$ converges to $z_n$
in mean square as $R\to\infty$. We use this convergence to pass to the limit in the
expected potential energy. Let $x_n$ denote the position
component of $z_n$. Since $\nabla U(0)=0$ and the gradient
is globally Lipschitz, we have
\begin{align*}
|U(x)-U(y)|\leq L(|x|+|y|)|x-y|.
\end{align*}
Consequently, by the Cauchy-Schwarz-Bunyakovsky inequality, we have
\begin{align*}
\left|\E U(X_n^R/R)-U(x_n)\right|
&\leq
L\left\{\E\left(|X_n^R/R|+|x_n|\right)^2\right\}^{1/2}
\left\{\E|X_n^R/R-x_n|^2\right\}^{1/2}.
\end{align*}
For fixed $n$, the first factor is bounded uniformly over
$R\geq1$ by
\eqref{nesterov_sharma_eq_hb_deterministic_limit}, while
the second factor tends to zero as $R \rightarrow \infty$. Hence $
\E U(X_n^R/R)\longrightarrow U(x_n)$. Applying \eqref{nesterov_sharma_eq_hb_expected_energy}
to the sampler with potential $U_R$ and initial state $Rz$,
and using $U_R(X_n^R)/R^2=U(X_n^R/R)$, we obtain
\begin{align}\label{nesterov_sharma_eqn_4.12}
\E U(X_n^R/R)
\leq
2LC^2e^{-2cn/\sqrt\kappa}|z|^2
+\frac{4L(2C^2+1)A_N^2}{R^2}.
\end{align}
The constants are independent of $R$. Therefore, as $R\to\infty$ we have $U(x_n)
\leq
2LC^2e^{-2cn/\sqrt\kappa}|z|^2$, $ n\geq0$. 
Since $U(x_n)\geq m|x_n|^2/2$, we get
\begin{align}
|x_n|
\leq
2C\sqrt\kappa\,e^{-cn/\sqrt\kappa}|z|,
\qquad n\geq0.
\label{nesterov_sharma_eq_hb_deterministic_distance_eqn_4.13}
\end{align}
Therefore, the assumed accelerated convergence for the sampler implies
an accelerated convergence bound for its deterministic
position sequence. When the noise vectors are zero, i.e. no random input, this sequence satisfies
\begin{align}
x_{n+1} =x_n+\eta(x_n-x_{n-1})-\alpha\nabla U(x_n).
\label{nesterov_sharma_eq_hb_deterministic_recurrence}
\end{align}
By the assumption on initialization, every initial pair $(x_{-1},x_0)$ can be represented by an initial state $z$. Therefore, \eqref{nesterov_sharma_eq_hb_deterministic_distance_eqn_4.13}
holds for every initial pair of this heavy-ball recurrence.  We now apply the deterministic non-acceleration result of \citep{goujaud2025provable}.
For every fixed choice of $\alpha,\eta$, either a quadratic potential admits a trajectory with asymptotic distance contraction factor at least $1-C_0/\kappa$, where $C_0$ is an absolute constant (i.e. independent of the potential function, dimension, \(m,L,\kappa\), and the parameters \(\alpha,\eta\)) or a smooth strongly convex potential admits a nonconvergent heavy-ball trajectory.

\textbf{Case I.} First, suppose that a deterministic trajectory does not
converge to $0$. Then, there exists $r>0$ such that $|x_n|\geq r$ for arbitrarily large $n$. Since the right-hand
side of \eqref{nesterov_sharma_eq_hb_deterministic_distance_eqn_4.13}
tends to zero, we can choose one such sufficiently large $n$ so
that \begin{align*}
2C\sqrt\kappa\,|z|e^{-cn/\sqrt\kappa}<r\leq|x_n|.
\end{align*}
This contradicts \eqref{nesterov_sharma_eq_hb_deterministic_distance_eqn_4.13}, and therefore, the assumed sampling bound cannot hold when the deterministic recurrence admits such a trajectory.

\textbf{Case II.} In the second case, the obstruction is a quadratic trajectory
whose asymptotic rate satisfies $
\limsup_{n\to\infty}|x_n|^{1/n}
\geq1-\frac{C_0}{\kappa}$. 
Here convergence to the origin is not itself a contradiction.
We need to show that the trajectory cannot converge at the
rate required by
\eqref{nesterov_sharma_eq_hb_deterministic_distance_eqn_4.13}.
For fixed $\kappa$ and a nonzero initial state $z$, taking
$n$th roots in that estimate gives
\begin{align*}
|x_n|^{1/n}
\leq
\bigl(2C\sqrt\kappa\,|z|\bigr)^{1/n}
e^{-c/\sqrt\kappa}.
\end{align*}
Since $\bigl(2C\sqrt\kappa\,|z|\bigr)^{1/n} =
\exp\left\{\frac{\log(2C\sqrt\kappa\,|z|)}{n}\right\}
\longrightarrow 1$,
we obtain
\begin{align*}
\limsup_{n\to\infty}|x_n|^{1/n}
\leq e^{-c/\sqrt\kappa}.
\end{align*}
Therefore, even though the prefactor may depend on $\kappa$,
it does not change the asymptotic rate for fixed $\kappa$.
However, since $c>0$ and $C_0$ are independent of $\kappa$,
\begin{align*}
1-\frac{C_0}{\kappa}>e^{-c/\sqrt\kappa}
\end{align*}
for sufficiently large $\kappa$. The deterministic lower
bound and the rate implied by the sampling estimate are incompatible. As a consequence, this trajectory also
violates \eqref{nesterov_sharma_eq_hb_deterministic_distance_eqn_4.13}
at some finite $n$.

In both cases, the violation occurs at a finite time and 
after choosing this time, we can take $R$ sufficiently
large that $\E U(X_n^R/R)$ is arbitrarily close to $U(x_n)$ in \eqref{nesterov_sharma_eqn_4.12} and the additive term
$4L(2C^2+1)A_N^2/R^2$ in \eqref{nesterov_sharma_eqn_4.12} is arbitrarily small.
Hence the strict violation of the deterministic energy
bound also contradicts the rescaled expected energy bound.
Here we have used convergence as $R\to\infty$ only for the chosen finite $n$, and does not require the noisy trajectory to remain close to the deterministic trajectory for all time.
Therefore the assumed uniform accelerated sampling bound   cannot hold.
\end{proof}

The same argument gives a positional version, without
requiring convergence of the auxiliary variables.

\section{Numerical experiments}\label{nesterov_sharma_numerical_exp_sec_5}
\subsection{Verification of convergence result of Theorem~\ref{nesterov_sharma_theorem_2}.}
\label{sec:archive_order_check}
\begin{figure}[htbp!]
    \centering
    \includegraphics[width=0.55\linewidth, height = 0.2\textheight]{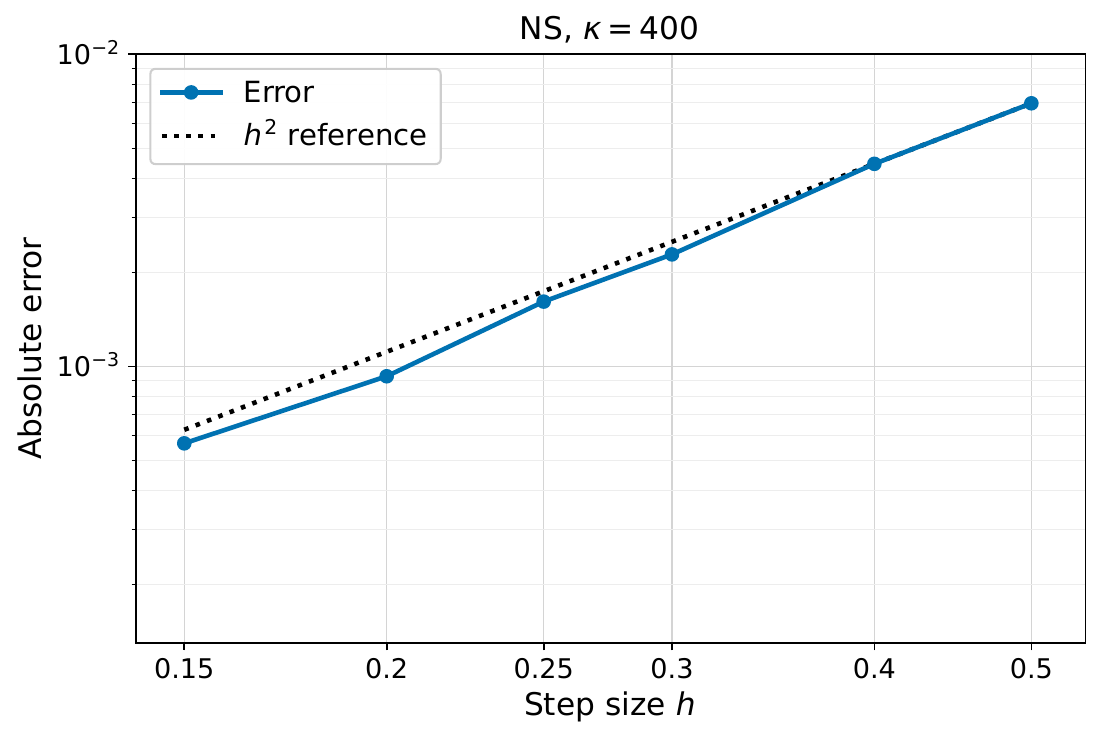}
    \caption{Plot confirming $\mathcal{O}(h^2)$ order of convergence (i.e., the result of Theorem~\ref{nesterov_sharma_theorem_2}) of the proposed sampling scheme \eqref{nesterov_sharma_position_update}-\eqref{nesterov_sharma_velocity_update} approximating kinetic Langevin dynamics in long time simulation. Table~\ref{nesterov_sharma_order_results} provides additional details.}
    \label{second_order_convergence_result}
\end{figure}

Consider the potential function
\begin{align}
U(x)&=\frac m2x^2+\frac{L-m}{a^2}\log\cosh(ax),
\label{eq:archive_smooth_potential}
\end{align}
with $m=0.0025$, $ L=1$, $a=1$ and $\gamma=1$. As can be verified, $m\leq U''(x)\leq L$ and $\kappa=400$. The observable is
$\varphi(x)=\operatorname{sech}^2(ax)$, with Gibbs expectation
$0.500621240859$ (computed using quadrature rule with error tolerance
$7.56\times10^{-13}$).
Every trajectory starts from $Y_{-1}=Y_0=0$, and the term $\theta_1\xi_0$ in the first step uses a standard Gaussian $\xi_0$ drawn independently of $(\xi_n)_{n\geq1}$. For each step size, we simulated $2^{25}=33554432$ independent trajectories. The estimates are
endpoint ensemble averages, without time averaging. The
final time is 100. The results are reported in Table~\ref{nesterov_sharma_order_results} and plotted in Figure~\ref{second_order_convergence_result} which verify the second order convergence in terms of discretization $h$.  
\begin{table}[t]
\centering
\caption{Ensemble averaging estimate for proposed sampler \eqref{nesterov_sharma_position_update}-\eqref{nesterov_sharma_velocity_update} The standard error is the Monte Carlo standard error
of the ensemble mean.}
\label{nesterov_sharma_order_results}
\begin{tabular}{r r r r r}
\toprule
$h$ & Steps & Estimate & Absolute error & MC standard error\\
\midrule
0.50 & 200 & 0.4936580633 & $6.9632\times10^{-3}$ & $6.081\times10^{-5}$\\
0.40 & 250 & 0.4961655227 & $4.4557\times10^{-3}$ & $6.086\times10^{-5}$\\
0.30 & 334 & 0.4983366632 & $2.2846\times10^{-3}$ & $6.092\times10^{-5}$\\
0.25 & 400 & 0.4990082324 & $1.6130\times10^{-3}$ & $6.095\times10^{-5}$\\
0.20 & 500 & 0.4996913044 & $9.2994\times10^{-4}$ & $6.096\times10^{-5}$\\
0.15 & 668 & 0.5000542106 & $5.6703\times10^{-4}$ & $6.097\times10^{-5}$\\
\bottomrule
\end{tabular}
\end{table}

\subsection{Comparison using strongly log-concave measure}\label{nesterov_sha_subsection4.1}
\begin{figure}[htbp!]
    \centering
    \includegraphics[width=0.75\linewidth,height = 0.3\textheight]{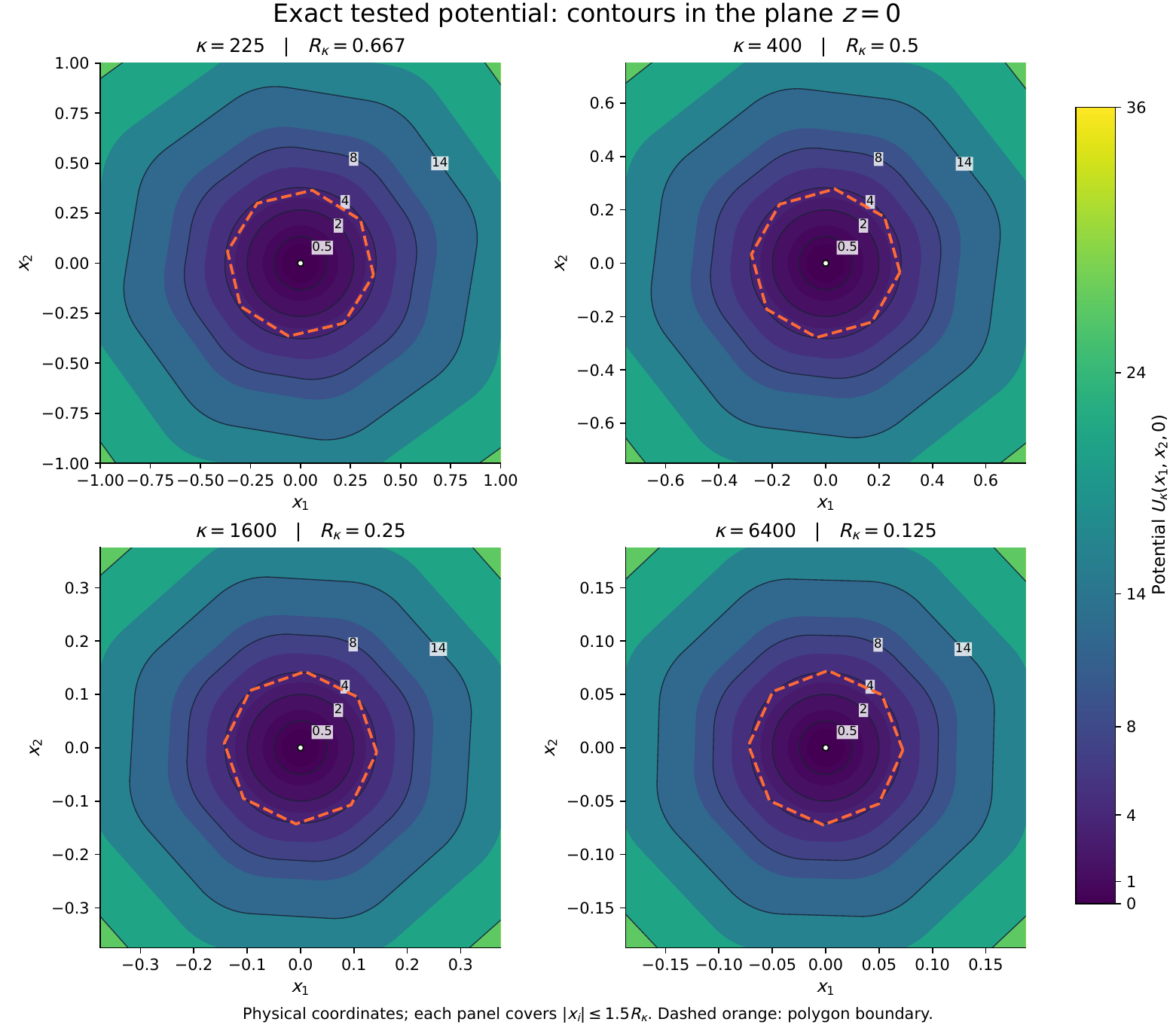}
    \caption{\textbf{Experiment of Subsection~\ref{nesterov_sha_subsection4.1}} Contour plot of potential function \eqref{neste_eq:supp_five_potential} with $z = 0$ for different values of $\kappa$. Note the different axis limits in the subplots for each $\kappa$.}
    \label{nesterov_sharma_contour_plot_poenfunc}
\end{figure}

We consider the three-dimensional Gibbs measure, $
\mu_\kappa(\der y) =Z_\kappa^{-1}\exp\{-U_\kappa(y)\}\der y$, $y=(x,z)\in\mathbb R^2\times\mathbb R$, parameterized by $\kappa$, with potential function
\begin{align}
U_\kappa(x,z)
&=\frac{R_\kappa^2}{4}\psi_\kappa(x/R_\kappa)
+\frac{z^2}{8}+\log\cosh(z/2),
\qquad R_\kappa=\frac{10}{\sqrt\kappa}.
\label{neste_eq:supp_five_potential}
\end{align}
 The construction of $\psi_\kappa(x)$, detailed below, is borrowed from \citep[Theorem~3.5]{goujaud2025provable} (see also \citep{bou2026provable}). We add the scalar potential $z^2/8+\log\cosh(z/2)$ to introduce
a non-Gaussian direction whose curvature remains between $1/4$
and $1/2$, independently of $\kappa$. Here
$\cosh(t)=(e^t+e^{-t})/2$ denotes the hyperbolic cosine. 
This component allows us to assess exploration of a broad
direction when the step size is constrained by the stiff
polygon component.
Let $Q\in\mathbb R^{2\times2}$ denote the matrix corresponding
to a counterclockwise rotation by $\pi/4$, i.e. 
\begin{align}
Q=\frac{1}{\sqrt{2}}
\begin{pmatrix}
1 & -1\\
1 & 1
\end{pmatrix}
\end{align} 
Let $e_j=(\cos(j\pi/4),\sin(j\pi/4))^\top$,  $j=0,\ldots,7$, 
and define
\begin{align}\label{nesterov_sharma_eqn4.2}
\eta_{\mathrm{geom},\kappa}
&=\frac{\sqrt\kappa-1}{\sqrt\kappa+1},
\qquad
M_\kappa=
\frac{(1+\eta_{\mathrm{geom},\kappa}-\kappa^{-1})I_2
-Q-\eta_{\mathrm{geom},\kappa}Q^\top}{1-\kappa^{-1}},
\end{align}
\begin{align}
C_\kappa&=\operatorname{conv}\{M_\kappa e_j:0\leq j\leq7\},
\qquad
\psi_\kappa(u)=\frac\kappa2|u|^2
-\frac{\kappa-1}{2}\operatorname{dist}(u,C_\kappa)^2.
\label{eq:supp_five_polygon}
\end{align}
Here $\operatorname{conv}$ denotes the convex hull and
$\operatorname{dist}$ is Euclidean distance.  Figure~\ref{nesterov_sharma_contour_plot_poenfunc} shows the contour plots of the potential function for $z = 0$.

We measure coupling distance by the root mean square distance (see \eqref{eq:supp_five_coupling_nesterov_sharma_metric})
between coupled chains ($8192$ pairs per starting point),
normalized by its initial value and maximized over $3$
starting points. We obtain the  starting points of the reference chains by running NS \eqref{nesterov_sharma_ns-position-recurrence}, OBABO
and BAOAB until certain equilibration checks pass, and by rejection
sampling from the Gibbs distribution for HMC and MALA.
Each coupled pair is then driven by the same random variables. 

\begin{figure}[t]
\centering
\begin{subfigure}{\linewidth}
    \centering
    \includegraphics[width=\linewidth, height=0.25\textheight]{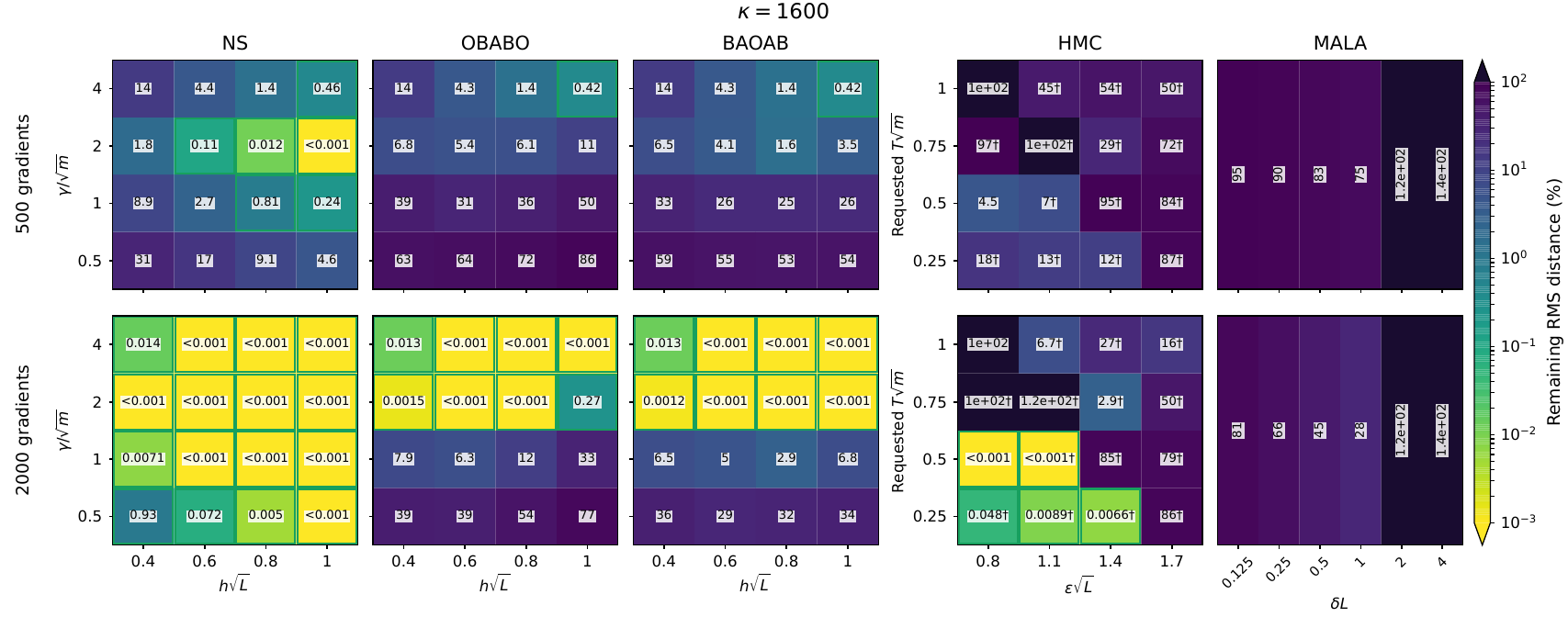}
    \phantomsubcaption
    \label{fig:heatmaps_kappa1600}
\end{subfigure}
\par\nointerlineskip
\vspace{-1em}
\begin{subfigure}{\linewidth}
    \centering
    \includegraphics[width=\linewidth, height=0.25\textheight]{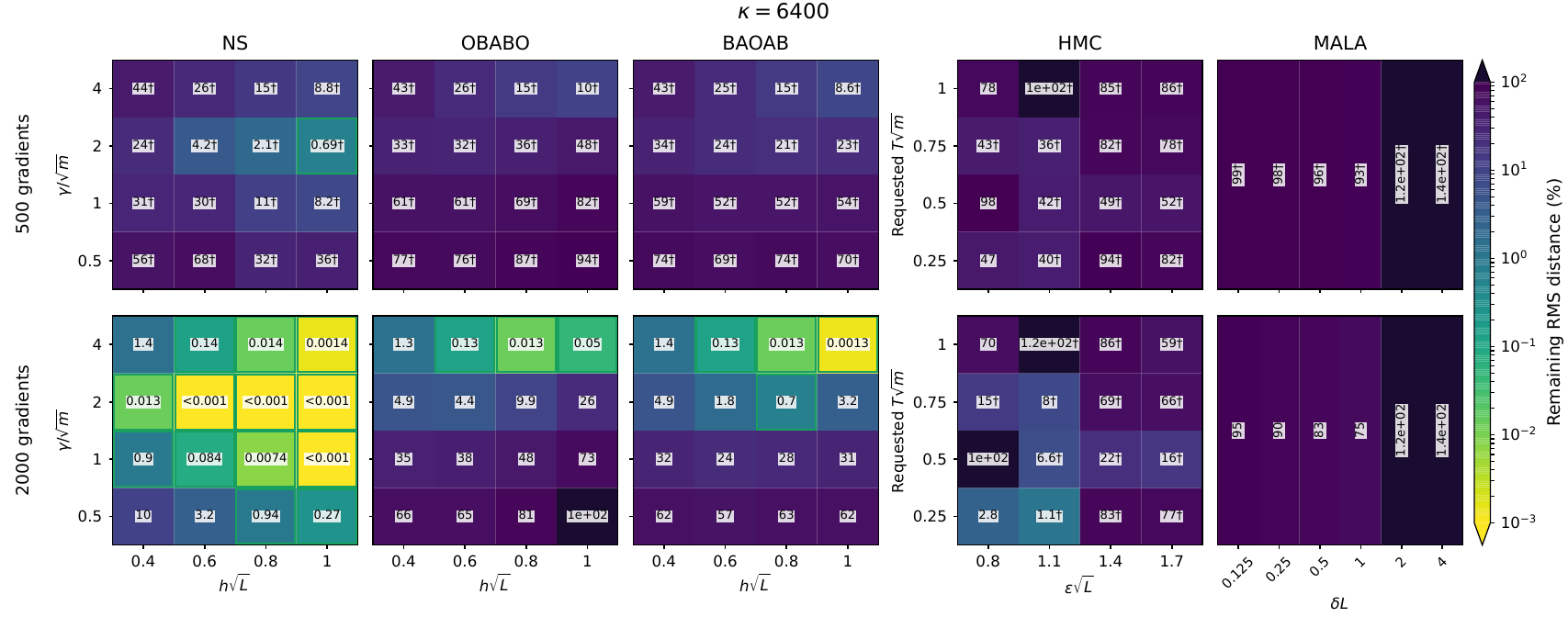}
    \phantomsubcaption
    \label{fig:heatmaps_kappa6400}
\end{subfigure}

\captionsetup{font=small}\caption{\textbf{Experiment of Subsection~\ref{nesterov_sha_subsection4.1}}. Normalized root mean square synchronous coupling
distance, expressed as a percentage, for $\kappa=1600$
and $6400$ from top to bottom.
Green outlines indicate distances at or below $1\%$.
Dagger indicate that the displayed value was achieved before the stated gradient budget is reached.}
\label{nestorov_sharma_figure_coupling_heatmaps_selected_kappas}
\end{figure}

Within a budget of 500 gradient evaluations per path, NS \eqref{nesterov_sharma_ns-position-recurrence} reduces
the normalized root mean square coupling distance below $1\%$
at $6$ of the 16 tested parameter settings for $\kappa=1600$ as compared to one each for OBABO and BAOAB, and none for HMC and MALA.
At $\kappa=6400$, NS is the only sampler to achieve this in coupling distance
within the tested grids, with a remaining distance of $0.69\%$
at $h L^{1/2}=1$ and $\gamma/m^{1/2}=2$. The detailed results are presented in Figure~\ref{nestorov_sharma_figure_coupling_heatmaps_selected_kappas}. For HMC, a wider pilot search was conducted to select
the HMC parameter grid which we have shown in Figure~\ref{nestorov_sharma_figure_coupling_heatmaps_selected_kappas}.

\begin{figure}[htbp!]
\centering
\includegraphics[width=\linewidth, height = 0.40 \textheight]{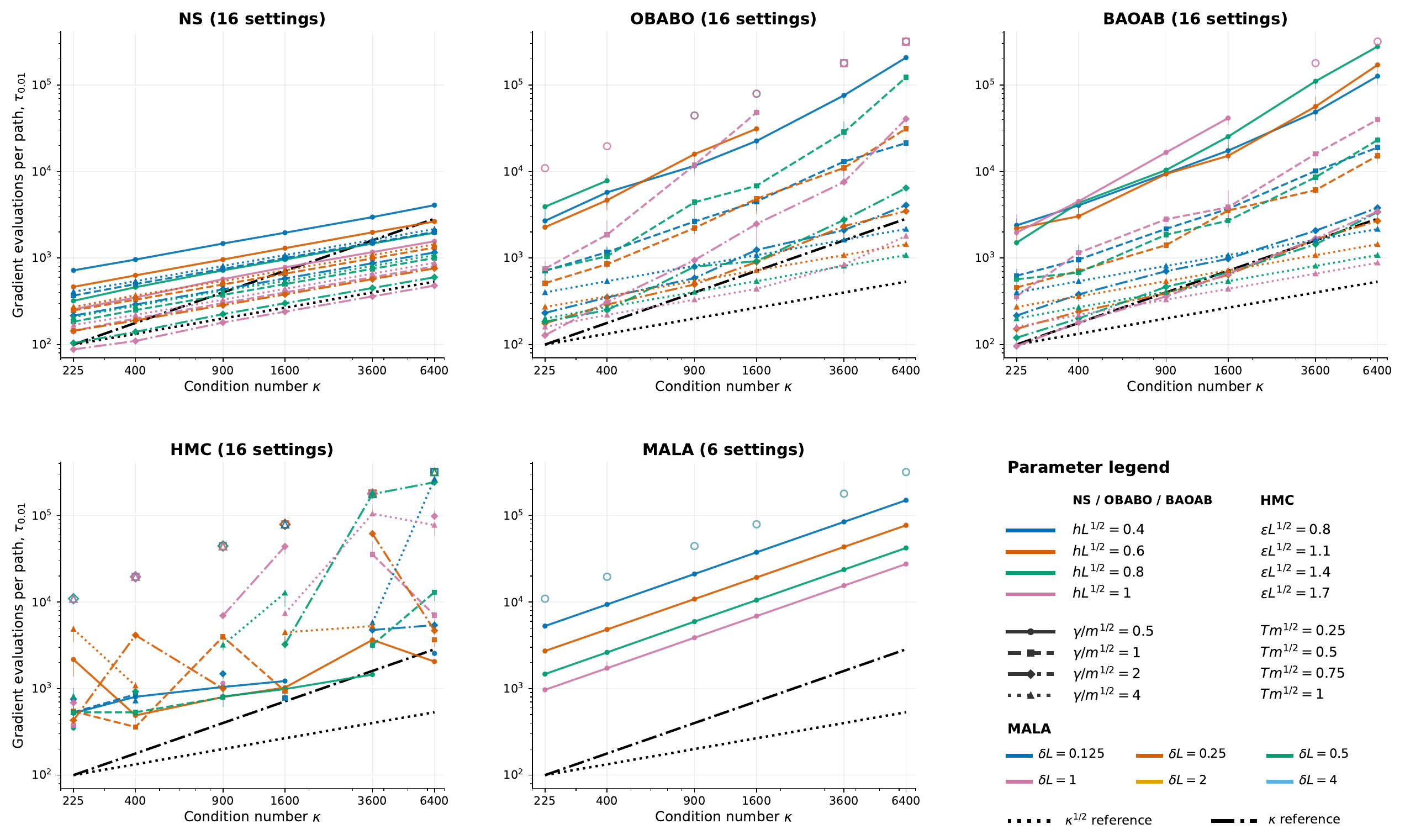}
\captionsetup{font=small}
\caption{\textbf{Experiment of Subsection~\ref{nesterov_sha_subsection4.1}}. Relaxation cost $\tau_{0.01}$ against condition number
for all tested parameter settings of the five samplers.
Both axes are logarithmic. Black dotted and dot-dashed lines
have slopes $1/2$ and $1$, respectively.
For unresolved thresholds, open markers show the last recorded
cost allowing a full persistence check.
}
\label{neste_shar_supp_five_nesterov_sharma_scaling_all}
\end{figure}

   We denote by $\tau_{0.01}$ the first recorded gradient cost
from which this distance stays at or below $1\%$ at every saved
observation for at least $\lceil400(\kappa/400)^{1/2}\rceil$
further gradient evaluations per path.  We observe square-root growth of $\tau_{0.01}$
with $\kappa$ for proposed Nesterov's method based sampler for all $16$ settings, compared with OBABO and BAOAB which showed square root growth for 3 and 4 of the $16$ settings, respectively. MALA scaled linearly for 4 of 6 settings while other 2 remained unresolved (see Figure~\ref{neste_shar_supp_five_nesterov_sharma_scaling_all}).
Accelerated relaxation is known for continuous-time kinetic
Langevin dynamics \citep{cao2023explicit,eberle2026non, borkowski2026fenchel} under suitable assumptions, but this does not directly establish acceleration for a discretization.
Our construction, and therefore this paper, addresses the
question posed in the introduction of obtaining a
Nesterov-type sampler with second order invariant accuracy. Sublinear sampling costs in $\kappa$ have been established for randomized midpoint schemes
\citep{altschuler2026shifted}, while square-root dependence has been obtained for proximal bouncy particle sampling
from a warm start \citep{chen2026accelerated} and
Metropolis-adjusted kinetic Langevin sampling under an
additional Lipschitz Hessian assumption
\citep{bou2026persistence}.
Whether our sampler also admits an acceleration guarantee
beyond the present example remains to be established.

\begin{figure}[htbp!]
    \centering
\includegraphics[width=0.95\linewidth, height = 0.35\textheight]{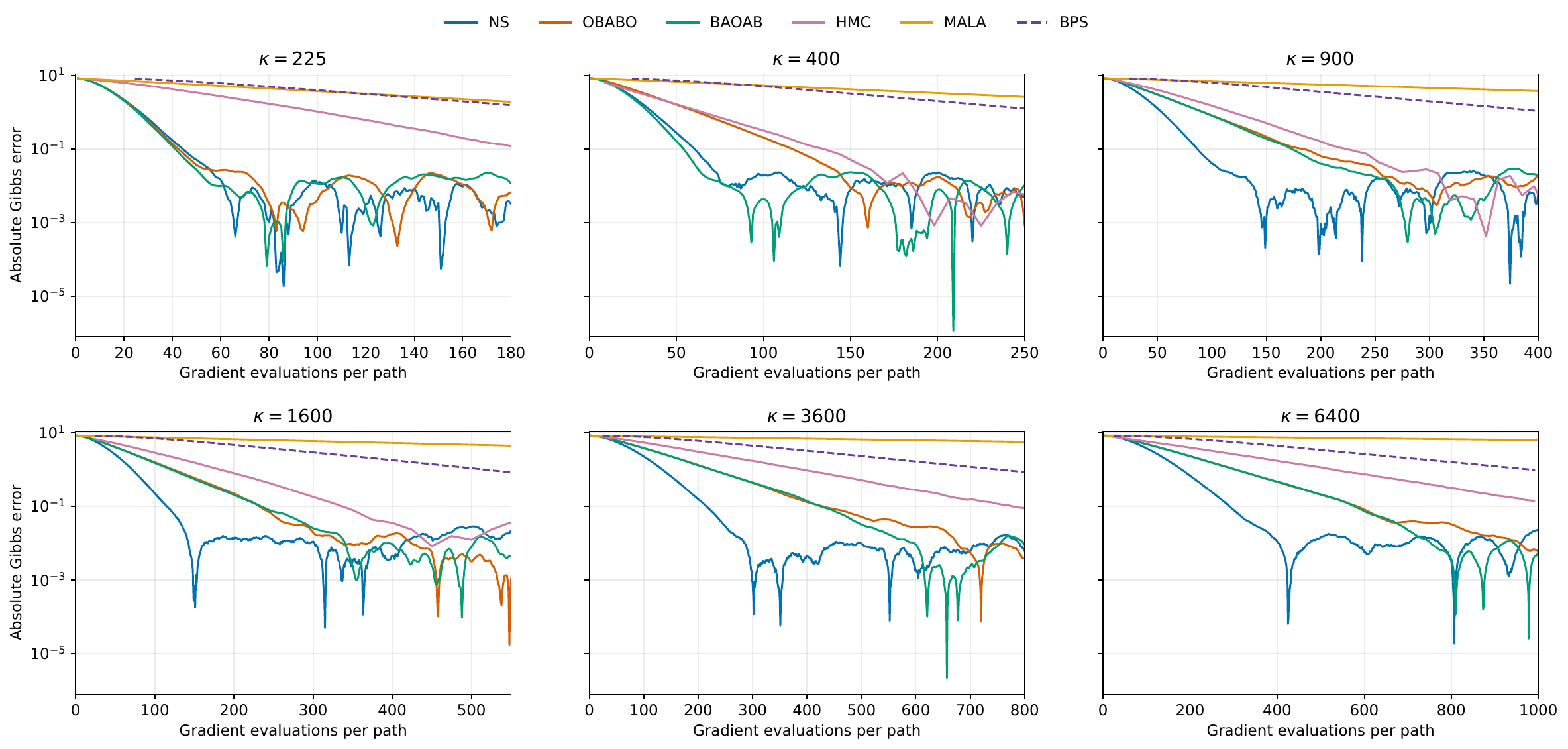}
  \captionsetup{font=small}  \caption{\textbf{Experiment of Subsection~\ref{nesterov_sha_subsection4.1}}. Absolute error of the instantaneous ensemble estimate
    of the Gibbs expectation of $\varphi(x,z)=(z/2)^2$ using     8192 trajectories for each method and condition number. The horizontal axes show gradient evaluations per path (for BPS, we average this count over all trajectories). The vertical axes are logarithmic. }
\label{fig:six_sampler_scalar_gibbs_error_nesterov_sharma}
\end{figure}

We also compare the instantaneous ensemble estimates of test function
$\varphi(x,z)=(z/2)^2$ with its Gibbs expectation
$\mu_\kappa(\varphi)$, the target ergodic limit, in
Figure~\ref{fig:six_sampler_scalar_gibbs_error_nesterov_sharma}. In this comparison, we also include bouncy particle sampler \citep{peters2012rejection,bouchard2018bouncy}. Before Monte Carlo fluctuations and discretization bias interfere,
we observe the faster convergence of OBABO, BAOAB and NS at
$\kappa=225$, comparable performance of BAOAB and NS at
$\kappa=400$, and a clearly faster convergence of NS at $\kappa=900,1600,3600,6400$.  For NS, OBABO and BAOAB, we first choose the friction coefficient which gives
the lowest resolved coupling cost for each method and then we
select a common step size minimizing the largest of the
three coupling costs.

\subsection{Comparison using a smoothed total-variation potential function}
\label{nesterov_sharma_tv_results}
In this section, we compare NS, OBABO and BAOAB on a nonquadratic
strongly convex potential in dimension $d=256$. This potential is based on smoothed total-variation used in the image-denoising model of \citet{rudin1992nonlinear},
using the square-root approximation studied in
\citet{vogel1996iterative,chan1999nonlinear}.
We identify $y\in\R^{256}$ with a $16\times16$ array and define
\[
(D_1y)_{ij}=y_{i,j+1}-y_{ij},\qquad
(D_2y)_{ij}=y_{i+1,j}-y_{ij},
\]
with $(D_1y)_{i,16}=(D_2y)_{16,j}=0$.
For $\delta=0.05$, let
\begin{align}
U_\kappa(y)
&=\frac12|y-b|^2+\frac{\kappa-1}{B_\delta}V_\delta(y),
\label{nesterov_sharma_tv_target}\\
V_\delta(y)
&=\sum_{i,j=1}^{16}
\left\{\delta^2+(D_1y)_{ij}^2+(D_2y)_{ij}^2\right\}^{1/2},
\qquad B_\delta=\frac{8\cos^2(\pi/32)}{\delta}.
\label{nesterov_sharma_tv_regularizer}
\end{align}
The quadratic term keeps the array close to $b$, while
$V_\delta$ penalizes differences between neighbouring entries.
We take
\[
b_{ij}=0.3+\mathbf 1_{\{j\geq9\}}
+0.3\mathbf 1_{\{i\geq9\}}+0.1G_{ij},
\]
where the entries of $G$ are independent standard Gaussians,
and we fix them for the experiment. The potential  function \eqref{nesterov_sharma_tv_target} satisfies
$I_d\preceq\nabla^2U_\kappa\preceq\kappa I_d$.
Both bounds are attained at constant arrays, giving
$m=1$ and $L=\kappa$.

\begin{figure}[htbp!]
\centering
\includegraphics[width=\linewidth, height = 0.30\textheight]{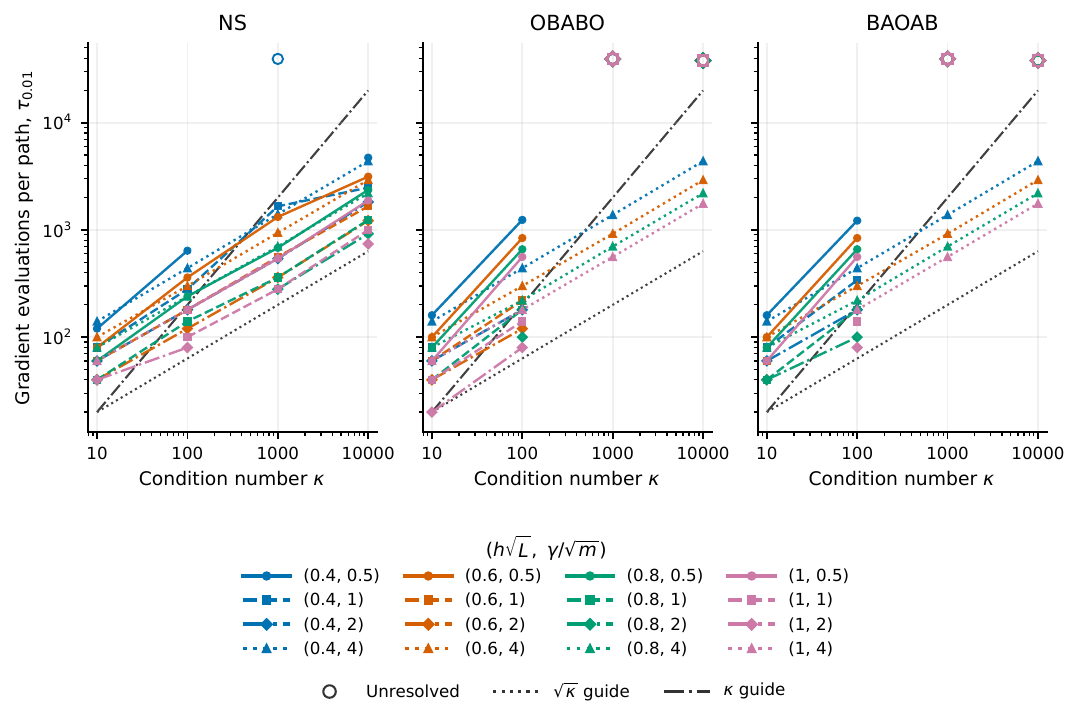}
\caption{\textbf{Experiment of Subsection~\ref{nesterov_sharma_tv_results}}. Relaxation cost $\tau_{0.01}$ for the smoothed
total-variation potential \eqref{nesterov_sharma_tv_target},
using all 16 parameter choices. Both axes are logarithmic.
Open markers show the last recorded cost allowing a full
persistence check for unresolved thresholds.
Dotted and dot-dashed lines have slopes $1/2$ and $1$,
respectively.}
\label{nesterov_sharma_tv_all_settings}
\end{figure}

We use $\kappa\in\{10,100,1000,10000\}$ and the same
16 parameter choices
\eqref{eq:supp_five_nesterov_sharma_kinetic_grid}
for all three methods. Computational details are given in
Appendix~\ref{nesterov_sharma_tv_comparison}.
Figure~\ref{nesterov_sharma_tv_all_settings} shows the
relaxation costs $\tau_{0.01}$  within a budget of 40000
new gradient evaluations per path.
At $\kappa=10$ and $100$, every setting with a qualified
reference reaches the threshold.
At $\kappa=1000$ and $10000$, NS reaches it at 14 and 15
settings, respectively, whereas OBABO and BAOAB each reach
it at only four settings out of $16$. All resolved settings for the two
splitting methods at these condition numbers use
$\gamma/\sqrt m=4$.
At matching parameters where all three methods resolve
the threshold, their costs are similar. NS also resolves
the threshold for many choices with smaller friction.
A setting is unresolved when its reference is qualified (i.e. the reference simulations pass the prescribed checks
for stationarity) but the normalized root mean square coupling distance does not remain below
$0.01$ for the required interval within the computational budget.  The recorded distances for unresolved settings level off above this threshold and this describes separation
under synchronous coupling but does not imply that the individual chains fail to converge in distribution.
For the single unresolved NS setting, this occurs for the start at the minimizer, whereas for other two starting points NS reach the threshold (note that we measure coupling distance  by maximizing over three different starting points).

Table~\ref{nesterov_sharma_tv_best} gives the lowest observed
costs among the tested settings.
At $\kappa=1000$ and $10000$, the lowest costs for OBABO
and BAOAB are approximately 2 and 2.38 times the
corresponding NS costs.

\begin{table}[t]
\centering
\caption{\textbf{Experiment of Subsection~\ref{nesterov_sharma_tv_results}}. Lowest observed relaxation costs among qualified,
resolved settings in
Figure~\ref{nesterov_sharma_tv_all_settings}, in terms of gradient
evaluations per path. Parameters are selected separately
for each method and condition number over the tested grid.}
\label{nesterov_sharma_tv_best}
\begin{tabular}{r r r r}
\toprule
$\kappa$ & NS & OBABO & BAOAB\\
\midrule
10    & 40  & 20   & 40\\
100   & 80  & 80   & 80\\
1000  & 280 & 560  & 560\\
10000 & 740 & 1760 & 1760\\
\bottomrule
\end{tabular}
\end{table}
\subsection{Bayesian logistic regression for mortgage applications}\label{nesterov_sha_subsection4.2}

We consider home-purchase applications in North Carolina from
the 2025 Home Mortgage Disclosure Act data \citep{ffiec2026hmda}. The application denial rate in original data is approximately $13.4\%$.
We keep all denied applications and randomly select an equal
number of undenied applications separately in each set, giving
$38516$ training and $9630$ test observations with a $50\%$
accuracy baseline. We set the response to $1$ for denial and $0$ otherwise. We represent categories by zero-one variables, and transform and scale numerical predictors to put them on comparable scales.
Missing values are replaced by their medians
in the training set.  We model the denial probability by
$\operatorname{logit}(p_i)=x_i^\top\beta+b_{\ell(i)}$, where
$\ell(i)$ identifies the lender.
The model keeps 32 features from the data and
add 1003 lender effects.
Independent centred Gaussian priors have variances $25$,
$6.25$ and $1$ for the intercept, remaining regression
coefficients and lender effects, respectively.

\begin{figure}[!htbp]
\centering
\includegraphics[width=0.9\linewidth, height = 0.18\textheight]{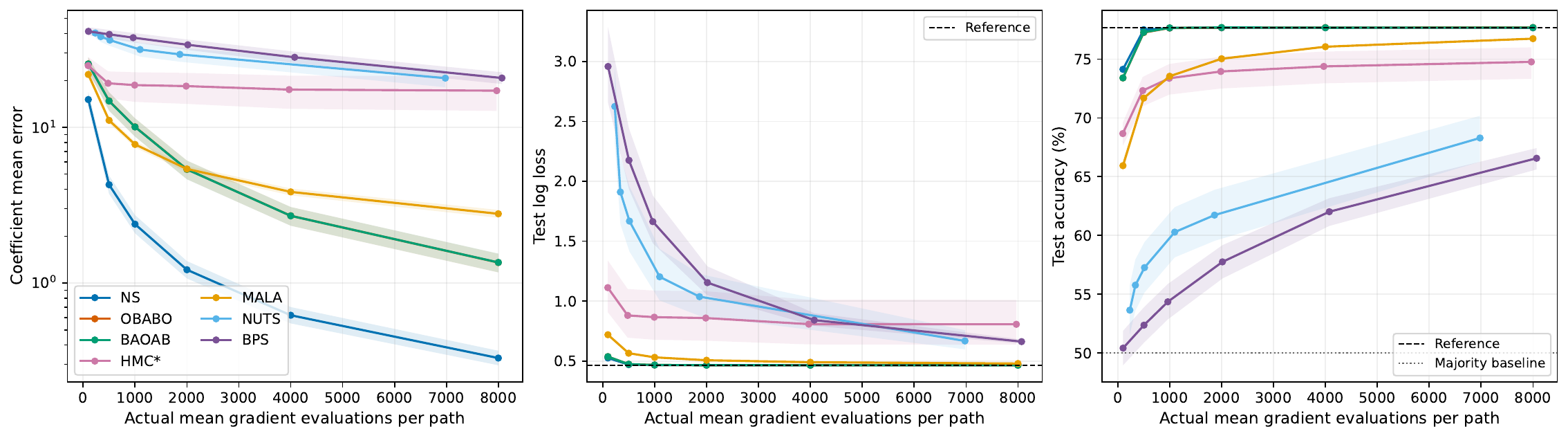}
 \captionsetup{font=small} \caption{\textbf{Experiment of Subsection~\ref{nesterov_sha_subsection4.2}}. Coefficient mean error, test log loss and accuracy
against average gradient cost, starting at a budget of 100 gradients. Shaded region highlights $95\%$
trajectory-bootstrap intervals conditional on the reference
and fixed settings. 
NUTS costs include gradient evaluations on inactive batch entries.}
\label{fig:hmda_seven_samplers_nesterov_sharma}
\end{figure}

In this experiment, we also include the No-U-Turn sampler (NUTS)
\citep{hoffman2014no}. 
All seven samplers use the same $64$ initial positions sampled
from the prior distribution.
For NS, OBABO and BAOAB, we use
$h=0.0030733$ and $\gamma=1$.
MALA uses $\delta=3.2795\times10^{-5}$, and HMC uses
$\epsilon=h$ with 32 leapfrog steps.
NUTS uses step size $0.0040469$ and maximum tree depth 12.
BPS uses windowed thinning \citep{lu2026windowed} with window
length $0.040977$ and Gaussian velocity refreshment at rate $32$.
The parameters used for each method are selected by grid search,
with common $h$ and $\gamma$ imposed for NS, OBABO and BAOAB.
We compute the posterior reference estimates by averaging $8\times 1024$ samples from eight NUTS chains, after discarding 1000 warmup iterations per chain  (the maximum $\widehat R$ is $1.0054$). The references' computations used 8,307,480 gradient evaluations.
We divide computed coefficient differences from the
reference means by their reference posterior standard deviations
and report the root mean square error  for each sampler.
We also report test log loss and classification accuracy,
using each trajectory's averaged predictive probabilities (see Figure~\ref{fig:hmda_seven_samplers_nesterov_sharma}).

At 8000 gradient evaluations budget, the coefficient errors are
$0.330$, $1.356$ and $1.357$ for NS, OBABO and BAOAB,
respectively.
Their test accuracies are approximately $77.7\%$,
close to the reference. Therefore, NS gives smaller coefficient errors
for averages including the initial transient, with similar
predictive performance. NUTS is plotted up to a mean cost of 6979 gradients, with 12 iterations chosen (based on a pilot run) beforehand  for $64$ trajectories to target 8000 gradients.

\section*{Declaration of the use of generative AI and AI-assisted technologies}
The research ideas and main proof strategies are the author's.
OpenAI's ChatGPT was used under the author's guidance,
primarily for coding and to assist with proofs for routine calculations, and with
exposition and preparation of the \LaTeX\ source. All outputs generated by the LLM were verified by the author,
who takes full responsibility for the content of the paper.

\section*{Acknowledgment}
The numerical experiments were enabled by resources provided by Chalmers e-Commons at Chalmers as well as the National Academic Infrastructure for Supercomputing in Sweden (NAISS), partially funded by the Swedish Research Council through grant no. 2022-06725.

{\small 
\bibliographystyle{plainnat}
\bibliography{reference}
}

\appendix
\counterwithin{theorem}{section}

\section{Proofs of Section~\ref{nesterov_sharma_section_3}}
\subsection{Proof of Theorem~\ref{nesterov_sharma_theorem_geometric_ergodicity}}

\begin{proof}[Proof of Lemma~\ref{nesterov_lem_uniform-moment-bounds_sharma}]
Let
$m=\gamma^2\beta(2-\beta)/(4(1-\beta))>0$. We have
\begin{align}
0
\leq U(y)
\leq C(1+|y|^2), \qquad
( y \cdot F(y))
\leq-m|y|^2 +2b,
\qquad
|F(y)-F(z)|\leq L|y-z|.
\label{eq:moment-potential-bounds}
\end{align}
Throughout the proof, constants denoted by $C$ may change
from place to place but are independent of $h$, $n$ and the current state. Let us consider the following Lyapunov function \citep{mattingly2002ergodicity}:
\begin{align}
\mathcal{V}(y,v)
=1+U(y)+\frac14|v|^2
+\frac14|v+\gamma y|^2.
\label{eq:moment-energy}
\end{align}
Since $U\geq0$, the quadratic terms in
$\mathcal V$ control both $y$ and $v$.
There exist constants $c_1,c_2>0$ such that
\begin{align}
c_1(1+|y|^2+|v|^2)
\leq\mathcal V(y,v)
\leq c_2(1+|y|^2+|v|^2).
\label{eq:moment-energy-comparison}
\end{align}
Moreover, $\nabla\mathcal V$ is globally Lipschitz and therefore
\begin{align}
|\nabla\mathcal{V}(y,v)|
\leq C(1 + |y| + |v|).
\label{eq:moment-energy-gradient}
\end{align}
Consequently,
\begin{align}
\mathcal{V}(y+\delta y,v+\delta v)
\leq &
\mathcal{V}(y,v)
+(\nabla_y\mathcal{V}(y,v)\cdot\delta y)
+(\nabla_v\mathcal{V}(y,v)\cdot\delta v)
\nonumber\\
&+C\bigl(|\delta y|^2+|\delta v|^2\bigr).
\label{eq:moment-taylor-bound}
\end{align}
Using \eqref{eq:moment-potential-bounds}, we have
\begin{align}
&\big(\nabla_y\mathcal{V}(y,v)\cdot v\big)
+ \big(\nabla_v\mathcal V(y,v)\cdot\{F(y)-\gamma v\}\big)
=-\frac{\gamma}{2}|v|^2
+\frac{\gamma}{2}( y \cdot F(y))
\nonumber\\
&\qquad
\leq-\frac{\gamma}{2}|v|^2
-\frac{\gamma m}{2}|y|^2 + \gamma b.
\label{nesterov_sampler_eqn_moment_drift}
\end{align}
The noise coefficients satisfy $\eta\theta_0+\theta_1
=(1-\eta)(1-\theta_0)$. 
Define $\sigma_h$ by
\begin{align}
h\sigma_h
=\sqrt{\frac{2h}{\gamma}}\,
(1-\eta)(1-\theta_0).
\label{eq:moment-velocity-noise}
\end{align}
For fixed $\gamma>0$ and sufficiently small $h$, the coefficient
formulas implies the following bounds:
\begin{align}
|1-\eta-\gamma h| \leq Ch^2,
\qquad |\theta_0|\leq C\sqrt h,
\qquad |\sigma_h|\leq C\sqrt h.
\label{eq:moment-coefficient-bounds}
\end{align}
Condition on $(Y_n,V_n)=(y,v)$. The updates can be written as
\begin{align}
Y_{n+1}-y
&=hv+\sqrt{\frac{2h}{\gamma}}\,
\theta_0\xi_{n+1},
\nonumber\\
V_{n+1}-v
&=h\big(F(y)-\gamma v\big)
+\sigma_h\xi_{n+1}+R_{n+1},
\label{eq:moment-increments}
\end{align}
where
\begin{align}
R_{n+1}
={}&(1-\eta-\gamma h)
\left(\frac{F(y)}{\gamma}-v\right)
+\frac{1-\eta^2}{\gamma}
\{F(Y_{n+1})-F(y)\}.
\label{eq:moment-increment-remainder}
\end{align}
By Lipschitz continuity, we have
\begin{align}
|F(Y_{n+1})-F(y)|
&\leq
L\left(
h|v|+\sqrt{\frac{2h}{\gamma}}\,
|\theta_0|\,|\xi_{n+1}|
\right) \leq Ch\bigl(|v|+|\xi_{n+1}|\bigr).
\label{eq:moment-force-increment}
\end{align}
Together with \eqref{eq:moment-coefficient-bounds} and
$|F(y)|\leq|F(0)|+L|y|\leq C(1+|y|)$, we obtain
\begin{align}
|R_{n+1}|
&\leq Ch^2\bigl(1+|y|+|v|+|\xi_{n+1}|\bigr),
\label{eq:moment-remainder-bound}\\
|Y_{n+1}-y|^2+|V_{n+1}-v|^2
&\leq Ch^2(1+|y|^2+|v|^2)
+Ch|\xi_{n+1}|^2.
\label{eq:moment-squared-increments}
\end{align}
Let $\mathcal{V}_n := \mathcal{V}(Y_n,V_n)$ and define
\begin{align}
M_{n+1} =\sqrt{\frac{2h}{\gamma}}
\left\{ \theta_0\nabla_y\mathcal V(Y_n,V_n)
+\frac{(1-\eta)(1-\theta_0)}{h}
\nabla_v\mathcal V(Y_n,V_n)
\right\}\cdot\xi_{n+1}.
\label{eq:moment-martingale-increment}
\end{align}
Substituting \eqref{eq:moment-increments} into
\eqref{eq:moment-taylor-bound}, and using
\eqref{nesterov_sampler_eqn_moment_drift},
\eqref{eq:moment-remainder-bound} and
\eqref{eq:moment-squared-increments}, yield
\begin{align}
\mathcal{V}_{n+1} \leq{}&
\mathcal V_n -h\left(
\frac{\gamma}{2}|V_n|^2
+\frac{\gamma m}{2}|Y_n|^2
\right) +M_{n+1} + \gamma b h
\nonumber\\
&+Ch^2({1+}|Y_n|^2+|V_n|^2)
+Ch|\xi_{n+1}|^2.
\label{eq:moment-energy-increment}
\end{align}
We choose $h_0\leq1$ sufficiently small and there is a
constant $c>0$ (see \eqref{eq:moment-energy-comparison}) such that
\begin{align}
\mathcal V_{n+1}
\leq
(1-ch)\mathcal{V}_n
+M_{n+1}+Ch(1+|\xi_{n+1}|^2)
\label{eq:moment-pathwise-drift}
\end{align}
for $0<h\leq h_0$. We require $ch_0\leq1$. Let $\mathcal F_n$ be the sigma-algebra generated by
the initial state and the Gaussian noise up tp step $n$. The coefficient of $\xi_{n+1}$ in
\eqref{eq:moment-martingale-increment} is
$\mathcal F_n$-measurable. The vector $\xi_{n+1}$ is independent of $\mathcal F_n$ and has mean zero. Hence
\begin{align}
\E[M_{n+1}\mid\mathcal F_n]=0,
\label{eq:moment-centering}
\end{align}
and 
boundedness of Gaussian moments as well as
\eqref{eq:moment-energy-gradient}-\eqref{eq:moment-coefficient-bounds} yield, for every
integer $j\geq1$,
\begin{align}
\E\bigl[|M_{n+1}|^j\mid\mathcal F_n\bigr]
\leq C_jh^{j/2}\mathcal{V}_n^{j/2}.
\label{eq:moment-martingale-bound}
\end{align}
Taking conditional expectations in
\eqref{eq:moment-pathwise-drift}, we ascertain
\begin{align}
\E[\mathcal{V}_{n+1}\mid\mathcal F_n]
\leq(1-ch)\mathcal{V}_n+Ch.
\label{eq:moment-first-power}
\end{align}
For an integer $p\geq2$, take $p$-th power of
\eqref{eq:moment-pathwise-drift}.
Its right-hand side is nonnegative, since 
$\mathcal V_{n+1}\geq1$.
In the binomial expansion, the term linear in
$M_{n+1}$ has zero conditional expectation.
For the remaining terms, we have
\begin{align}
&\E[|M_{n+1}+Ch(1+|\xi_{n+1}|^2)|^j
\mid\mathcal F_n] \leq C_j\bigl(h^{j/2}\mathcal V_n^{j/2}+h^j\bigr),
\qquad 2\leq j\leq p,
\label{nesterov_sampler_eq:moment-binomial_bound}
\end{align}
where $C_j$ depends on $p$ but independent of $h$. 
Since $h\leq1$ and $\mathcal V_n\geq1$, multiplication
by $\mathcal V_n^{p-j}$ bounds these terms by
$C_ph\mathcal V_n^{p-1}$.
The term linear in $Ch(1+|\xi_{n+1}|^2)$ satisfies
the same bound. Therefore, we get
\begin{align}
\E[\mathcal V_{n+1}^p\mid\mathcal F_n]
\leq
(1-ch)^p\mathcal V_n^p
+C_ph\mathcal V_n^{p-1}.
\label{eq:moment-power-drift}
\end{align}

Using $(1-ch)^p\leq1-ch$ and Young's inequality (i.e. $C_px^{p-1}\leq(c/2)x^p+K_p$,  $x\geq1$), 
we obtain
\begin{align}
\E[\mathcal V_{n+1}^p\mid\mathcal F_n]
\leq
\left(1-\frac c2h\right)\mathcal V_n^p+K_ph.
\label{eq:moment-lyapunov-drift}
\end{align}
For $p=1$, the same bound follows from \eqref{eq:moment-first-power} with $K_1=C$.  Taking expectations on both sides and iterating gives
\begin{align}
\E\mathcal{V}_n^p
&\leq
\left(1-\frac c2h\right)^n\E\mathcal V_0^p
+K_ph\sum_{j=0}^{n-1}
\left(1-\frac{c}{2}h\right)^j
\nonumber\\
&\leq
e^{-cnh/2}\E\mathcal V_0^p+\frac{2K_p}{c}.
\label{eq:moment-energy-bound}
\end{align}
Using \eqref{eq:moment-energy-comparison}, we obtain \eqref{eq:moment-estimate}.
\end{proof}

\begin{proof}[Proof of Lemma~\ref{nesterov_sharma_lem_minorization}]
Since a single noise vector is imparted into two $\R^d$-valued variables, we need to consider two steps to establish the minorization. To this end, fix $(Y_0,V_0)=(y,v)$ and consider the following two-step procedure:
\begin{align}
\begin{aligned}
Y_1&=y+hv+r_1\xi_1,\\
hV_1&=\eta hv+\alpha\{(1+\eta)F(Y_1)-\eta F(y)\}+r_1r_2\xi_1,\\
Y_2&=Y_1+hV_1+r_1\xi_2,\\
hV_2&=\eta hV_1+\alpha\{(1+\eta)F(Y_2)-\eta F(Y_1)\}+r_1r_2\xi_2,
\end{aligned}
\label{nesterov_sharma_eq_minorization_two_steps}
\end{align}
where $\eta=e^{-\gamma h}$, $\alpha=h(1-\eta)/\gamma$ and we have denoted
\begin{align*}
r_1=\left(\frac{2h}{\gamma}\right)^{1/2}\theta_0,
\qquad
r_2=\frac{\eta\theta_0+\theta_1}{\theta_0}.
\end{align*}
Note that $r_1>0$ since $\theta_0>0$. Let us also fix a final state $(Y_2,V_2)=(y',v')$. Then \eqref{nesterov_sharma_eq_minorization_two_steps} is a system of four equations in the four unknowns $\xi_1$, $\xi_2$, $Y_1$ and $V_1$. If there were no nonlinearity in $Y_1$ due to $F(Y_1)$, this system would be linear and we could directly apply the change of variables formula. Therefore, we first eliminate $\xi_1$, $\xi_2$ and $V_1$, which gives
\begin{align}
\mathcal T(Y_1)=w,
\label{nesterov_sharma_eq_minorization_recover_position}
\end{align}
where $\mathcal{T}(x)=a_1x-b_1F(x)$ with
$a_1=r_2(1+r_2-\eta)$ and
$b_1=\alpha\{\eta^2-(1+\eta)r_2\}$,
and $w$ depends only on the initial and final states:
\begin{align*}
w=(\eta-r_2)\{(\eta-r_2)hv-r_2y-\alpha\eta F(y)\}-hv'+r_2y'+\alpha(1+\eta)F(y').
\end{align*}
We first show that $\mathcal T$ is a bijection with continuously differentiable inverse. For this, we use the facts that $F=-\nabla U$ is continuously differentiable and $L$-Lipschitz, where $L$ is a Lipschitz constant of $\nabla U$, so that $\|\nabla^2U(x)\|\leq L$ for all $x\in\R^d$. For small $h$, the definitions of $\theta_0$ and $\theta_1$ show that $r_2$ behaves like $(\gamma h)^{1/2}$. Consequently, $a_1$ is of order $\gamma h$, and $b_1$ is $\mathcal{O}(h^2)$. Therefore, we have $a_1-L|b_1|>0$ for all sufficiently small $h>0$. For any $w\in\R^d$, the equation $\mathcal T(x)=w$ is equivalent to
\begin{align*}
x=\frac{w+b_1F(x)}{a_1},
\end{align*}
whose right-hand side has Lipschitz constant $L|b_1|/a_1<1$. The contraction mapping theorem therefore gives a unique solution, so $\mathcal T$ is a bijection. Moreover,
\begin{align*}
D\mathcal T(x)=a_1I_d+b_1\nabla^2U(x)\succeq(a_1-L|b_1|)I_d,
\end{align*}
so $D\mathcal T(x)$ is invertible for every $x$, and the inverse function theorem shows that $\mathcal T^{-1}$ is continuously differentiable.
Consequently, $Y_1=\mathcal T^{-1}(w)$, and using \eqref{nesterov_sharma_eq_minorization_two_steps}, we recover
\begin{align}
hV_1&=r_2Y_1+\alpha(1+\eta)F(Y_1)+(\eta-r_2)hv-r_2y-\alpha\eta F(y), \nonumber\\
\xi_1&=\frac{Y_1-y-hv}{r_1},
\qquad
\xi_2=\frac{y'-Y_1-hV_1}{r_1}. \label{nesterov_sharma_xi1_xi2_minorization}
\end{align}
Therefore, for each initial state, the map $(\xi_1,\xi_2)\mapsto(Y_2,V_2)$ is a continuously differentiable bijection of $\R^{2d}$ with continuously differentiable inverse $(y',v')\mapsto(\xi_1,\xi_2)$. By the change of variables formula, $(Y_2,V_2)$ has the density
\begin{align*}
p_h^{(2)}((y,v),(y',v'))
=(2\pi)^{-d}\exp\left(-\frac{|\xi_1|^2+|\xi_2|^2}{2}\right)|\det J|,
\end{align*}
where $J$ is the Jacobian matrix of the map $(y',v')\mapsto(\xi_1,\xi_2)$. Since $J$ is the inverse of the Jacobian matrix of
$(\xi_1,\xi_2)\mapsto(Y_2,V_2)$, we have $\det J\neq0$.
Hence, the density is  positive everywhere. It is also
continuous in $(y,v)$ and $(y',v')$, since $\xi_1$, $\xi_2$ (see \eqref{nesterov_sharma_xi1_xi2_minorization})
and $J$ are continuous in these variables.  Let $K$ be nonempty and let $B$ be the closed unit ball in $\R^{2d}$. Since $K\times B$ is compact, we have
\begin{align*}
\varepsilon_{h,K} =\operatorname{Leb}(B)
\inf_{(y,v)\in K,\,(y',v')\in B}
p_h^{(2)}((y,v),(y',v'))>0.
\end{align*}
With $\rho(A)=\operatorname{Leb}(A\cap B)/\operatorname{Leb}(B)$, we obtain, for all $(y,v)\in K$ and $A\in\mathcal B(\R^{2d})$,
\begin{align*}
P_h^2((y,v),A)
\geq\int_{A\cap B}p_h^{(2)}((y,v),(y',v'))\der y'\der v'
\geq\varepsilon_{h,K}\rho(A),
\end{align*}
which proves \eqref{nesterov_sharma_eq_minorization}.
\end{proof}

\subsection{Proof of Theorem~\ref{nesterov_sharma_theorem_2}}\label{nesterov_sharma_appendic_A}
For the convenience of the reader, we recall the sampler
\begin{align}
Y_1 &=Y_0+hV_0+ \sqrt{\frac{2h}{\gamma}}\,\theta_0\xi_1,
\label{nes_sh_eq:position-update}\\
hV_1  & =h\eta V_0
+ \frac{h(1-\eta)}{\gamma}\bigl\{(1+\eta)F(Y_1)-\eta F(Y_0)\bigr\} + \sqrt{\frac{2h}{\gamma}}\,(\eta\theta_0+\theta_1)\xi_1,
\label{nesh_sh_eq:velocity-update}
\end{align}
where $\eta = e^{-\gamma h}$ and
\begin{align}
\theta_0 =\frac{1-\eta}{2}
+\sqrt{\frac{\eta(1-\eta^2)}{2}},
\quad
\theta_1=\frac{1-\eta}{2}
-\sqrt{\frac{\eta(1-\eta^2)}{2}}.
\label{eq:coefficients}
\end{align}
Let $\nu_h$ represent the invariant measure for \eqref{nes_sh_eq:position-update}-\eqref{nesh_sh_eq:velocity-update}, and let $\mu_h$ denote its position marginal. In what follows, we need the following notation:
\begin{align}
a=\frac{(1-\eta)(\eta+\frac12)}{\gamma},
\qquad
hb= \frac{2(1-\eta)(1-\theta_0)^2}{\gamma(1+\eta)}.
\label{eq:reference-parameters}
\end{align}
For sufficiently small $h$, $a=\mathcal{O}(h)$ and $|\theta_0|+|\theta_1|=\mathcal{O}(\sqrt{h})$. It is not difficult to verify that the following holds: 
\begin{align*}
\theta_0+\theta_1&=1-\eta,
\theta_0\theta_1
=\frac{(1-\eta)^2}{4}
-\frac{\eta(1-\eta^2)}{2}.
\end{align*}
Substituting above in \eqref{eq:reference-parameters} provides
$\theta_0^2/\gamma+(h/2)b=a$,
$h\eta b+2(1-\eta)\theta_0(1-\theta_0)/\gamma=(1-\eta)(1-\gamma a)/\gamma$ and
\begin{align}
2ha-\frac{(1-\eta)^2(1-\gamma a)^2}{\gamma^2b}=\frac{2h}{\gamma}\theta_1^2.
\label{eq:identity-three}
\end{align}
Consider the following probability measure:
\begin{align}
\overline\nu_h(\der y,\der v)
=\mu(\der y)\,\mathcal{N}\left(aF(y),bI_d\right)(\der v).
\label{nesterov_shar_ref_measure}
\end{align}
Its position marginal is $\mu$. Let $(Y_0,V_0)\sim\bar{\nu}_h$,
independently of $\xi_1$, and define
\begin{align}
\delta_1^x
&:=h\{V_0-aF(Y_0)\}+\sqrt{\frac{2h}{\gamma}}\,\theta_0\xi_1,
\label{eq:delta-position}\\
\delta_2^v
&:=h\eta\{V_0-aF(Y_0)\}
+\sqrt{\frac{2h}{\gamma}}\,(\eta\theta_0+\theta_1)\xi_1.
\label{eq:delta-velocity}
\end{align}
The pair $(\delta_1^x,\delta_2^v)$ is centered Gaussian and
independent of $Y_0$. The corresponding covariance matrices are
\begin{align*}
\E\!\left[\delta_1^x(\delta_1^x)^\top\right]&=2haI_d,\qquad
\E\!\left[\delta_2^v(\delta_2^v)^\top\right]
=h^2bI_d,\qquad 
\E\left[\delta_1^x(\delta_2^v)^\top\right]
=\frac{h(1-\eta)(1-\gamma a)}{\gamma}I_d.
\end{align*}
Consequently, conditioning on $\delta_2^v=hv$ gives
\begin{align}
\E[\delta_1^x\mid\delta_2^v=hv]
&=\frac{(1-\eta)(1-\gamma a)}{\gamma b}\,v,
\label{eq:conditional-mean}\\
\E\!\left[\delta_1^x(\delta_1^x)^\top
\mid\delta_2^v=hv\right]
&=\frac{(1-\eta)^2(1-\gamma a)^2}{\gamma^2b^2}\,vv^\top
+\frac{2h}{\gamma}\theta_1^2I_d.
\label{eq:conditional-second-moment}
\end{align}
 The  update rules \eqref{nes_sh_eq:position-update}-(\ref{nesh_sh_eq:velocity-update}) can now be written as
\begin{align}
Y_1
&=Y_0+haF(Y_0)+\delta_1^x,
\label{eq:shifted_position}\\
h\{V_1-aF(Y_1)\}
&=\delta_2^v+\frac{h(1-\eta)(1-\gamma a)}{\gamma}F(Y_0) +\frac{h(1-\eta)}{2\gamma}\{F(Y_1)-F(Y_0)\}.
\label{nesterov_shifted-velocity-update}
\end{align}
Therefore, conditional on $Y_0=y$ and $\delta_2^v=hv$,
\begin{align*}
Y_1&=y+haF(y)+\delta_1^x,\\
V_1-aF(Y_1)
&=v+\frac{(1-\eta)(1-\gamma a)}{\gamma}F(y)
+\frac{1-\eta}{2\gamma}\{F(Y_1)-F(y)\}.
\end{align*}
 Using Tyalor's expansion, we have
\begin{align}
F(Y_1)
& = F(y)+DF(y)(Y_1-y)
+\mathcal{O}\left(
|Y_1-y|^2(1+|y|+|Y_1-y|)^r
\right).
\label{nesterov_sha_eqn_force_taylor}
\end{align}
Let $f$ be a smooth test function with derivatives having polynomial growth.
Applying Taylor's formula to $f$ and
using \eqref{eq:conditional-mean}-\eqref{eq:conditional-second-moment}, yield
\begin{align}
&\E\left[
f(Y_1,V_1-aF(Y_1))
\mid Y_0=y,\ \delta_2^v=hv
\right]
\nonumber\\
=&f
+\frac{(1-\eta)(1-\gamma a)}{\gamma}
\left\{\frac{1}{b}(v\cdot\nabla_yf)+ (F\cdot\nabla_vf)\right\}
 +ha(F\cdot\nabla_yf)
\nonumber \\  & +\frac{(1-\eta)^2(1-\gamma a)}{2\gamma^2b}\big((DF\,v)\cdot\nabla_vf\big)
\nonumber\\
& + \frac{(1-\eta)^2(1-\gamma a)^2}{2\gamma^2b^2}
\sum_{i,j=1}^d v_iv_j\partial_{y_i}\partial_{y_j}f
+ \frac{h}{\gamma}\theta_1^2\Delta_yf
+ \frac{(1-\eta)^2(1-\gamma a)^2}{\gamma^2b}
\sum_{i,j=1}^d v_iF_j\partial_{y_i}\partial_{v_j}f
\nonumber\\
&+ \frac{(1-\eta)^2(1-\gamma a)^2}{2\gamma^2}
\sum_{i,j=1}^d F_iF_j\partial_{v_i}\partial_{v_j}f
+R_hf.
\label{nestorov_sharma_taylor_expansion}
\end{align}
All functions and derivatives on the right are evaluated at
$(y,v)$ and the remainder term satisfies
\begin{align}
|R_hf|
\leq C h^3(1+|y|+|v|)^r,
\label{nesterov_remainder}
\end{align}
since the Taylor remainder contains products
of three increments of  $\mathcal{O}(h)$. The terms from
the expansion \eqref{nesterov_sha_eqn_force_taylor} which contribute $\mathcal{O}(h^3)$ are also included in reminder, since they are
multiplied by $1-\eta=O(h)$.
For ease in analysis below, we
introduce the differential operator
\begin{align}
\mathcal H_bf
=\frac{1}{b}(v\cdot\nabla_yf) + (F\cdot\nabla_vf) 
\label{nestorov_hamiltonian_operat}
\end{align}
and its square is
\begin{align}
\mathcal H_b^2f
 &=\frac1{b^2}\sum_{i,j=1}^d v_iv_j
\partial_{y_i}\partial_{y_j}f
+\frac2b\sum_{i,j=1}^d v_iF_j
\partial_{y_i}\partial_{v_j}f
\nonumber\\
& \quad +\sum_{i,j=1}^d F_iF_j
\partial_{v_i}\partial_{v_j}f
+\frac{1}{b}(F\cdot\nabla_yf)
+\frac{1}{b}(DF\,v\cdot\nabla_vf).
\label{eqnhamiltonian-square_nesterov}
\end{align}
Using \eqref{eq:identity-three}, \eqref{nestorov_hamiltonian_operat}-\eqref{eqnhamiltonian-square_nesterov}, the Taylor expansion \eqref{nestorov_sharma_taylor_expansion} becomes
\begin{align}
&\E\big[
f(Y_1,V_1-aF(Y_1))
\mid Y_0=y,\ \delta_2^v=hv\big]
\nonumber\\
& = f
+ \frac{(1-\eta)(1-\gamma a)}{\gamma}\mathcal H_bf
+\frac{(1-\eta)^2(1-\gamma a)^2}{2\gamma^2}\mathcal H_b^2f +\frac{h}{\gamma}\theta_1^2\bigl((F\cdot\nabla_yf)+\Delta_yf\bigr)
\nonumber\\
& \quad +\frac{(1-\eta)^2a(1-\gamma a)}{2\gamma b}
((DF\,v)\cdot\nabla_vf)
+R_hf(y,v).
\label{nestorov_sharma_eqn_S25}
\end{align}
We integrate \eqref{nestorov_sharma_eqn_S25} against
$\mu(\der y)\,\mathcal{N}(0,bI_d)(\der v)$. The operator $\mathcal H_b$
generates the deterministic flow $\dot y=\frac{v}{b}$ and $ \dot v=F(y)$.
This flow preserves volume and the energy Hamiltonian $U(y)+\frac{|v|^2}{2b}$. 
It therefore preserves $\mu\otimes\mathcal{N}(0,bI_d)$. We write $\mathcal N_b=\mathcal N(0,bI_d)$.
Differentiating its invariant expectation once and twice also provides
\begin{align}
\int\mathcal H_bf\der \mu\,\der\mathcal{N}_b=0,
\qquad
\int\mathcal H_b^2f\der\mu\der\mathcal{N}_b=0,
\label{eq:hamiltonian-cancellation}
\end{align}
Using integration by parts, we have
\begin{align}
\int\bigl((F\cdot\nabla_yf)+\Delta_yf\bigr)
\der \mu=0.
\label{eq:position-cancellation}
\end{align}
The first three terms after $f$  in
\eqref{nestorov_sharma_eqn_S25} are zero on integrating with respect to measure $\mu\otimes\mathcal N(0,bI_d)$, the distribution of $(Y_0,V_0-aF(Y_0))$ under $\overline\nu_h$. The remaining
displayed coefficient satisfies
\begin{align}
\frac{(1-\eta)^2a(1-\gamma a)}{2\gamma b}= \mathcal{O}(h^3).
\label{eq:remaining-coefficient}
\end{align}
Therefore, we have
\begin{align}
&\mathbb{E}_{\,\overline{\nu}_{h}}
\bigl[f(Y_1,V_1-aF(Y_1))\bigr] =
\int f(y,v)\,\mu(\der y)\,\mathcal{N}(0,bI_d)(\der v)
+\mathcal{O}(h^3).
\label{eq:shifted-stationarity}
\end{align}
Taking $f(y,v)=g(y,v+aF(y))$ we get the orignal observable function $g$. If $P_h$
denotes the numerical transition operator for Markov chain \eqref{nes_sh_eq:position-update}-\eqref{nesh_sh_eq:velocity-update}, we obtain
\begin{align}
\left|\overline\nu_h(P_hg-g)\right|
\leq C_gh^3,
\label{sharma_nestorov_stationry_approx}
\end{align}
where $C_g$ is independent of $h$ but depends on finitely many polynomially growing derivatives of $g$. In \eqref{sharma_nestorov_stationry_approx}, we have bounded
the change in expectations under $\bar{\nu}_h$ after one
step $P_h$ of the sampler.
It remains to compare $\overline{\nu}_h$ (see \eqref{nesterov_shar_ref_measure}) with $\nu_h$.
The corresponding generator of kinetic Langevin dynamics is
\begin{align}
\mathcal{L} \psi
= (v\cdot\nabla_y \psi)
+\big((F-\gamma v)\cdot\nabla_v \psi\big)
+\gamma\Delta_v \psi,
\label{eq:kinetic-generator}
\end{align}
with invariant measure
\begin{align}
\nu(\der y,\der v)=\mu(\der y)\mathcal{N}(0, I_d)(\der v).
\label{eq:kinetic-gibbs}
\end{align}
Let $S_t$ denote its semigroup. We will use weighted derivative norms denoted as
\begin{align}
\|g\|_{k,r}
=\max_{|\ell|\leq k}\sup_{y,v \in \R^d}
(1+|y|+|v|)^{-r}|D^\ell g(y,v)|,
\label{eq:weighted-norm}
\end{align}
Using Proposition~2.15 of \cite{kopec2015weak}, we get
\begin{align}
\|S_t g-\nu(g)\|_{k,r}
\leq Ce^{-\lambda t}
\|g-\nu(g)\|_{k+d+1,q},
\label{neste_sharma_eqn_deri_decay}
\end{align}
for a sufficiently large output weight $r \geq 2q$. 
Here $\nu(g)=\int g\der\nu$.
The smoothness and strong convexity assumptions ensure
the required confinement hypotheses in \cite{kopec2015weak}. Therefore,
\begin{align}
\psi=-\int_0^\infty\{S_tg - \nu(g)\}\der t
\label{eq:poisson-representation}
\end{align}
is well defined, satisfies
\begin{align}
\mathcal{L}\psi=g-\nu(g),
\label{sharma_nestorov_poissoneqn}
\end{align}
and has polynomially weighted derivative bounds controlled
by those of $g$. A second order Taylor expansion gives
\begin{align}
P_h\psi=\psi+h\mathcal{L} \psi+h\mathcal R_h \psi,
\label{eq:transition-remainder_nestorov_shar}
\end{align}
where
\begin{align}
\|\mathcal R_hg\|_{k,r'}
\leq C\sqrt h\,\|g\|_{k+3,r}.
\label{eq:transition-remainder_nestorov_shar-bound}
\end{align}
Here $r^{\prime}$ is sufficiently large compared to $r$ and is independent of $h$.
This bound uses the facts that $1-\eta=\gamma h+\mathcal{O}(h^2)$ and $|\theta_0|+|\theta_1|\leq C\sqrt{h}$. We also use $\eta\theta_0+\theta_1=(1-\eta)(1-\theta_0)$, which gives velocity noise covariance $2\gamma hI_d+\mathcal O(h^{3/2})$. Using \eqref{sharma_nestorov_poissoneqn}, we get
\begin{equation}\label{nesterov_sharma_P_h_action}
    P_h \psi= \psi + h(g - \nu(g))  + h\mathcal R_h \psi.
\end{equation}
Integrating \eqref{nesterov_sharma_P_h_action} with respect to $\nu_h$, we get
\begin{equation}
    \nu_h(g) - \nu(g) = -\nu_h(\mathcal{R}_h\psi). \end{equation}
The polynomial moment bounds justify these identities for
the possibly unbounded function $\psi$. 
Subtracting the identity obtained by integrating \eqref{nesterov_sharma_P_h_action} with respect to $\overline\nu_h$ and using \eqref{sharma_nestorov_stationry_approx}, we obtain
\begin{align}
(\nu_h-\overline\nu_h)(g)
=-(\nu_h-\overline{\nu}_h)(\mathcal{R}_h\psi)+ \mathcal{O}(h^2).
\label{eq:bootstrap-identity}
\end{align}
We apply \eqref{eq:bootstrap-identity} recursively. Set
$g_0=g$ and, for $j=0,1,2,3$, define
\begin{align}
\mathcal{L}\psi_j=g_j-\nu(g_j),
\qquad
g_{j+1}=\mathcal R_h\psi_j,
\label{eq:bootstrap-functions}
\end{align}
where $\psi_j$ is given by
\eqref{eq:poisson-representation} with $g$ replaced by $g_j$.
The semigroup derivative bounds \eqref{neste_sharma_eqn_deri_decay} and
\eqref{eq:transition-remainder_nestorov_shar-bound} yield
\begin{align}
g_j= \mathcal{O}(h^{j/2}),
\qquad
\psi_j= \mathcal{O}(h^{j/2}),
\label{eq:bootstrap-bounds}
\end{align}
for $0\leq j\leq4$ and $0\leq j\leq3$, respectively,
in the weighted derivative norms required at each stage.
The assumptions on $U$ and $g$ provide the required derivative bounds which we used here. Applying \eqref{eq:bootstrap-identity} to $g_j$ and using
\eqref{eq:bootstrap-bounds}, we obtain
\begin{align}
(\nu_h-\overline{\nu}_h)(g_j)
=-(\nu_h-\overline{\nu}_h)(g_{j+1})
+ \mathcal{O}(h^{2+j/2}).
\label{eq:bootstrap-step}
\end{align}
Four successive applications give
\begin{align}
(\nu_h-\overline{\nu}_h)(g)
=(\nu_h-\overline{\nu}_h)(g_4) + \mathcal{O}(h^2).
\label{eq:bootstrap-iteration}
\end{align}
Since $g_4 = \mathcal{O}(h^2)$ in a polynomially weighted supremum
norm, the uniform moment bounds for $\nu_h$ (see Lemma~\ref{nesterov_lem_uniform-moment-bounds_sharma}) and
$\overline{\nu}_h$ imply
\begin{align}
\left|(\nu_h-\overline{\nu}_h)(g_4)\right|
\leq C_gh^2.
\label{eq:bootstrap-final-bound}
\end{align}
Combining these estimates yields
\begin{align}
\left|\nu_h(g)-\overline{\nu}_h(g)\right|
\leq C_gh^2.
\label{eq:phase-comparison}
\end{align}
We take $g(y,v)=\varphi(y)$. Since the position
marginals of $\overline{\nu}_h$ and $\nu_h$ are $\mu$
and $\mu_h$, respectively,
\begin{align}
\overline{\nu}_h(g)=\mu(\varphi),
\qquad
\nu_h(g)=\mu_h(\varphi).
\label{eq:position-marginals}
\end{align}
Equation~\eqref{eq:phase-comparison} therefore gives
\begin{align}
\left|\mu_h(\varphi)-\mu(\varphi)\right|
\leq C_\varphi h^2.
\label{eq:conclusion}
\end{align}
Therefore, the theorem is proved.

\begin{remark}
Observables, depending on velocity, computed from $V_n$ generally would not have
$\mathcal O(h^{2})$ bias, since
$b=1-2(\gamma h)^{1/2}+\mathcal O(h)$ in
\eqref{eq:reference-parameters}. However, we can obtain second order
accuracy by using
\[
\widetilde V_n=\frac{V_n-aF(Y_n)}{\sqrt b},
\]
which requires no additional gradient evaluation.
Indeed, if $(Y,V)\sim\overline\nu_h$, then
$(Y,(V-aF(Y))/\sqrt b)\sim\nu$.
Since $a=\mathcal O(h)$ and $b\to1$ as $h \to 0$, the polynomial
derivative bounds of the transformed observable
$g(y,(v-aF(y))/\sqrt b)$ are uniform for sufficiently
small $h$. Applying \eqref{eq:phase-comparison} therefore gives
\[
\big|\E_{\nu_h}g(Y,\widetilde V)-\nu(g)\big|
\leq C_g h^2
\]
for every smooth $g$ with polynomially growing derivatives.
\end{remark}

\subsection{Proof of Theorem~\ref{nesterov_sharma_contraction_thrm}}

\begin{proof}
Let $\Delta{Y}_n = Y_n -\widetilde{Y}_{n}$, $\Delta V_n=V_n-\widetilde V_n$ and $\Delta P_n=P-\widetilde{P}_n$. Using \eqref{nesterov_sharma:eq:transformed_velocity}, we write $\Delta P_n$ as 
 \begin{align}\label{neste_p_v_relation_delta}
     \Delta{P}_n =  h\Delta V_n-\frac{h(1-\eta)(1+\eta)}{\gamma}\big(F(Y_n)-F(\widetilde{Y}_n)\big).
 \end{align}
From Taylor's expansion, we have
\begin{align}
\nabla U(Y_n)-\nabla U(\widetilde{Y}_n)
&=H_n\Delta Y_n,
&
H_n&=\int_0^1
\nabla^2U(\widetilde Y_n+r\Delta{Y}_n)\der r.
\label{nesterov_sharma_eqn_average_hessian}
\end{align}
Due to our assumption on $U$, for all $n \geq 0$, the following holds:
\begin{align}
mI_d\preceq H_n\preceq L I_d.
\end{align} Since $F=-\nabla U$, from
\eqref{nesterov_sharmaeqn_transformed_updates_y}-\eqref{nesterov_sharmaeqn_transformed_updates_p}, we have
\begin{align}
\Delta{Y}_{n+1}
&=\{I_d-\alpha(1+\eta)H_n\}\Delta{Y}_n+\Delta P_n,
\label{nesterov_sharma:eq:difference_updates_Y}\\
\Delta P_{n+1}
&=\eta\Delta P_{n}-\alpha\eta^2H_n\Delta Y_{n},
\label{nesterov_sharma:eq:difference_updates_P}
\end{align}
where, recall $\alpha = h(1- \eta)/\gamma$. Since $H$ is symmetric, choose an orthonormal basis
$e^1,\ldots,e^d$ with
\[
H_n e^i_n=\lambda^i_n e^i_n,\qquad m\leq\lambda^i_n\leq L.
\]
We work in the eigenbasis as it makes the analysis easier. We denote the coordinates of the two quantities in the norm \eqref{nesterov_sharma_norm} by
\begin{align}
\widehat{X}^i_n
&=\left( e^i_n \cdot \Delta Y_n+\frac{\Delta P_n}{1-\eta}\right),
&
\widehat{P}^i_n
&=\left( e^i_n \cdot \frac{\Delta P_n}{\eta(1-\eta)}\right).
\label{nesterov_sharma:eq:difference_coordinates}
\end{align}
Using \eqref{neste_p_v_relation_delta}, we get
\begin{align*}
\widehat{X}_n^i
&=\left( e_n^i  \cdot 
\Big( Y_n-\widetilde Y_n+\frac{h}{1-\eta}(V_n-\widetilde{V}_n)
-\frac{h(1+\eta)}{\gamma}\{F(Y_n)-F(\widetilde{Y}_n)\}
\Big)\right),\\
\widehat{P}^i_n
&=\left( e_n^i \cdot\,
\Big(\frac{h}{\eta(1-\eta)}(V_n-\widetilde V_n)
-\frac{h(1+\eta)}{\gamma\eta}\{F(Y_n)-F(\widetilde{Y}_n)\}
\Big)\right).
\end{align*}
For a fixed $n$, we express the differences at steps $n$ and
$n+1$ in the same orthonormal eigenbasis of $H_n$. Therefore, for
$j\in\{n,n+1\}$, we write
\begin{align*}
\widehat X_j^i
&=e_n^i\cdot\left(
\Delta Y_j+\frac{\Delta P_j}{1-\eta}\right), \quad 
\widehat P_j^i =e_n^i\cdot\frac{\Delta P_j}{\eta(1-\eta)}.
\end{align*}
These coordinates are used for this fixed step. Using \eqref{nesterov_sharma:eq:difference_coordinates}, we can write
\[
e_n^i\cdot\Delta Y_n
=\widehat X_n^i-\eta\widehat P_n^i.
\]
Taking components of
\eqref{nesterov_sharma:eq:difference_updates_Y}--\eqref{nesterov_sharma:eq:difference_updates_P}
and using $\alpha/(1-\eta)=h/\gamma$, we obtain
\begin{align}
\widehat X_{n+1}^i
&=\widehat X_n^i-\frac{h\lambda_n^i}{\gamma}
(\widehat X_n^i-\eta\widehat P_n^i) = \left(1-\frac{h\lambda_n^i}{\gamma}\right)\widehat X_n^i
+\frac{h\eta\lambda_n^i}{\gamma}\widehat P_n^i, 
\label{nesterov_sharma_eqn_3.20} \\
\widehat P_{n+1}^i 
&=\eta\widehat P_n^i-\frac{h\eta\lambda_n^i}{\gamma}
(\widehat X_n^i-\eta\widehat P_n^i)= -\frac{h\eta\lambda_n^i}{\gamma}\widehat X_n^i
+\left(\eta+\frac{h\eta^2\lambda_n^i}{\gamma}\right)
\widehat P_n^i. \label{nesterov_sharma_eqn_3.21}
\end{align}
Due to condition \eqref{nesterov_sharma:eq:step_condition}, we have 
\begin{align}
0<\frac{h\lambda_n^i}{\gamma}
\leq\frac{hL}{\gamma}
\leq\frac{1-\eta}{4}. \label{eqn_3.22_nesterov_sharma}
\end{align} 
We next show that the sum of the squared coordinates decreases
in each direction. Using \eqref{nesterov_sharma_eqn_3.20}-\eqref{nesterov_sharma_eqn_3.21}, we get
\begin{align}
\left(1-\frac{h\lambda_n^i}{\gamma}\right)
&\left\{(\widehat X_n^i)^2 +(\widehat P_n^i)^2\right\}
-(\widehat X_{n+1}^i)^2-(\widehat P_{n+1}^i)^2
=
\left\{\frac{h\lambda_n^i}{\gamma}
-(1+\eta^2)\frac{h^2(\lambda_n^i)^2}{\gamma^2}\right\}
(\widehat X_n^i)^2
\nonumber\\
&\qquad
+\left\{1-\eta^2-(1+2\eta^3)\frac{h\lambda_n^i}{\gamma}
-\eta^2(1+\eta^2)\frac{h^2(\lambda_n^i)^2}{\gamma^2}\right\}
(\widehat P_n^i)^2
\nonumber\\
&\qquad
-2\eta\frac{h\lambda_n^i}{\gamma}
\left\{1-\eta-(1+\eta^2)\frac{h\lambda_n^i}{\gamma}\right\}
\widehat X_n^i\widehat P_n^i.
\label{nesterov_sharma_scalar_square_3.22}
\end{align}
We now deal with each term appearing on right hand side of \eqref{nesterov_sharma_scalar_square_3.22} one by one.  For the coefficient of $(\widehat X_n^i)^2$, we use
$1+\eta^2\leq2$ and
$h\lambda_n^i/\gamma\leq1/4$ to get
\begin{align}
\frac{h\lambda_n^i}{\gamma}
-(1+\eta^2)\frac{h^2(\lambda_n^i)^2}{\gamma^2}
&=\frac{h\lambda_n^i}{\gamma}
\left\{1-(1+\eta^2)\frac{h\lambda_n^i}{\gamma}\right\}
\geq\frac{h\lambda_n^i}{2\gamma}.
\label{nesterov_sharma_first_coefficient_bound}
\end{align}
For the coefficient of $(\widehat P_n^i)^2$, we use
$\eta^3\leq\eta$,
$\eta^2(1+\eta^2)\leq2$ and
\eqref{eqn_3.22_nesterov_sharma} to get
\begin{align*}
(1+2\eta^3)\frac{h\lambda_n^i}{\gamma}
&\leq\frac{(1+2\eta)(1-\eta)}4,\quad 
\eta^2(1+\eta^2)\frac{h^2(\lambda_n^i)^2}{\gamma^2}
\leq\frac{(1-\eta)^2}{8}.
\end{align*}
Therefore, we get
\begin{align}
1-\eta^2-(1+2\eta^3)\frac{h\lambda_n^i}{\gamma}
&-\eta^2(1+\eta^2)\frac{h^2(\lambda_n^i)^2}{\gamma^2} \geq
(1-\eta)(1+\eta)
-\frac{(1+2\eta)(1-\eta)}4
-\frac{(1-\eta)^2}{8}
\nonumber\\
&\quad=
\frac58(1-\eta)(1+\eta)
\geq\frac58(1-\eta).
\label{nesterov_sharma_second_coefficient_bound}
\end{align}
It remains to obtain a bound on the third term. By
\eqref{eqn_3.22_nesterov_sharma}, $0\leq
1-\eta-(1+\eta^2)\frac{h\lambda_n^i}{\gamma}
\leq1-\eta. $ Since $\eta\leq1$, using Young's inequality, we ascertain 
\begin{align*}
2\frac{h\lambda_n^i}{\gamma}(1-\eta)
|\widehat X_n^i\widehat P_n^i|
&\leq
\frac{h\lambda_n^i}{2\gamma}(\widehat X_n^i)^2
+2\frac{h\lambda_n^i}{\gamma}(1-\eta)^2
(\widehat P_n^i)^2.
\end{align*}
Again, by \eqref{eqn_3.22_nesterov_sharma}, $
2\frac{h\lambda_n^i}{\gamma}(1-\eta)^2
\leq\frac{(1-\eta)^3}{2}
\leq\frac{1-\eta}{2}$. As a consequence, we have 
\begin{align}
2\frac{h\lambda_n^i}{\gamma}(1-\eta)
|\widehat X_n^i\widehat P_n^i|
\leq
\frac{h\lambda_n^i}{2\gamma}(\widehat X_n^i)^2
+\frac{1-\eta}{2}(\widehat P_n^i)^2.
\label{nesterov_sharma_mixed_term_bound}
\end{align}
Using
\eqref{nesterov_sharma_first_coefficient_bound},
\eqref{nesterov_sharma_second_coefficient_bound} and
\eqref{nesterov_sharma_mixed_term_bound} in
\eqref{nesterov_sharma_scalar_square_3.22}, we obtain
\begin{align*}
&\left(1-\frac{h\lambda_n^i}{\gamma}\right)
\left\{(\widehat X_n^i)^2+(\widehat P_n^i)^2\right\}
-(\widehat X_{n+1}^i)^2-(\widehat P_{n+1}^i)^2 \geq\frac{1-\eta}{8}(\widehat P_n^i)^2
\geq0.
\end{align*}
Hence, for each $i$, the following holds:
\begin{align}
(\widehat X_{n+1}^i)^2+(\widehat P_{n+1}^i)^2
&\leq
\left(1-\frac{h\lambda_n^i}{\gamma}\right)
\left\{(\widehat X_n^i)^2+(\widehat P_n^i)^2\right\}
\leq
\left(1-\frac{mh}{\gamma}\right)
\left\{(\widehat X_n^i)^2+(\widehat P_n^i)^2\right\},
\label{nesterov_sharma_scalar_contraction}
\end{align}
where the second inequality uses $\lambda_n^i\geq m$. We now sum over the eigenvector directions. Using the orthonormal basis $e_n^1,\ldots,e_n^d$, we have
\[
\sum_{i=1}^d
\left\{(\widehat X_j^i)^2+(\widehat P_j^i)^2\right\}
=\|(\Delta Y_j,\Delta P_j)\|_h^2,
\qquad j\in\{n,n+1\}.
\]
Thus \eqref{nesterov_sharma_scalar_contraction} gives
\begin{align}
\|(\Delta Y_{n+1},\Delta P_{n+1})\|_h^2
\leq
\left(1-\frac{mh}{\gamma}\right)
\|(\Delta Y_n,\Delta P_n)\|_h^2.
\label{nesterov_sharma_one_step_contraction}
\end{align}

The matrix $H_n$ and its eigenvectors may change at the next
step, but the norm \eqref{nesterov_sharma_norm} remains the
same. Since the bounds on $H_n$ hold for every $n$, we can
repeat this one-step argument using the corresponding
eigenbasis at each step. Iterating
\eqref{nesterov_sharma_one_step_contraction}, we get
\[
\|(\Delta Y_n,\Delta P_n)\|_h^2
\leq
\left(1-\frac{mh}{\gamma}\right)^n
\|(\Delta Y_0,\Delta P_0)\|_h^2.
\]
Since
$0<mh/\gamma\leq(1-\eta)/4<1$,
this proves the contraction result.
\end{proof}

\subsection{Proof of Theorem~\ref{nesterov_sharma_thrm_energy_decay}}

\begin{proof}
Recall that we can rewrite sampler \eqref{nesterov_sharma_position_update}-\eqref{nesterov_sharma_velocity_update} as \eqref{nesterov_sharmaeqn_transformed_updates_y}-\eqref{nesterov_sharmaeqn_transformed_updates_p} using $P_n=hV_n-\alpha(1+\eta)F(Y_n)$. We first express the sampler in the coordinates
\eqref{nesterov_sharma:eq:energy-coordinates}.
The transformed updates are
\begin{align}
\mathsf X_{n+1}
&=Y_n+\alpha F(Y_n)
-\sqrt{\frac{2h}{\gamma}}\,
\frac{\theta_1}{\eta}\xi_{n+1},
\nonumber\\
\mathsf Z_{n+1}
&=(1-\tau)\mathsf Z_n+\tau Y_n+\frac\tau mF(Y_n) +
\sqrt{\frac{2h}{\gamma}}
\left(
\theta_0+\frac{\eta\theta_0+\theta_1}{\eta\tau}
\right)\xi_{n+1},
\label{nesterov_sharma:eq:energy-updates}
\end{align}
where, recall, $\tau = \sqrt{m \alpha}$.  We omit the Gaussian terms and represent the resulting
updates by $\widehat{\mathsf X}_{n+1}$ and
$\widehat{\mathsf Z}_{n+1}$. The update of $\widehat{\mathsf X}_{n+1}$ is a gradient step
from $Y_n$. Since $F=-\nabla U$, using the Lipschitz continuity of $\nabla U$, we get
\begin{align}\label{eqn_3.24_neste_sharma}
U(\widehat{\sX}_{n+1})
&\leq U(Y_n)
+\alpha (\nabla U(Y_n) \cdot F(Y_n))
+\frac{L\alpha^2}{2}|F(Y_n)|^2\nonumber \\
&=U(Y_n)-\alpha\left(1-\frac{L\alpha}{2}\right)|F(Y_n)|^2 \leq U(Y_n)-\frac\alpha2|F(Y_n)|^2,
\end{align}
where the last inequality follows from $\alpha L\leq1/4$. For the squared-distance term, i.e. the last term in \eqref{nesterov_sharma_eqn_nesterov-energy}, we use the identity $
|(1-\tau)a+\tau b|^2
=(1-\tau)|a|^2+\tau|b|^2-\tau(1-\tau)|a-b|^2$,  with $a=\mathsf Z_n-x_\star$ and
$b=Y_n-x_\star$ to obtain 
\begin{align}
\frac m2|\widehat{\mathsf Z}_{n+1}-x_\star|^2
&=(1-\tau)\frac m2|\sZ_n-x_\star|^2
+\frac{m\tau}{2}|Y_n-x_\star|^2 -\frac{m\tau(1-\tau)}2|\mathsf Z_n-Y_n|^2 \nonumber \\
&\quad+\tau(1-\tau)
(\sZ_n-x_\star \cdot F(Y_n)) +\tau^2( Y_n-x_\star \cdot F(Y_n))
+\frac{\tau^2}{2m}|F(Y_n)|^2. \label{nesterov_sharma_3.24_eqn}
\end{align}
Since $\tau^2/m=\alpha$, the last term cancels the
negative squared-force term in the potential estimate. We use strong convexity of $U$ to bound the remaining terms. Splitting
$\mathsf Z_n-x_\star=(\mathsf Z_n-Y_n)+(Y_n-x_\star)$, 
the two force cross terms in \eqref{nesterov_sharma_3.24_eqn} become
$\tau(1-\tau)(\sZ_n-x_\star\cdot F(Y_n))
+\tau^2( Y_n-x_\star \cdot F(Y_n))=
\tau ( Y_n-x_\star \cdot F(Y_n))
+\tau(1-\tau)(\sZ_n-Y_n \cdot F(Y_n))$. 
Thanks to strong convexity of $U$, we have
\begin{equation*}
U(x_\star)
\geq U(Y_n)+ (\nabla U(Y_n)\cdot x_\star-Y_n)
+\frac m2|Y_n-x_\star|^2,
\end{equation*}
and, therefore $
\frac{m}{2}|Y_n-x_\star|^2
+(Y_n-x_\star\cdot F(Y_n))
\leq-\{U(Y_n)-U(x_\star)\}$.
Multiplying this inequality by $\tau$ and combining \eqref{eqn_3.24_neste_sharma}-\eqref{nesterov_sharma_3.24_eqn}, we obtain
\begin{align*}
&U(\widehat{\mathsf X}_{n+1})-U(x_\star)
+\frac m2|\widehat{\mathsf Z}_{n+1}-x_\star|^2 \leq(1-\tau)\Big[
U(Y_n)-U(x_\star)+\frac{m}{2}|\sZ_n-x_\star|^2\\
&\hspace{34mm}
+\tau (F(Y_n)\cdot\mathsf Z_n-Y_n)
-\frac{m\tau}{2}|\mathsf Z_n-Y_n|^2\Big].
\end{align*}
A second application of strong convexity at $Y_n$ gives $
U(\sX_n)
\geq U(Y_n)
+\tau (F(Y_n)\cdot \sZ_n-Y_n)
+\frac{m\tau^2}{2}|\sZ_n-Y_n|^2$ which means
\begin{align*}
&U(Y_n)+\tau ( F(Y_n) \cdot\mathsf Z_n-Y_n )
-\frac{m\tau}{2}|\mathsf{Z}_n-Y_n|^2 \leq U(\mathsf X_n)
-\frac{m\tau(1+\tau)}2|\mathsf{Z}_n-Y_n|^2
\leq U(\mathsf X_n).
\end{align*}
Therefore, the update without noise satisfies
\[
U(\widehat{\sX}_{n+1})-U(x_\star)
+\frac m2|\widehat{\sZ}_{n+1}-x_\star|^2
\leq(1-\tau)\mathcal{E}_n.
\]
We now include the Gaussian noise terms. Let us denote their coefficients
in \eqref{nesterov_sharma:eq:energy-updates} by
\begin{align*}
k_{\mathsf X}=-\sqrt{\frac{2h}{\gamma}}\frac{\theta_1}{\eta},
\qquad
k_{\mathsf Z}=\sqrt{\frac{2h}{\gamma}}
\left(\theta_0+\frac{\eta\theta_0+\theta_1}{\eta\tau}\right).
\end{align*}
Conditionally on $(Y_n,P_n)$, the quantities
$\widehat{\mathsf X}_{n+1}$ and $\widehat{\mathsf Z}_{n+1}$
are fixed, and $\xi_{n+1}$ is an independent standard
Gaussian vector. Due to Lipschitz condition on the gradient vector field, we have the following bound:
\begin{align*}
U(\sX_{n+1}) \leq U(\widehat{\sX}_{n+1})
+k_{\sX}(\nabla U(\widehat{\sX}_{n+1}) \cdot \xi_{n+1}) +\frac{L}{2}k_{\mathsf{X}}^2|\xi_{n+1}|^2.
\end{align*}
Since $\mathbb E|\xi_{n+1}|^2=d$, taking conditional expectations yields 
\begin{equation}\label{nesterov_sharma_eqn_3.26}
\mathbb E[U(\mathsf X_{n+1})\mid Y_n,P_n]
\leq U(\widehat{\mathsf X}_{n+1})+\frac{Ld}{2}k_{\mathsf X}^2.
\end{equation}
Similarly, $
\frac{m}{2}|\sZ_{n+1}-x_\star|^2
=\frac{m}{2}|\widehat{\mathsf Z}_{n+1}-x_\star|^2 +mk_{\mathsf Z}(\widehat{\mathsf Z}_{n+1}-x_\star \cdot \xi_{n+1})
+\frac m2k_{\mathsf{Z}}^2|\xi_{n+1}|^2$. 
The cross term again has conditional expectation zero, so we have
\begin{align}\label{nesterov_sharma_eqn_3.27}
\mathbb{E}\left[\frac m2|\mathsf Z_{n+1}-x_\star|^2
\,\middle|\,Y_n,P_n\right]
=\frac{m}{2}|\widehat{\sZ}_{n+1}-x_\star|^2
+\frac{md}{2}k_{\mathsf Z}^2.
\end{align}
Using \eqref{nesterov_sharma_eqn_3.26}-\eqref{nesterov_sharma_eqn_3.27}, we have
\[
\mathbb E[\mathcal E_{n+1}\mid Y_n,P_n]
\leq(1-\tau)\mathcal{E}_n+K_0,
\]
where 
$K_0 = \frac d2(Lk_{\mathsf X}^2+mk_{\mathsf Z}^2)=K_0$.
Taking expectations and applying this inequality repeatedly gives
\begin{align*}
\mathbb E\mathcal E_n
&\leq(1-\tau)^n\mathcal E_0
+K_0\sum_{j=0}^{n-1}(1-\tau)^j =(1-\tau)^n\mathcal E_0
+\frac{K_0}{\tau}\{1-(1-\tau)^n\},
\end{align*}
which proves \eqref{nesterov_sharma_eq:strongly-convex-energy_resu}.

It remains to bound the additive term. Using
$h/\gamma=\alpha/(1-\eta)$,
$\eta\theta_0+\theta_1=(1-\eta)(1-\theta_0)$,
$1-\eta=2\tau/(1+\tau)$ and $\alpha=\tau^2/m$, we obtain
\[
\frac{K_0}{\tau}
=\frac{d(1+\tau)}2
\left[
\frac{\kappa\theta_1^2}{\eta^2}
+\left(\theta_0+\frac{2(1-\theta_0)}{1-\tau}\right)^2
\right],\qquad \kappa=L/m.
\]
Here $\tau\leq1/2$, $\eta\geq1/3$,
$0<\theta_0<1$ and $\theta_1^2\leq2(1-\eta)$.
The first term is therefore bounded by $
\frac{d(1+\tau)}2\frac{\kappa\theta_1^2}{\eta^2}
\leq\frac{d(1+\tau)}2\,9\kappa\,
\frac{4\tau}{1+\tau}
=18d\kappa\tau$. 
For the second term, note that $
0<\theta_0+\frac{2(1-\theta_0)}{1-\tau}
\leq\theta_0+4(1-\theta_0)\leq 4$.
Since $(1+\tau)/2\leq3/4$, this term can be at most $12d$. Therefore, $
\frac{K_0}{\tau}\leq18d\kappa\tau+12d$, and for $\alpha=1/(4L)$, we have
$\tau=1/(2\sqrt\kappa)$ and hence
\[
\frac{K_0}{\tau}
\leq 9d\sqrt\kappa+12d
\leq 21d\sqrt{\kappa}.
\]
Moreover, $(1-\tau)^n\leq e^{-n\tau}
=e^{-n/(2\sqrt\kappa)}$, 
which proves
\eqref{nesterov_sharma:eq:accelerated-strongly-convex-energy}.
\end{proof}

\subsection{Proof of Theorem~\ref{nesterov_acc:thm:gaussian}}

\begin{proof}
Since the force is linear,
the recursion started from a fixed point has a Gaussian
distribution. Therefore, the first aim here is to obtain decaying bounds on mean and covariance. Recall that $u_i$ denotes the eigenvector with eigenvalue $a_i$. Let $\widetilde{Q}_n^i$ and $\widetilde P_n^i$ denote the components of $Y_n- x_\star$ and $V_n$ in eigen basis, i.e.  
\begin{align}\label{q_p_neste_eqn_3.17}
\widetilde Q_n^i= ( u_i \cdot Y_n-x_\star),
\qquad
\widetilde{P}_n^i= ( u_i \cdot V_n ).
\end{align}
 We multiply the position coordinate by
$\sqrt{a_i}$ so that summing the square of initial coordinates is equal to $E_0$. Substituting $F(y)=-A(y-x_\star)$ into \eqref{nesterov_sharma_position_update}-\eqref{nesterov_sharma_velocity_update} and using the coordinates \eqref{q_p_neste_eqn_3.17}, we get 
\begin{align}
\begin{pmatrix}
\sqrt{a_i}\,\widetilde Q_{n+1}^i\\
\widetilde P_{n+1}^i
\end{pmatrix}
&=
M_i \begin{pmatrix}
\sqrt{a_i}\,\widetilde Q_n^i\\
\widetilde P_n^i
\end{pmatrix}
+g_i\xi_{n+1}^i,
\end{align}
where
\begin{align}
M_i
&=
\begin{pmatrix}
1&h\sqrt{a_i}\\
-\alpha\sqrt{a_i}/h&\eta-(1+\eta)\alpha a_i
\end{pmatrix}, \;\;
g_i
=
\sqrt{\frac{2h}{\gamma}}
\begin{pmatrix}
\sqrt{a_i}\,\theta_0\\
\{[\eta-(1+\eta)\alpha a_i]\theta_0+\theta_1\}/h
\end{pmatrix},
\label{nestero_sharma_matrix_recursion_linear_A.95}
\end{align}
 $\xi_n^i=\langle u_i,\xi_n\rangle$ are independent
standard Gaussian variables, and $\theta_0$ and $\theta_1$ are from \eqref{eq:first-noise-coefficients}.
The recursion above is linear with constant coefficients and therefore, in what follows, we take the standard approach to deal with difference equations \citep{perko2013differential,elaydi2005introduction}. Its mean and covariance are determined by the powers $M_i^k$, the characteristic polynomial of $M_i$,
\begin{align*}
z^2-(1+\eta)(1-\alpha a_i)z+\eta(1-\alpha a_i).
\end{align*}
The roots of this polynomial are the eigenvalues of $M_i$. If both roots have modulus less than one, then $M_i^k$ tends to zero as $k$ grows. 
We first show that the two roots are not real. Let $z=1+w$. Then the characteristic polynomial becomes
\begin{align*}
w^2+2\sqrt{\alpha a_i}\,b_iw+\alpha a_i,
\end{align*}
where, we have denoted, $b_i= \big((1-\eta)+(1+\eta)\alpha a_i\big)/\big(2\sqrt{\alpha a_i}\big)$.
The discriminant of corresponding quadratic equation is $4\alpha a_i(b_i^2-1)$, so the roots are complex when $b_i<1$. As a function of $\sqrt{\alpha a_i}$, $b_i$ is of the form $u + 1/u$ and therefore it is convex function.  The maximum value is attained at $a_i=m$ or at $a_i=L$ over $a_i\in[m,L]$, which is either $\frac{\gamma}{2\sqrt m}+h\sqrt m$ or $h\sqrt L+\frac{\gamma}{2\sqrt L}$. Hence $b_i\leq\vartheta<1$ (see \eqref{nsacc:vartheta}) for all $i=1,\dots,d$.

Since the coefficients of the characteristic polynomial are real, its two roots are complex conjugates and have the same modulus. Let $z_i$ be the root with positive imaginary part. Using the fact that product of two roots of quadratic polynomial is the constant term appearing in the polynomial, we have
\begin{align}
|z_i|^2=\eta(1-\alpha a_i)\leq\eta=e^{-\gamma h}. \label{nesterov_sharma_eqn_3.22}
\end{align}
The modulus $|z_i|$ shows how fast $M_i^k$ decays in the long run, but it does not by itself give a bound for each $k$, because $M_i$ not only rotates vectors but also changes the magnitude. To separate these two effects, we change coordinates using a transformation matrix $\mathbb{T}_i$. To construct the transformation matrix $\mathbb{T}_i$, we use an eigenvector of $M_i$ for $z_i$ of the form $e=(1,e_2)$. A short computation with the first row of $M_i$ gives
\begin{align*}
e_2=\frac{z_i-1}{h\sqrt{a_i}}=\frac{\sqrt\alpha}h\Bigl(-b_i+  \mathrm{i}\sqrt{1-b_i^2}\Bigr),
\end{align*}
where we recall that $z_i$ is the root of the characteristic polynomial with positive imaginary part, given by
$z_i=1+\sqrt{\alpha a_i}\,(-b_i+\mathrm{i}\sqrt{1-b_i^2})$. 
The real and imaginary parts of $e = (1, e_2)$ form the columns of
\begin{align*}
\mathbb T_i=
\begin{pmatrix}
1&0\\
-\dfrac{\sqrt{\alpha}}h b_i&
\dfrac{\sqrt{\alpha}}h\sqrt{1-b_i^2}
\end{pmatrix}.
\end{align*}  
In these coordinates, we can write $\mathbb T_i^{-1}M_i\mathbb T_i=|z_i|R_i$ for an orthogonal matrix $R_i$. Therefore $M_i^k=|z_i|^k\,\mathbb T_iR_i^k\mathbb T_i^{-1}$, and since the first row of $\mathbb T_i$ is $(1,0)$ (as we are only interested in the $Y$ marginal) and $R_i^k$ is rotation and does not change length of the vector, we get, for $w\in\R^2$ and $k\geq0$,
\begin{align*}
|(1,0)M_i^kw|^2
\leq\|\mathbb T_i^{-1}\|^2|z_i|^{2k}|w|^2.
\end{align*}
A direct computation gives $\|\mathbb T_i^{-1}\|^2
\leq\frac{2h^2}{\alpha(1-b_i^2)}$. 
Since $h^2/\alpha = h\gamma / (1 -e^{-\gamma h}) \leq e^{1/2}$ for $\gamma h \leq 1$, and  since $b_i\leq\vartheta$, we obtain
\begin{align}
|(1,0)M_i^kw|^2
\leq Ce^{-\gamma hk}|w|^2,
\qquad
C=\frac{2e^{1/2}}{1-\vartheta^2}.
\label{nsacc:matrix-bound_A.98}
\end{align}
Iterating \eqref{nestero_sharma_matrix_recursion_linear_A.95} $n$ times, we ascertain
\begin{align*}
\begin{pmatrix}
\sqrt{a_i}\,\widetilde Q_n^i\\
\widetilde P_n^i
\end{pmatrix}
=M_i^n
\begin{pmatrix}
\sqrt{a_i}\,\widetilde Q_0^i\\
\widetilde P_0^i
\end{pmatrix}
+\sum_{k=0}^{n-1}M_i^kg_i\,\xi_{n-k}^i.
\end{align*}
Taking the first component, we have
\begin{align}\label{sharma_nesterov_A.99}
\sqrt{a_i}\,\widetilde Q_n^i
=(1,0)M_i^n
\begin{pmatrix}
\sqrt{a_i}\,\widetilde Q_0^i\\
\widetilde P_0^i
\end{pmatrix}
+\sum_{k=0}^{n-1}\bigl((1,0)M_i^kg_i\bigr)\,\xi_{n-k}^i.
\end{align}
The first term is deterministic, since the initial state is deterministic, and $\xi_1^i,\ldots,\xi_n^i$ are independent standard Gaussian variables. Applying
\eqref{nsacc:matrix-bound_A.98} to the initial state and summing
over $i$ gives
\begin{align}
\sum_{i=1}^d a_i
\bigl(\mathbb E\widetilde Q_n^i\bigr)^2
\leq CE_0e^{-\gamma hn}.
\label{nsacc:mean-bound}
\end{align}
We next compare the variance of $\widetilde Q_n^i$ with its
stationary variance. In the similar manner, we have (see \eqref{nestero_sharma_matrix_recursion_linear_A.95} and \eqref{sharma_nesterov_A.99})
\begin{equation}
a_i\operatorname{Var}(\widetilde Q_n^i)
= \sum_{k=0}^{n-1}(e_1^\top M_i^k g_i)^2,
\qquad e_1=(1,0)^\top.
\end{equation}
 Our bound on the powers of $M_i$ gives $
|e_1^\top M_i^k g_i|^2
\leq
\|\mathbb{T}_i^{-1}\|^2|g_i|^2|z_i|^{2k}$ (see \eqref{nsacc:matrix-bound_A.98}). 
Since $|z_i|<1$, we can tend $n \rightarrow \infty$ to get the stationary variance:
\[
a_i\sigma_{\infty,i}^2
=
\sum_{k=0}^{\infty}(e_1^\top M_i^k g_i)^2.
\]
Subtracting the finite sum gives
\begin{align*}
a_i\{\sigma_{\infty,i}^2
-\operatorname{Var}(\widetilde Q_n^i)\}
=\sum_{k=n}^{\infty}(e_1^\top M_i^k g_i)^2 \leq
\|\mathbb{T}_i^{-1}\|^2|g_i|^2
\sum_{k=n}^{\infty}|z_i|^{2k} =
\frac{\|\mathbb{T}_i^{-1}\|^2|g_i|^2}
{1-|z_i|^2}|z_i|^{2n}.
\end{align*}
The difference is nonnegative because every term in the
remaining sum is a square.

To estimate the factor multiplying $|z_i|^{2n}$, we need
to bound $|g_i|^2$ in terms of $1-\eta$ and $\alpha a_i$,
since the denominator $1-|z_i|^2$ also tends to zero as
$h$ decreases. Using
$\eta\theta_0+\theta_1=(1-\eta)(1-\theta_0)$ in
\eqref{nestero_sharma_matrix_recursion_linear_A.95}, we get
\begin{align*}
|g_i|^2
&=\frac{2ha_i}{\gamma}\theta_0^2
+\frac{2}{\gamma h}
\bigl\{(1-\eta)(1-\theta_0)
-(1+\eta)\alpha a_i\theta_0\bigr\}^2.
\end{align*}
Since $0<\theta_0<1$, $\theta_0^2\leq2(1-\eta)$ and 
$\frac{1-\eta}{\gamma h} \leq 1$,
applying $(u-v)^2\leq2u^2+2v^2$, we ascertain
\begin{align*}
|g_i|^2 &\leq
4\alpha a_i
+\frac{4(1-\eta)^2}{\gamma h}
+\frac{8(1+\eta)^2(\alpha a_i)^2(1-\eta)}
{\gamma h} \leq
4(1-\eta+\alpha a_i)+32(\alpha a_i)^2.
\end{align*}
Also, using
\eqref{nesterov_sharma_eqn_3.22}, we have $
1-|z_i|^2
=1-\eta+\eta\alpha a_i
\geq\eta(1-\eta+\alpha a_i)$. 
Combining these two estimates, we obtain
\begin{align*}
\frac{|g_i|^2}{1-|z_i|^2}
\leq \frac{1}{\eta}
\left\{4+\frac{32(\alpha a_i)^2}{1-\eta+\alpha a_i}
\right\} \leq
\frac{4+32\alpha a_i}{\eta}
\leq\frac{12}{\eta},
\end{align*}
where the last inequality follows from
$\alpha a_i\leq h^2a_i\leq1/4$.
Our assumptions also give $\gamma h\leq1$, and therefore
$\eta^{-1}=e^{\gamma h}\leq e$.
Together with the bound on $\|\mathbb T_i^{-1}\|$
obtained above, this shows that
\[
\frac{\|\mathbb{T}_i^{-1}\|^2|g_i|^2}{1-|z_i|^2} \leq C,
\]
where $C$ depends only on $\vartheta$.
Since $|z_i|^{2n}\leq\eta^n=e^{-\gamma hn}$, we conclude
\begin{align}
0 \leq
a_i\{\sigma_{\infty,i}^2
-\operatorname{Var}(\widetilde Q_n^i)\} \leq
\frac{\|\mathbb T_i^{-1}\|^2|g_i|^2}
{1-|z_i|^2}|z_i|^{2n}
\leq Ce^{-\gamma hn}.
\label{nsacc:variance-bound}
\end{align}

The estimate above gives the difference from the stationary
variance of the sampler.To compare with the target variance $1/a_i$, we also
include the nonnegative stationary error
$\epsilon_i=a_i\sigma_{\infty,i}^2-1$.
Writing $\sigma_{n,i}^2=\operatorname{Var}(\widetilde Q_n^i)$,
we obtain
\begin{align*}
a_i\sigma_{n,i}^2
&=1+\epsilon_i
-a_i(\sigma_{\infty,i}^2-\sigma_{n,i}^2)
\geq1-Ce^{-\gamma hn},\\
|a_i\sigma_{n,i}^2-1|
&\leq\epsilon_i+Ce^{-\gamma hn}.
\end{align*}

We now combine the estimates for the mean and variance.
The Gaussian increments are independent across the
eigenvector directions, so the covariance is diagonal in
this basis. The formula for relative entropy between
Gaussian measures gives
\begin{align}
\mathrm{KL}\bigl(\operatorname{Law}(Y_n)\mid\mu\bigr)
=
\frac12\sum_{i=1}^d
\left\{
a_i\bigl(\mathbb E\widetilde Q_n^i\bigr)^2
+a_i\sigma_{n,i}^2-1-\log(a_i\sigma_{n,i}^2)
\right\}.
\label{nsacc:kl-gauss}
\end{align}
For $n$ such that $Ce^{-\gamma hn}\leq1/2$, we have
$a_i\sigma_{n,i}^2\geq1/2$.
Using $t-1-\log t\leq(t-1)^2$ for $t\geq1/2$,
followed by $(a+b)^2\leq2a^2+2b^2$, we get
\begin{align*}
&\frac12\sum_{i=1}^d
\{a_i\sigma_{n,i}^2-1-\log(a_i\sigma_{n,i}^2)\}\leq
\frac12\sum_{i=1}^d(a_i\sigma_{n,i}^2-1)^2 \leq
\sum_{i=1}^d\epsilon_i^2+Cd\,e^{-2\gamma hn}.
\end{align*}
Combining the above estimate with \eqref{nsacc:mean-bound}
in \eqref{nsacc:kl-gauss} proves
\eqref{nsacc:gauss-relax}.

For the stated choice of $h$ and $\gamma$, we have
$\gamma h=1/(4\sqrt\kappa)$.
Choosing $n$ so that $
Cd\,e^{-n/(2\sqrt\kappa)}
+CE_0e^{-n/(4\sqrt\kappa)}
\leq\delta $ 
bounds the relative entropy by
$\sum_{i=1}^d\epsilon_i^2+\delta$.
This provides the mentioned gradient cost, since each iteration  requires one new gradient evaluation.
\end{proof}

\section{Details of numerical experiment with nonlinear strongly convex family}
\label{sec:supp_five_sampler_experiments}
Here, we provide additional details and results of the experiment with non-linear strongly convex potential functions. 
We first compare relaxation towards numerical equilibrium using
synchronous coupling, and then examine errors of observable
expectations relative to the Gibbs target. The coupling study
includes NS \eqref{nes_sh_eq:position-update}-(\ref{nesh_sh_eq:velocity-update}), OBABO, BAOAB \citep{leimkuhler2015molecular}, HMC \citep{duane1987hybrid} and MALA \citep{roberts1996exponential}. We also perform a separate comparison
with Bouncy Particle Sampler (BPS) based on observable errors and gradient cost.

\subsection{Potential function and the sampling methods for comparison}
\label{sec:supp_five_target}
We recall here the three-dimensional Gibbs measure considered in the paper
\begin{align}
\mu_\kappa(\der q)
&=Z_\kappa^{-1}\exp\{-U_\kappa(q)\}\der q,
\qquad q=(x,z)\in\mathbb R^2\times\mathbb R,
\label{eq:supp_five_target}\\
U_\kappa(x,z)
&=\frac{R_\kappa^2}{4}\psi_\kappa(x/R_\kappa)
+\frac{z^2}{8}+\log\cosh(z/2),
\qquad R_\kappa=\frac{10}{\sqrt\kappa}.
\label{eq:supp_five_potential_nesterov_sharma}
\end{align}

Let  $\Pi_{C_\kappa}$ denote the Euclidean projection. The gradeint is given by
\begin{align}
\nabla\psi_\kappa(u)
&=u+(\kappa-1)\Pi_{C_\kappa}(u),
\nonumber\\
\nabla U_\kappa(x,z)
&=\begin{pmatrix}
(R_\kappa/4)\{u+(\kappa-1)\Pi_{C_\kappa}(u)\}\\
z/4+(1/2)\tanh(z/2)
\end{pmatrix},
\qquad u=x/R_\kappa.
\label{eq:supp_five_gradient}
\end{align}
The potential is $m$-strongly convex with an $L$-Lipschitz gradient,
where
\begin{align}
m=\frac14,\qquad L=\frac\kappa4,\qquad
\frac{\partial^2 U_\kappa}{\partial z^2}
=\frac14+\frac14\operatorname{sech}^2(z/2)\in[1/4,1/2].
\label{eq:supp_five_curvature}
\end{align}
The scale $R_\kappa$ does not change the curvature bounds of the
polygon component. In the coordinates
\begin{align}\label{nesterov_sharma_normalized_coordi}
\mathcal{T}_\kappa(x,z)&=(x/R_\kappa,z/2)=(u,w),
\end{align}
we get
\begin{align}
U_\kappa(\mathcal T_\kappa^{-1}(u,w))
&=\frac{25}{\kappa}\psi_\kappa(u)+\frac{w^2}{2}+\log\cosh(w).
\label{eq:supp_five_normalized_potential}
\end{align}
In particular, the quadratic coefficient in the interior of the
normalized polygon is independent of $\kappa$. In experiments, projection is evaluated directly, without smoothing the potential. 

Consequently, $U_\kappa$ is continuously differentiable with a globally
Lipschitz gradient, but is not globally twice continuously
differentiable.

The coupling study compares Nesterov based sampler (NS), OBABO, BAOAB, Hamiltonian Monte
Carlo (HMC) and the Metropolis-adjusted Langevin algorithm (MALA). For the sake of convenience of the reader, we recall HMC and MALA here. Let $F : = -\nabla U$.

\textbf{Hamiltonian Monte Carlo.}  For HMC, we independently draw $p_0\sim\mathcal N(0,I_3)$ at every
proposal and set $q_0=X_n$. Let $\epsilon$ denote the leapfrog step size and let $\ell$ denote
number of leapfrog steps. Then, perform the following iterations
\begin{align}
p_{j+1/2}=p_j+(\epsilon/2)F(q_j),
\quad
q_{j+1}=q_j+\epsilon p_{j+1/2}, \quad 
p_{j+1}&=p_{j+1/2}+(\epsilon/2)F(q_{j+1}).
\label{eq:supp_five_hmc_leapfrog}
\end{align}
With $H_\kappa(q,p)=U_\kappa(q)+|p|^2/2$, accept $q_\ell$ with
probability $
a_{\mathrm{HMC}}
=1\wedge\exp\{H_\kappa(q_0,p_0)-H_\kappa(q_\ell,p_\ell)\},$ 
otherwise remain at current state $X_n$. Momentum is fully refreshed at the next
proposal. The trajectory length is fixed within each parameter
setting. A proposal costs $\ell$ new gradient evaluations, including
when it is rejected.

\textbf{Metropolis adjusted Langevin algorithm (MALA).} In MALA, a proposal computed as
\begin{align}
X^*&=X_n+\delta F(X_n)+\sqrt{2\delta}\,\xi,
\end{align}
with probability $
a_{\mathrm{MALA}}(X_n,X^*)
=1\wedge\exp\bigl\{
U_\kappa(X_n)-U_\kappa(X^*)
+\frac{|X^*-X_n-\delta F(X_n)|^2
-|X_n-X^*-\delta F(X^*)|^2}{4\delta}
\bigr\}$, 
 otherwise the sampler remains at the
current position. The proposal covariance is
$2\delta I_3$, so $\delta$ is the only scalar tuning parameter for MALA.
Each proposal requires one new gradient evaluation. HMC and MALA
preserve $\mu_\kappa$ thanks to the Metropolis correction, whereas
the three unadjusted methods (NS, OBABO and BAOAB) have $\mathcal{O}(h^2)$ systematic bias.

\subsection{Parameter selection and numerical reference measures}
\label{sec:nesterov_sharma_supp_five_tuning}
We consider the following values of $\kappa$:
\begin{align}
\kappa\in\{225,400,900,1600,3600,6400\}.
\label{eq:supp_five_kappas_nesterov_sharma}
\end{align}
\textbf{Parameter selection.} For \textbf{NS}, \textbf{OBABO} and
\textbf{BAOAB}, the common parameter grid is
\begin{align}
h\sqrt L\in\{0.4,0.6,0.8,1.0\},\qquad
\gamma/\sqrt m\in\{0.5,1,2,4\}.
\label{eq:supp_five_nesterov_sharma_kinetic_grid}
\end{align}
These dimensionless parameters express the step size and friction
relative to the curvature scales of the target.
To see this, consider the deterministic kinetic dynamics
for a quadratic potential with curvature $\lambda$:
\begin{align*}
\ddot q+\gamma\dot q+\lambda q=0,
\qquad \lambda\in[m,L].
\end{align*}
The undamped frequency is $\sqrt\lambda$, so the shortest
characteristic time scale is $1/\sqrt L$. Therefore, $h\sqrt L$
measures the numerical step relative to the fastest motion.
For a direction with curvature $m$, rescaling time by
$\tau=\sqrt m\,t$ yields
\begin{align*}
\frac{\der^2q}{\der\tau^2}
+\frac{\gamma}{\sqrt m}\frac{\der q}{\der\tau}+q=0.
\end{align*}
Hence $\gamma/\sqrt m$ measures damping relative to the slow
curvature scale. The chosen grid includes underdamped, critically
damped and overdamped regimes for this quadratic mode.
Since $m=1/4$ and $L=\kappa/4$ in
\eqref{eq:supp_five_potential_nesterov_sharma}, a fixed value of
$h\sqrt L$ gives $h\propto\kappa^{-1/2}$.
At corresponding grid points, all three methods use the same
physical step size and friction.

For \textbf{HMC}, we choose the integration step
and trajectory duration through a wider pilot search: 
\begin{align}
\epsilon\sqrt L
&\in\{0.2,0.4,0.6,0.8,1.0,1.2,1.4,1.6,1.8,2\},
\nonumber\\
T\sqrt m
&\in\{0.1,0.25,0.5,0.75,1.0,1.25,1.5,2,3,4\}.
\label{nesterov_sharma_eq:supp_five_hmc_pilot_grid}
\end{align}
We also included one-step trajectories.
For each chosen duration $T$, we rounded
$T/\epsilon$ to the nearest positive integer,
$\ell=\max\{1,\lfloor T/\epsilon+1/2\rfloor\}$, giving
$\ell$ leapfrog steps and an actual duration
$T_{\mathrm{actual}}=\ell\epsilon$.
We remove duplicate $(\epsilon,\ell)$ pairs
which leaves 106 settings to test.
For each setting, we used 512 coupled pairs per starting point with a budget of 8000 gradient
evaluations per path.
Based on this pilot study, we use the following HMC grid:
\begin{align}
\epsilon\sqrt L\in\{0.8,1.1,1.4,1.7\},\qquad
T\sqrt m\in\{0.25,0.5,0.75,1.0\}.
\label{eq:supp_five_hmc_grid_nesterov_sharma}
\end{align}
\textbf{MALA} uses
\begin{equation}
\delta L\in\{0.125,0.25,0.5,1,2,4\},
\label{eq:nesterov_sharma_supp_five_mala_grid}
\end{equation}
which includes smaller proposals and larger proposals for which
Metropolis rejection can substantially affect the dynamics.

These grids give 16 settings for each of NS,
OBABO, BAOAB and HMC, and six for MALA, hence 70 settings
at each of the six condition numbers, therefore we test 420 combinations
in total.
We compare the methods using the same target,
starting positions, coupling statistic and gradient budgets.

\textbf{Computing numerical references for each sampler.}
The separation into $x\in\mathbb R^2$ and $z\in\mathbb R$
in \eqref{eq:supp_five_normalized_potential} allows us to
generate independent Gibbs positions by rejection sampling.
These serve directly as stationary references for HMC and MALA.

For each run of NS, OBABO and BAOAB, we evolve
two independent groups of 8192 chains, one from Gibbs positions
and the other from $(4R_\kappa,R_\kappa,-6)$, using the method
and parameters being tested.
NS starts with zero displacement and independent Gaussian realizations, and OBABO and BAOAB with zero
velocity. We initially run each chain for
$B_{0,\kappa}=\lceil20000\sqrt{\kappa/400}\rceil$ steps and,
if needed, continue to $2B_{0,\kappa}$, $4B_{0,\kappa}$
and $8B_{0,\kappa}$ total steps. We compare position and velocity statistics
between the groups and between successive simulation lengths.
 These checks assess whether the results still
depend on the initial positions or change with further simulation. After two consecutive successful checks, we keep
the first 4096 final states from the group initialized from
Gibbs positions, one per coupled pair. 
All computed references for 420 settings passed the applicable checks.

\subsection{Synchronous coupling and comparison based on RMS distance}
\label{nesterov_sharma_sec:supp_five_metric}

We use three starting positions in normalized coordinates \eqref{nesterov_sharma_normalized_coordi}:
\begin{align}
\mathcal T_\kappa(q_0^{(1)})=(0,0,0),\qquad
\mathcal T_\kappa(q_0^{(2)})=(1,0,3),\qquad
\mathcal T_\kappa(q_0^{(3)})=(4,0,-3).
\label{nesterov_sharma_initializations_for_samplers}
\end{align}
NS starts with $Y_{-1}=Y_0$, OBABO and BAOAB with zero
velocity, and HMC with standard Gaussian momentum.
For each parameter setting, we perform two independent runs
of $N=4096$ coupled pairs per starting position, giving
8192 pairs per inital point. We divide each run into 32 independent groups of 128 coupled
pairs and estimate uncertainty by bootstrap resampling of
these groups with replacement. For each pair index, the three starts
share a reference trajectory and future random variables.
We keep this dependence when resampling the groups. In particular, each chain and its reference use the same Gaussian vectors,
including both OBABO vectors and the HMC momentum. For NS, the initial stored noise is also shared, while the reference position history is kept. For HMC, the two copies share the freshly drawn momentum
and the uniform variable used for acceptance, while each
uses its own acceptance probability
\citep[Algorithm~1]{heng2019unbiased}. The two MALA chains also use the same uniform random variable to decide acceptance, each with its own acceptance probability.

Let $Q_{n,i}^{(r,s)}$ and $\overline Q_{n,i}^{(r)}$ denote
the chain and reference positions for run $r$ ($r = 1,2$), start $s$ ($s = 1,2,3$)
and pair $i$ ($i=1,\dots,4096$). Let $P_1(q) := (q_1,q_2)$ and $P_2(q) := q_3$. Define
\begin{align}
d_{s,j}(n)
&=\left\{\frac{1}{2N}\sum_{r=1}^2\sum_{i=1}^N
\left|P_j\left[\mathcal T_\kappa(Q_{n,i}^{(r,s)})
-\mathcal T_\kappa(\overline Q_{n,i}^{(r)})\right]\right|^2
\right\}^{1/2},
\nonumber\\
C_j(n)&=\max_{1\leq s\leq3}\frac{d_{s,j}(n)}{d_{s,j}(0)},
\qquad C(n)=\max\{C_1(n),C_2(n)\}.
\label{eq:supp_five_coupling_nesterov_sharma_metric}
\end{align}
Here, $d_{s,j}(0)$ is the initial RMS distance for start $s$ and component $j$. Thus, $100C(n)\leq1$ means that each of these RMS distances is at most $1\%$ of its initial value.

\begin{figure}[htbp!]
\centering
\includegraphics[width=\textwidth,height=0.2\textheight]{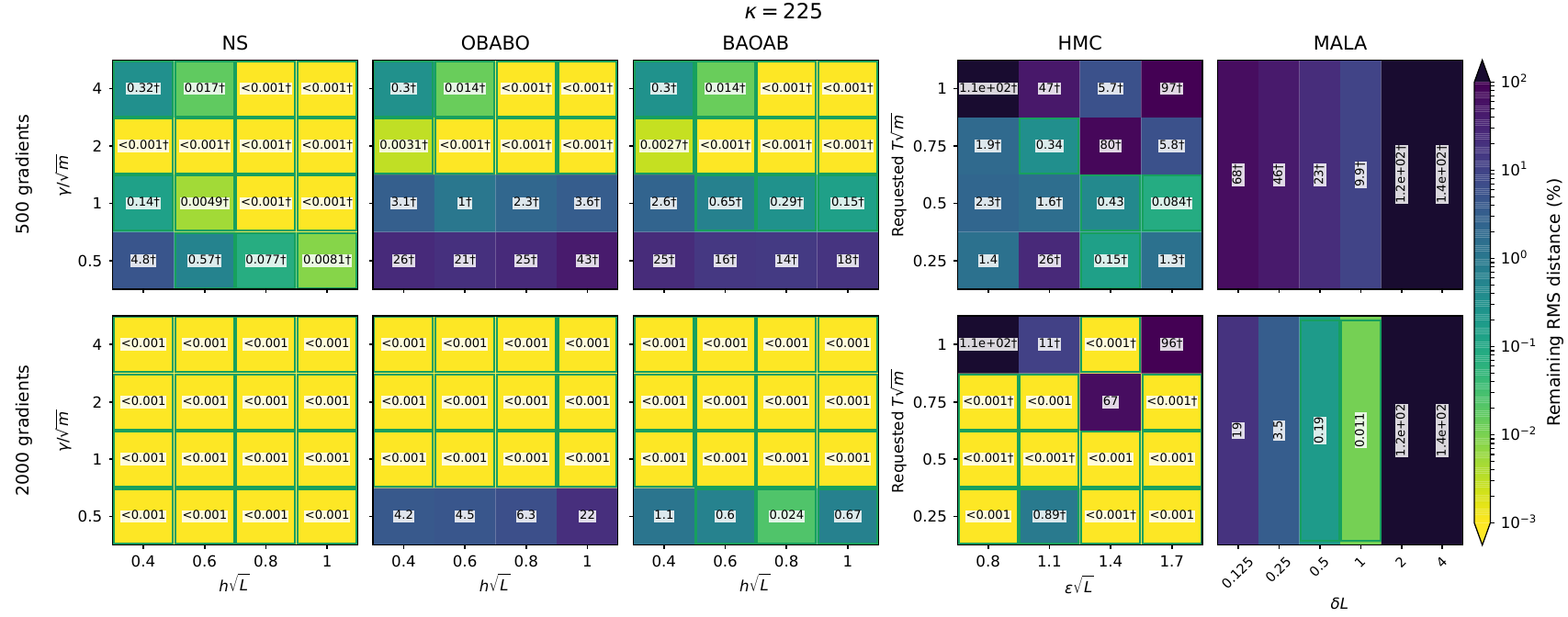}
\caption{Normalized synchronous RMS distance \eqref{eq:supp_five_coupling_nesterov_sharma_metric}, in percent,
at $\kappa=225$. Columns show NS, OBABO, BAOAB, HMC and
MALA and rows show budgets of 500 and 2000 gradients per path.
The grids are from \eqref{eq:supp_five_nesterov_sharma_kinetic_grid},
\eqref{eq:supp_five_hmc_grid_nesterov_sharma} and
\eqref{eq:nesterov_sharma_supp_five_mala_grid}.
Green outlines mark distances at or below $1\%$ and daggers
mark observations below the stated budget. For HMC,
`requested' denotes the duration before rounding the
number of integration steps.}
\label{fig:supp_five_heatmap225_nesterov_sharma}
\end{figure}

\begin{figure}[htbp!]
\centering
\includegraphics[width=\textwidth,height=0.15\textheight]{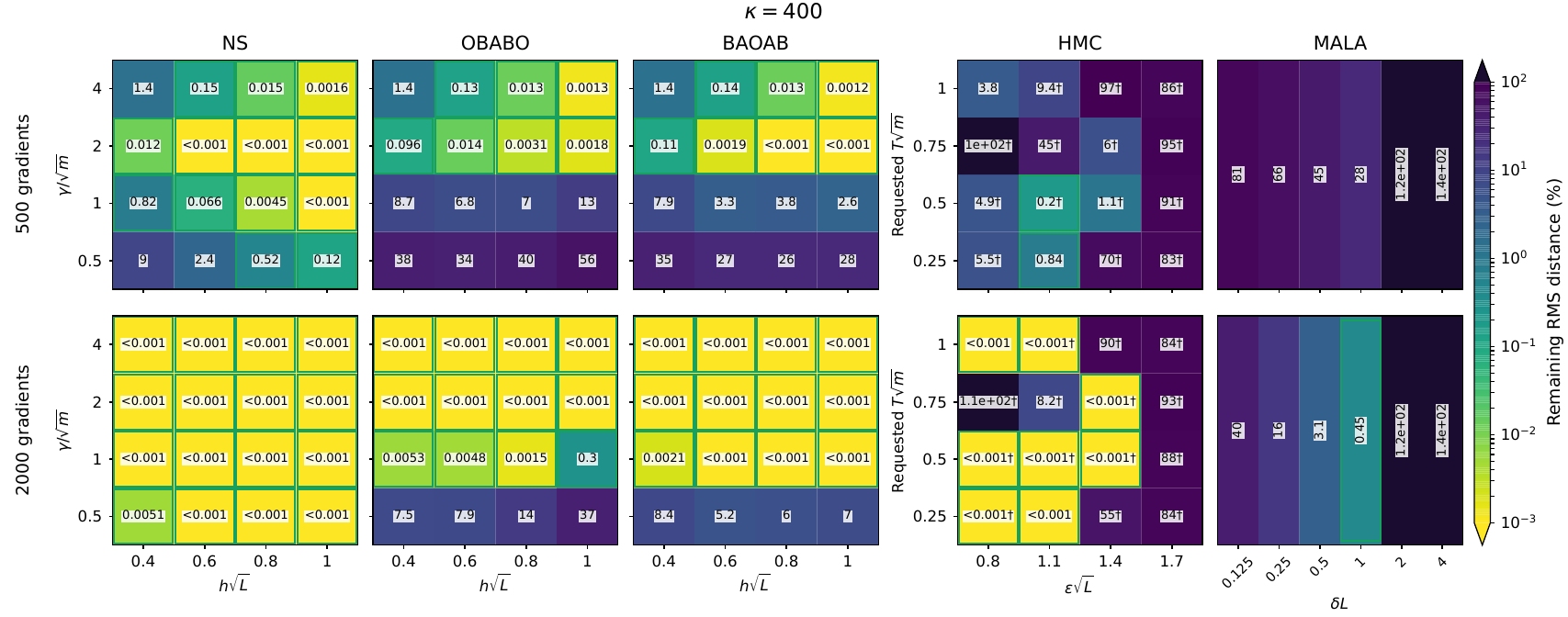}
\caption{Comparison at $\kappa=400$, with the same
conventions as Figure~\ref{fig:supp_five_heatmap225_nesterov_sharma}.}
\label{fig:nesterov_sharma_supp_five_heatmap400}
\end{figure}

\begin{figure}[htbp!]
\centering
\includegraphics[width=\textwidth,height=0.15\textheight]{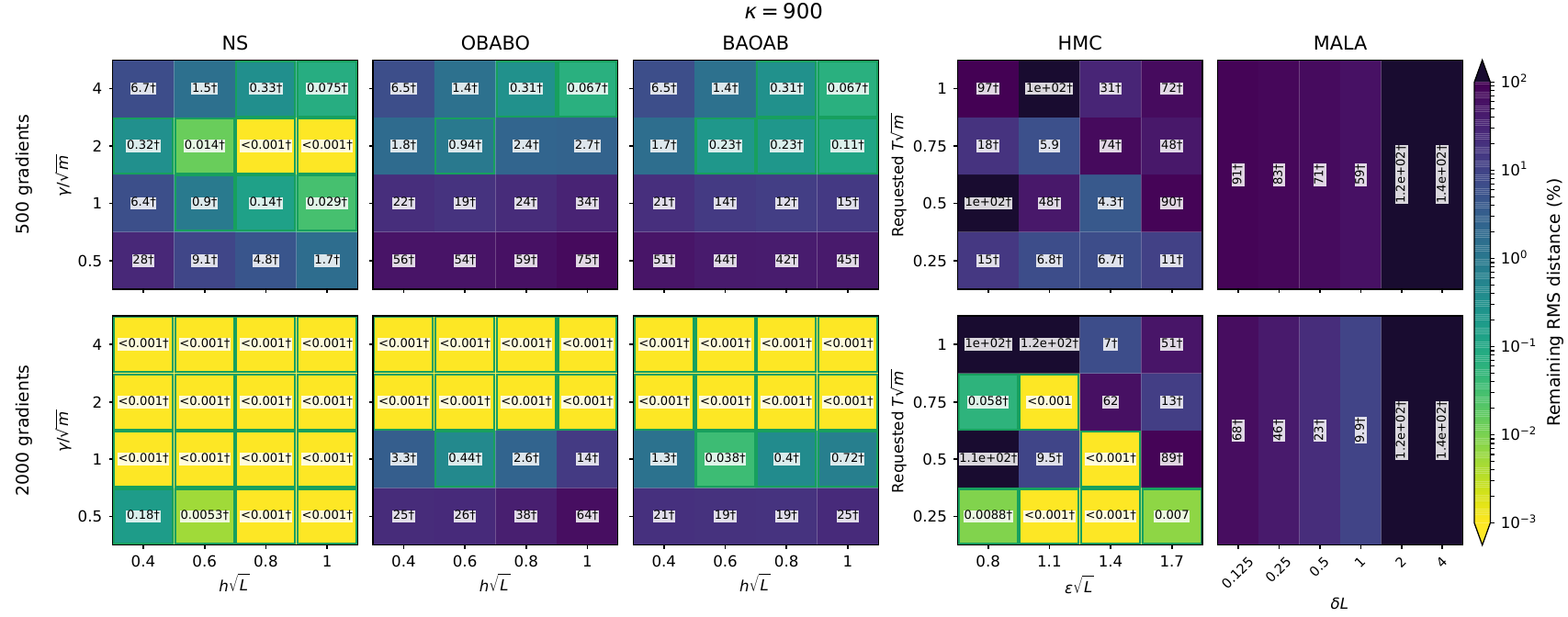}
\caption{Comparison at $\kappa=900$, with the same
conventions as Figure~\ref{fig:supp_five_heatmap225_nesterov_sharma}.}
\label{nesterov_sharma_fig:supp_five_heatmap900}
\end{figure}

\begin{figure}[htbp!]
\centering
\includegraphics[width=\textwidth,height=0.15\textheight]{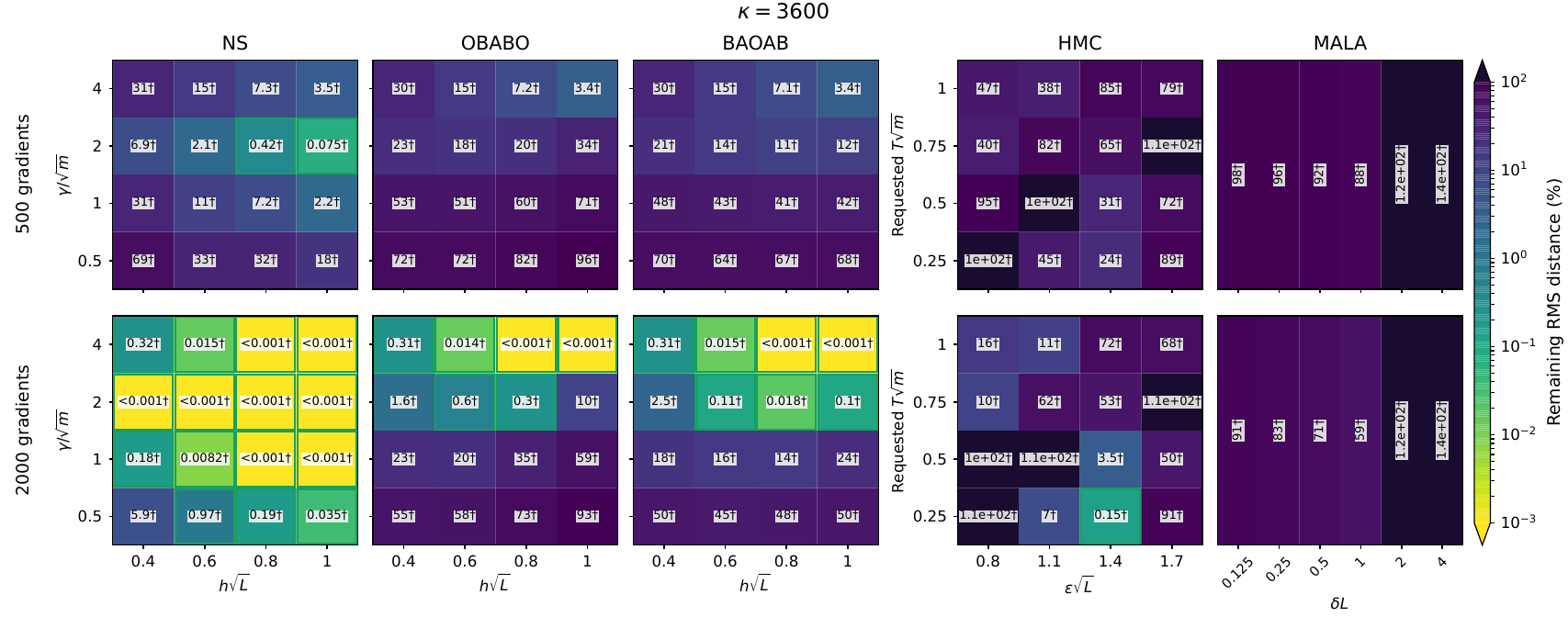}
\caption{Comparison at $\kappa=3600$, with the same
conventions as Figure~\ref{fig:supp_five_heatmap225_nesterov_sharma}.}
\label{fig:supp_five_heatmap3600_nesterov_sharma}
\end{figure}

Figures~\ref{fig:supp_five_heatmap225_nesterov_sharma}-\ref{fig:supp_five_heatmap3600_nesterov_sharma}
show $100C(g)$ at the latest recorded gradient cost $g$
not exceeding 500 or 2000 per path. 
A distance recorded below $1\%$ at one budget does not mean that it will remain below
it at later observations.

\begin{table}[htbp!]
\centering
\caption{Gradient cost $\tau_{0.01}$ of a sustained
reduction to $1\%$ of the initial normalized RMS distance \eqref{eq:supp_five_coupling_nesterov_sharma_metric}.
The parameter choices are: $h\sqrt L=0.8$ and $\gamma/\sqrt m=2$ for NS,
OBABO and BAOAB, $\epsilon\sqrt L=1.1$ and $T\sqrt m=0.5$
for HMC, and $\delta L=1$ for MALA.
Unresolved means that the required reduction and persistence
were not observed within the available budget.}
\label{tab:supp_five_costs_nesterov_sharma}
\begin{tabular}{r r r r r r}
\hline
$\kappa$ & NS & OBABO & BAOAB & HMC & MALA\\
\hline
225  & 104 & 184  & 120  & 546        & 968\\
400  & 140 & 250  & 200  & 360        & 1720\\
900  & 225 & 795  & 465  & 3976       & 3870\\
1600 & 300 & 920  & 680  & 936        & 6880\\
3600 & 450 & 2760 & 1440 & Unresolved & 15480\\
6400 & 600 & 6440 & 3400 & 3672       & 27520\\
\hline
\end{tabular}
\end{table}

We define $\tau_{0.01}$ as the first recorded gradient cost
from which $C(g)$ remains at or below $0.01$ at every recorded
observation over an interval of at least
$\lceil400\sqrt{\kappa/400}\rceil$ gradient evaluations per path.
The cost $\tau_{0.01}$ marks the beginning of this interval and 
the additional cost needed to check whether the closeness persists is not included.
If no such interval is observed, we report the threshold as
unresolved.
Under these rules, NS requires fewer gradients than OBABO
and BAOAB at every tested $\kappa$
(Table~\ref{tab:supp_five_costs_nesterov_sharma}).
To examine the dependence on $\kappa$, we fit $
\log\tau_{0.01}(\kappa)=a+p\log\kappa+\varepsilon_\kappa
$ by least squares, using all six condition numbers and then
the four values $\kappa\geq900$. We report a fit only when
all thresholds in that range are resolved, and an interval
only when all bootstrap fits are complete (see Table~\ref{nesterov_sharma_tab:supp_five_exponents}).

\begin{table}[htbp!]
\centering
\caption{Fitted exponents for the parameters in
Table~\ref{tab:supp_five_costs_nesterov_sharma}, with $95\%$
bootstrap intervals in parentheses.
The final column gives $R^2$ for the fit over all six condition numbers.}
\label{nesterov_sharma_tab:supp_five_exponents}
\begin{tabular}{l c c c}
\hline
Method & All condition numbers & $\kappa\geq900$ & $R^2$\\
\hline
NS & $0.525\;(0.523,0.528)$ & $0.500\;(0.500,0.542)$ & $0.999$\\
OBABO & $1.055\;(1.007,1.150)$ & $1.108\;(1.022,1.342)$ & $0.979$\\
BAOAB & $0.962\;(0.904,0.985)$ & $1.001\;(0.899,1.167)$ & $0.992$\\
HMC & Not fitted & Not fitted & ---\\
MALA & $1.000\;(0.974,1.001)$ & $1.000\;(0.997,1.002)$ & $1.000$\\
\hline
\end{tabular}
\end{table}

For the parameters in
Table~\ref{tab:supp_five_costs_nesterov_sharma}, the NS cost
grows approximately as $\sqrt\kappa$, whereas the OBABO,
BAOAB and MALA costs grow approximately as $\kappa$.
Across the grid tested for NS, all 96 parameter-condition-number
combinations resolve the threshold \eqref{tab:supp_five_nesterov_sharma_resolved}, and the 16 fitted exponents
range from $0.488$ to $0.543$ (see Table~\ref{tab:supp_five_nesterov_sharma_resolved}). 

Figure~\ref{neste_shar_supp_five_nesterov_sharma_scaling_all}
shows that relaxation costs $\tau_{0.01}$ vary with the parameters. As can be seen, NS reaches the $1\%$ threshold at all 16 settings for every
condition number.
The OBABO and BAOAB panels also show approximately
square-root growth for some parameter choices. However, this behaviour
under one fixed choice does not hold across the grids tested for OBABO and BAOAB while NS
shows square root growth behavior over all $16$ tested settings. For MALA, the four settings $\delta L\leq1$
resolve the threshold at every $\kappa$, with exponents between
$0.999$ and $1.002$, except the two larger settings
(see Table~\ref{tab:supp_five_nesterov_sharma_resolved}).
The only HMC parameter choice which is resolved at all six condition numbers
is $\epsilon\sqrt L=1.1$, $T\sqrt m=0.25$. The corresponding plot is nonmonotone, and the fitted exponent
$0.267$ has $R^2=0.211$ and therefore does not provide enough support to infer
a single power law.

\begin{table}[htbp!]
\centering
\caption{Number of settings with a resolved sustained
$1\%$ threshold. At each $\kappa$, we test 16 settings
for NS, OBABO, BAOAB and HMC, and six for MALA.}
\label{tab:supp_five_nesterov_sharma_resolved}
\begin{tabular}{l r r r r r r}
\hline
Method & 225 & 400 & 900 & 1600 & 3600 & 6400\\
\hline
NS & 16 & 16 & 16 & 16 & 16 & 16\\
OBABO & 15 & 15 & 14 & 14 & 12 & 12\\
BAOAB & 16 & 16 & 16 & 16 & 15 & 15\\
HMC & 13 & 9 & 10 & 10 & 10 & 11\\
MALA & 4 & 4 & 4 & 4 & 4 & 4\\
\hline
\end{tabular}
\end{table}

\begin{table}[htbp!]
\centering
\caption{Lowest recorded relaxation costs $\tau_{0.01}$ among
the 16 settings in Figure~\ref{neste_shar_supp_five_nesterov_sharma_scaling_all} for NS and BAOAB. The costs are for $h\sqrt L=1$ for both methods.
For NS, $\gamma/\sqrt m=2$ and for BAOAB, it is
$2$ at $\kappa=225,400$ and $4$ at $\kappa\geq900$.}
\label{tab:nesterov_sharma_lowest_coupling_costs}
\begin{tabular}{r r r}
\hline
$\kappa$ & NS & BAOAB\\
\hline
225  & 88  & 96\\
400  & 110 & 180\\
900  & 180 & 330\\
1600 & 240 & 440\\
3600 & 360 & 660\\
6400 & 480 & 880\\
\hline
\end{tabular}
\end{table}

\subsection{Computing ergodic limit}
\label{sec:nesterov_sharma_supp_observable_accuracy}
We compute ergodic limit for an observable with respect to the invariant measure \eqref{eq:supp_five_target} parameterized by $\kappa$ and, in this subsection,
we compare instantaneous ensemble estimates for NS, OBABO, BAOAB,
HMC, MALA and the bouncy particle sampler (BPS). We take the observable
\begin{align}
\varphi_{\mathrm{s}}(x,z)&=(z/2)^2,
\end{align}
and the corresponding ergodic limit is 
\begin{align}
\mu_\kappa(\varphi_{\mathrm{s}})
&=\frac{\int_{\mathbb R}w^2e^{-w^2/2}/\cosh(w)\,\der w}
{\int_{\mathbb R}e^{-w^2/2}/\cosh(w)\,\der w}
=0.591834 \quad (6 \text{ d.p.}).
\label{eq:supp_scalar_reference_nesterov_sharma}
\end{align}
This expectation is independent of $\kappa$ and is computed
by quadrature with absolute and relative tolerances $10^{-12}$.
For each method and condition number in 
\eqref{eq:supp_five_kappas_nesterov_sharma}, we use 8192
independent trajectories. For the first five methods (i.e. excluding BPS), we select parameters from the preceding
coupling experiment (see Table~\ref{tab:supp_nesterov_sharma_observable_parameters}).
For NS, OBABO and BAOAB, we first choose the
friction giving the lowest resolved $1\%$ coupling cost for each
method, then select the step minimizing the largest of these
three costs.  For HMC and MALA, we choose the settings with the
lowest resolved coupling costs.
The observable estimates are not used in this selection. The selected parameters are given in
Table~\ref{tab:supp_nesterov_sharma_observable_parameters}.

\begin{table}[htbp!]
\centering
\caption{Selected parameters for the first five methods, NS, OBABO, BAOAB, HMC and MALA.
The step $h$ is common to NS, OBABO and BAOAB.
For HMC, $\epsilon$ is the leapfrog step, $T$ is
the duration before rounded off to $\ell$ leapfrog steps, and
the actual duration is $\ell\epsilon$.
MALA uses $\delta L=1$.}
\label{tab:supp_nesterov_sharma_observable_parameters}
\begin{tabular}{r c c c c c c r}
\hline
& & NS & OBABO & BAOAB & \multicolumn{3}{c}{HMC}\\
\cline{6-8}
$\kappa$ & $h\sqrt L$
& $\gamma/\sqrt m$ & $\gamma/\sqrt m$ & $\gamma/\sqrt m$
& $\epsilon\sqrt L$ & $T\sqrt m$ & $\ell$\\
\hline
225  & 1.0 & 2 & 2 & 2 & 1.4 & 0.25 & 3\\
400  & 1.0 & 2 & 4 & 2 & 1.1 & 0.50 & 9\\
900  & 1.0 & 2 & 4 & 4 & 1.4 & 0.50 & 11\\
1600 & 1.0 & 2 & 4 & 4 & 0.8 & 0.50 & 25\\
3600 & 0.8 & 2 & 4 & 4 & 1.4 & 0.25 & 11\\
6400 & 0.8 & 2 & 4 & 4 & 1.1 & 0.25 & 18\\
\hline
\end{tabular}
\end{table}

For BPS, the position moves with constant velocity between
events. Bounces occur at rate
$[v^\top\nabla U_\kappa(q)]_+$ and change the velocity according to
\begin{align}
v\longmapsto
v-2\frac{v^\top\nabla U_\kappa(q)}
{|\nabla U_\kappa(q)|^2}\nabla U_\kappa(q).
\label{eq:supp_bps_nesterov_sharma_reflection}
\end{align}
Independent refreshment events replace the velocity by a
$\mathcal N(0,I_3)$ vector. We simulate bounce times by windowed
thinning using the Lipschitz gradient bound, with refreshment
rate and window length $\lambda_{\mathrm{BPS}}=\sqrt{3m}=\frac{\sqrt3}{2}$ and $w_{\mathrm{BPS}}=(3L)^{-1/2}
=\frac{2}{\sqrt{3\kappa}}.$
All methods start from $(R_\kappa,0,6)$.
NS starts with zero displacement and an independent stored
Gaussian realization, OBABO and BAOAB with zero velocity,
and BPS with velocity distributed as $\mathcal N(0,I_3)$.
HMC draws a fresh Gaussian momentum at each proposal. Let $Q_{n,i}$ denote the position of trajectory $i$ at the
$n$th recorded observation. In the results presented in the paper, we compute
\begin{align}
\widehat m_n(\varphi_{\mathrm{s}})
&=\frac1N\sum_{i=1}^N\varphi_{\mathrm{s}}(Q_{n,i}),
\qquad N=8192,
\nonumber\\
E_n(\varphi_{\mathrm{s}})
&=\left|\widehat m_n(\varphi_{\mathrm{s}})
-\mu_\kappa(\varphi_{\mathrm{s}})\right|.
\label{eq:supp_instantaneous_nesterov_sharma_gibbs_error}
\end{align}

\section{Details of numerical experiment with a smoothed
total-variation potential}
\label{nesterov_sharma_tv_comparison}

We use the parameter grid
\eqref{eq:supp_five_nesterov_sharma_kinetic_grid}, auxiliary
initialization, synchronous coupling and persistence criterion
described in
Appendices~\ref{sec:nesterov_sharma_supp_five_tuning}
and~\ref{nesterov_sharma_sec:supp_five_metric}. 
We use $\kappa\in\{10,100,1000,10000\}$.
Each setting uses two independent runs of 256 coupled
pairs per start, with observations every 20 new gradient
evaluations up to 40000. Costs exclude reference preparation
and the initial gradient.
The three starting positions are
\[
y_\star,\qquad y_\star+4\sqrt d\,u_1,\qquad
y_\star+4\sqrt d\,u_2,
\]
where $y_\star$ is the numerical minimizer,
$u_1=d^{-1/2}(1,\ldots,1)^\top$, and $u_2$ is a fixed
normalized Gaussian direction.
The coupling distance uses the full position vector,
without coordinate rescaling or separation into components.

For reference preparation, the two groups contain 2048
chains each, initialized from $\mathcal N(y_\star,I_d)$
and $\mathcal N(y_\star+4\sqrt d\,u_1,I_d)$.
The initial length is
$B_{0,\kappa}=\max\{1000,\lceil1000\sqrt{\kappa/100}\rceil\}$,
with successive doubling up to $8B_{0,\kappa}$.
We compare ten fixed projections and their squares,
hyperbolic tangents of position projections, and squared
position and velocity norms.
Normalized differences, including three estimated standard
errors, must not exceed $0.20$. Empirical distribution
differences for position projections must not exceed $0.08$. After two consecutive successful checks, we keep 256 states from the first group. Both runs must qualify and failed references are reported separately from unresolved thresholds.

\section{Details of Bayesian logistic regression for mortgage applications}
\label{sec:supp_hmda}

\textbf{Data and data processing.} We consider $179465$ home-purchase applications in North Carolina
from the 2025 HMDA data \citep{ffiec2026hmda}. We set the response to one for denial
($24073$ applications) and zero for loan issued  or approval
without acceptance. We split the dataset into $80\%$ training
and $20\%$ test observations. Within each set, we keep all
denials and randomly select an equal number of non-denials
without replacement. This gives $38516$ training and $9630$
test observations, with a $50\%$ baseline accuracy. The presented results are for this balanced data. For loan amount $a$, income $r$ and loan term $t$, measured in
dollars, thousands of dollars and months, respectively, we use transformation $
z_a=\log(a/100000)$, $z_r=\operatorname{arsinh}(r/100)$ and
$z_t=\log(t/360)$.

Nonpositive amounts and terms are treated as missing, and zero
and negative incomes are are not changed. For data processing, we use spline terms, scaling and replacement of missing
values by training medians by following standard regression and data processing methods
\citep{harrell2001regression,kuhn2013applied}.

\textbf{Model.} We model $\operatorname{logit}(p_i)=x_i^\top\beta+b_{\ell(i)}$,
where $\ell(i)$ corresponds to the lender. We keep $1003$ lender
effects with $32$ features, giving $1035$ parameters.  We use independent centred Gaussian
priors with variances $25$ for the intercept, $6.25$ for the
other coefficients and $1$ for lender effects. Writing
$q=(\beta^\top,b^\top)^\top$, $a_i^\top=(x_i^\top,e_{\ell(i)}^\top)$ and $D=\operatorname{diag}(0.04,0.16I_{31},I_{1003})$, where $e_\ell$ is the $\ell$th coordinate vector, the potential is
\begin{align}
U(q)=\sum_{i=1}^{38516}
\{\log(1+e^{a_i^\top q})-y_i a_i^\top q\}
+\tfrac12q^\top Dq.
\label{eq:supp_hmda_potential}
\end{align}
Every gradient uses all training observations. 

\textbf{Numerical reference.} We use eight NUTS chains with
the original slice-sampling version
\citep[Section~3.1]{hoffman2014no}, identity mass matrix and
maximum tree depth $12$. During the first $800$ of $1000$
warmup iterations, we adapt one step size for all eight
chains, targeting an acceptance statistic of $0.9$ using smoothening formula in
\citep[Algorithm~6]{hoffman2014no}. The resulting
step size, $0.0073805322$, is then kept fixed.
We discard all warmup samples and simulate $1024$ further
samples per chain. We use these $8192$ samples to estimate
coefficient means, standard deviations and predictive
probabilities. The maximum $\widehat R$ is $1.00544$, and
the minimum bulk and tail effective sample sizes are
$3223.6$ and $3411.2$, respectively. The total cost is
$8307480$ gradient evaluations, including $4121352$
during warmup.

\textbf{Parameter selection.} We select parameters based on  numerical reference computed using NUTS. Each search
trajectory starts with intercept
$\operatorname{logit}(19258/143572)$ and all other coordinates
zero.
Table~\ref{tab:supp_hmda_parameter_grids} provides the search grids. The initial search uses $16$, $40$, $160$ and $20$ trajectories
per setting for each kinetic method (i.e., NS, OBABO and BAOAB), MALA, NUTS and BPS,
respectively, with a target budget of $2000$ gradient
evaluations per trajectory. We divide squared coefficient
errors by the reference variances and average over coefficients,
trajectories and the recorded budgets $500$, $1000$ and $2000$.
We use this average to rank settings. A setting is eligible
if the reference checks pass, no divergences are recorded,
and the mean cost is within $25\%$ of the target budget.
The two best performing settings for each method are then
assessed using $64$ new trajectories at budgets
$2000$, $4000$ and $8000$.

\begin{table}[t]
\centering
\caption{Parameter grids for the initial and further searches.
For NUTS, $\epsilon_{\mathrm{pilot}}=0.0053959105$ is the
step size obtained during pilot adaptation with target
acceptance statistic $0.9$
\citep[Algorithm~6]{hoffman2014no}.
For BPS, $\lambda_{\mathrm{BPS}}$ is the refreshment rate and $w$ is the
window length.}
\label{tab:supp_hmda_parameter_grids}
\small
\begin{tabular}{p{0.17\linewidth} p{0.35\linewidth} p{0.35\linewidth}}
\hline
Method & Initial search & Further search\\
\hline
NS, OBABO, BAOAB
& $h\sqrt L\in\{0.2,0.4,0.8,1.2\}$,\newline
  $\gamma\in\{0.25,0.5,1,2,4\}$
& ---\\
MALA
& $\delta L\in\{0.125,0.25,0.5,1,$\newline $2,4,8,16\}$
& $\delta L\in\{2,3,4,5,6,7\}$\\
HMC
& ---
& $\epsilon\sqrt L\in\{0.4,0.8,1.2,1.6\}$,\newline
  $\ell\in\{8,32,128,512\}$\\
NUTS
& Target acceptance statistic $0.8$ or $0.9$
& $\epsilon/\epsilon_{\mathrm{pilot}}\in\{0.25,0.5,0.75,1\}$\\
BPS
& $\lambda_{\mathrm{BPS}}\in\{0.25,1,4,16\}$,\newline
  $w\sqrt L\in\{0.5,1,2,4\}$
& $\lambda_{\mathrm{BPS}}\in\{8,16,32,64\}$,\newline
  $w\sqrt L\in\{2,4,8,16\}$\\
\hline
\end{tabular}
\end{table}

The further search uses $32$ trajectories per setting at
nominal budgets $2000$ and $4000$. We complete all MALA,
HMC and NUTS settings and $14$ of the $16$ BPS settings,
without further independent runs to confirm the choices.  For MALA, HMC and BPS,
we select the completed setting with the smallest coefficient
error at the final recorded budget.

For the final comparison, NS, OBABO and BAOAB use the same
$h=0.003073$ and $\gamma=1$. 
MALA uses $\delta=3.2795495\times10^{-5}$; HMC uses
$\epsilon=h$ and $32$ leapfrog steps; NUTS uses
$\epsilon=0.0040469329$ and maximum tree depth $12$;
BPS uses $\lambda_{\mathrm{BPS}}=32$ and $w=0.0409771808$.

\textbf{Initialization and implementation.} We sample $64$ independent initial positions from the prior
and use the same positions for all methods.
NS starts with zero displacement and an independent stored
Gaussian vector, and OBABO and BAOAB start with zero velocity. The initial momentum for
HMC and NUTS is standard Gaussian, and BPS starts
with standard Gaussian velocity.

For BPS, we use windowed thinning \citep{lu2026windowed} to simulate bounce times.
A matrix upper bound on the Hessian, computed from the
training data, gives an upper bound on the bounce rate within
each time window. We update this bound at candidate events,
velocity refreshments and window boundaries to account for
the current position and velocity. The gradient at an accepted
bounce is reused to avoid evaluating it again.

We aim to compare the methods at budgets
$100,500,1000,2000,4000,8000$ gradient evaluations per trajectory.
Since the cost of a NUTS iteration or a BPS time interval is
random, we use separate preliminary runs to choose NUTS
iteration counts and BPS observation times that approximately
match these budgets. These are fixed before the comparison,
because stopping each trajectory at a fixed gradient count
can change the distribution of the recorded states.
We therefore plot the actual costs averaged over trajectories.

HMC uses complete proposals within each budget so that every
recorded state follows a full transition. For NUTS, the chosen
$1,3,4,6,7,12$ iterations give costs
$232,340,515,1097$, $1864,6979$. 

\textbf{Estimates and comparison.}
We compare how accurately the samplers estimate coefficients
and predict test outcomes within the given budgets.
For each trajectory $i$, we average the coefficient values
and predictive probabilities over all iterations up to the
recorded budget. For BPS, we average over physical time
because its position changes continuously.
We keep the initial part of each trajectory before the chains approach the posterior distribution.

Let $\overline\beta_{i,j}(G)$ be the average of coefficient $j$
at budget $G$, and let $\widehat\mu_j$ and $\widehat s_j$
be its reference mean and standard deviation.
Since the coefficients have different scales, we divide
each difference from the reference mean by its reference
standard deviation and report
\begin{align}
E(G)=\left\{
\frac1{64\times32}\sum_{i=1}^{64}\sum_{j=0}^{31}
\left(\frac{\overline\beta_{i,j}(G)-\widehat\mu_j}
{\widehat s_j}\right)^2
\right\}^{1/2}.
\label{eq:supp_hmda_coefficient_error}
\end{align}
This error includes the $32$ regression coefficients.
Lender effects are also sampled and used in predictions,
but are not included in this error.

To assess predictions, we compute test log loss and accuracy
from each trajectory's averaged predictive probabilities,
using $1/2$ as the threshold for predicting denial.
We then average these scores over the $64$ trajectories
and compare them with the reference log loss and accuracy,
$0.46491184$ and $77.66355\%$. We estimate uncertainty in our estimation by drawing $64$ whole trajectories
with replacement from the simulated trajectories and
repeating this $10000$ times. Each selected trajectory keeps
all its recorded values. We select matching trajectories
across methods so that results from the same initial position
remain together. For each resample, we recalculate the errors
and prediction scores. The $2.5$th and $97.5$th percentiles
give $95\%$ intervals at each recorded budget (see Table~\ref{table_hmda_final_results_nestorov_sharma}).

\begin{table}[t]
\centering
\caption{Results at the target budget of $8000$ gradient
evaluations per trajectory. Actual costs are averaged over
trajectories. Brackets give $95\%$ bootstrap intervals for
the coefficient error.}
\label{table_hmda_final_results_nestorov_sharma}
\small
\begin{tabular}{l r c r r}
\hline
Method & Actual cost & Coefficient error & Test log loss & Accuracy (\%)\\
\hline
NS & 8000 & $0.330\ [0.297,0.368]$ & 0.46493 & 77.664\\
OBABO & 8000 & $1.356\ [1.171,1.543]$ & 0.46464 & 77.690\\
BAOAB & 8000 & $1.357\ [1.173,1.547]$ & 0.46464 & 77.686\\
MALA & 8000 & $2.787\ [2.627,2.951]$ & 0.47883 & 76.741\\
HMC & 7969 & $17.206\ [12.754,21.297]$ & 0.80627 & 74.765\\
NUTS & 6979 & $20.708\ [18.110,23.254]$ & 0.66730 & 68.290\\
BPS & 8069.38 & $20.753\ [19.001,22.599]$ & 0.66176 & 66.566\\
\hline
\end{tabular}
\end{table}

Table~\ref{table_hmda_final_results_nestorov_sharma} shows that NS gives
smaller coefficient error than OBABO and BAOAB, while all
three attain similar prediction scores. Computations use TensorFlow 2.10.1, double precision and an NVIDIA A100 GPU.

\end{document}